\documentclass[12pt]{article}
\usepackage{amsmath}
\usepackage{amssymb}
\renewcommand\smallskip{\vskip\smallskipamount}
\renewcommand\medskip{\vskip\medskipamount}
\renewcommand\bigskip{\vskip\bigskipamount}

\newcommand{\qed}{\hfill $\Box$ \medskip}

\renewcommand{\thefootnote}{}

\begin{document}

\footnotetext{\hspace{-.6cm}\textsc{M.A. Khuri}\vspace{.01cm}}

\footnotetext{\hspace{-.6cm}Department of Mathematics, Stony Brook
University, Stony Brook, NY 11794\vspace{.01cm}}

\footnotetext{\hspace{-.6cm}(e-mail:
khuri@math.sunysb.edu)\vspace{.01cm}}
\renewcommand{\thefootnote}{\fnsymbol{footnote}}
\begin{large}
\noindent\textbf{Local Solvability of a Class of Degenerate
Monge-Amp\`{e}re Equations and Applications to
Geometry\footnote[1]{Research partially supported by an NSF
Postdoctoral Fellowship.} }\end{large}

\bigskip\smallskip

\noindent\textbf{Marcus A. Khuri}

\bigskip

\small

\noindent\textbf{Abstract.}  We consider two natural problems
arising in geometry which are equivalent to the local solvability
of specific equations of Monge-Amp\`{e}re type.  These are: the
problem of locally prescribed Gaussian curvature for surfaces in
$\mathbb{R}^{3}$, and the local isometric embedding problem for
two-dimensional Riemannian manifolds. We prove a general local
existence result for a large class of Monge-Amp\`{e}re equations
in the plane, and obtain as corollaries the existence of regular
solutions to both problems, in the case that the Gaussian
curvature possesses a nondegenerate critical point.\bigskip

\noindent\textit{Mathematics Subject Classification} (2000): 53B,
53A05, 35M10

\bigskip

\noindent\hrulefill
\normalsize
\bigskip\bigskip

\noindent\textbf{1. Introduction}\setcounter{equation}{0}
\setcounter{section}{1}
\bigskip

\noindent Let $K(u,v)$ be a function defined in a neighborhood of
a point in $\mathbb{R}^{2}$, say \linebreak $(u,v)=0$.  A
well-known problem is to ask, when does there exist a piece of a
surface \linebreak $z=z(u,v)$ in $\mathbb{R}^{3}$ having Gaussian
curvature $K$?\par
   The classical results on this problem may be found in [10],
[19], and [20].  They show that a solution always exists when $K$
is analytic or $K$ does not vanish at the origin.  In the case
that $K\geq 0$ and is sufficiently smooth, or $K(0)=0$ and $\nabla
K(0)\neq 0$, C.-S. Lin provides an affirmative answer in [15] and
[16] (see [4] for a simplified proof of [16]).  When $K\leq 0$ and
$\nabla K$ possesses a certain nondegeneracy, Han, Hong, and Lin
[8] show that a solution always exists.  Furthermore, if $K$
degenerates to arbitrary finite order on a single smooth curve,
then Q. Han and the author independently provide an affirmative
answer in [5] and [11] (see also [6] for improved regularity).
For an excellent survey of these results and related topics, see
[7]. In this paper we prove the following,\medskip

\noindent\textbf{Theorem 1.1.}  \textit{Suppose that the origin is
a nondegenerate critical point for $K$ and $K\in C^{l}$, $l\geq
100$.  Then there exists a piece of a $C^{l-98}$ surface in
$\mathbb{R}^{3}$ with Gaussian curvature $K$.}\medskip

   If a surface in $\mathbb{R}^{3}$ is given by $z=z(u,v)$, then
its Gaussian curvature is given by
\begin{equation}
 z_{uu}z_{vv}-z^{2}_{uv}=K(1+|\nabla z|^{2})^{2}.
\end{equation}
Therefore our problem is equivalent to the local solvability of
the above equation.\par
   Another well-known and related problem, is that of the local
isometric embedding of surfaces into $\mathbb{R}^{3}$.  That is,
if $(M^{2},ds^{2})$ is a two-dimensional Riemannian manifold, when
can one realize this, locally, as a small piece of a surface in
$\mathbb{R}^{3}$?  Suppose that $ds^{2}=Edu^{2}+2Fdudv+Gdv^{2}$ is
given in the neighborhood of a point, say $(u,v)=0$.  Then we must
find three function $x(u,v)$, $y(u,v)$, $z(u,v)$, such that
$ds^{2}=dx^{2}+dy^{2}+dz^{2}$.  The following strategy was first
used by J. Weingarten [25].  We search for a function $z(u,v)$,
with $|\nabla z|$ sufficiently small, such that $ds^{2}-dz^{2}$ is
flat in a neighborhood of the origin.  Suppose that such a
function exists, then since any Riemannian manifold of zero
curvature is locally isometric to Euclidean space (via the
exponential map), there exists a smooth change of coordinates
$x(u,v)$, $y(u,v)$ such that $dx^{2}+dy^{2}=ds^{2}-dz^{2}$.
Therefore, our problem is reduced to finding $z(u,v)$ such that
$ds^{2}-dz^{2}$ is flat in a neighborhood of the origin.  A
computation shows that this is equivalent to the local solvability
of the following equation,
\begin{equation}
(z_{11}-\Gamma_{11}^{i}z_{i})(z_{22}-\Gamma_{22}^{i}z_{i})-
(z_{12}-\Gamma_{12}^{i}z_{i})^{2}=K(EG-F^{2}-Ez_{2}^{2}-
Gz_{1}^{2}+2Fz_{1}z_{2}),
\end{equation}
where $z_{1}=\partial z/\partial u$, $z_{2}=\partial z/\partial
v$, $z_{ij}$ are second derivatives of $z$, and $\Gamma_{jk}^{i}$
are Christoffel symbols.  For this problem we obtain a similar
result to that of theorem 1.1.\medskip

\noindent\textbf{Theorem 1.2.}  \textit{Suppose that the origin is
a nondegenerate critical point for $K$ and $ds^{2}\in C^{l}$,
$l\geq 102$.  Then there exists a $C^{l-100}$ local isometric
embedding into $\mathbb{R}^{3}$.}\medskip

  We note that A. V. Pogorelov has constructed a $C^{2,1}$ metric
with no $C^{2}$ isometric embedding in $\mathbb{R}^{3}$.  Other
examples of metrics with low regularity not admitting a local
isometric embedding have also been proposed by Nadirashvili and
Yuan [17].  Furthermore, an alternate method for obtaining
\textit{smooth} examples of local nonsolvability, for equations
with similar structure, may be found in [12].\par
   Equations (1.1) and (1.2) are both two-dimensional Monge-Amp\`{e}re
equations.  With the goal of treating both problems
simultaneously, we will study the local solvability of the
following general Monge-Amp\`{e}re equation
\begin{equation}
\det(z_{ij}+a_{ij}(u,v,z,\nabla z))=Kf(u,v,z,\nabla z),
\end{equation}
where $a_{ij}(u,v,p,q)$ and $f(u,v,p,q)$ are smooth functions of
$p$ and $q$, $f>0$, and $a_{ij}(0,0,p,q)=
\partial^{\alpha}a_{ij}(0,0,0,0)=0$, for any multi-index $\alpha$ in
the variables $(u,v)$ satisfying $|\alpha|\leq 2$.  Clearly (1.1)
is of the form (1.3), and (1.2) is of the form (1.3) if
$\Gamma_{jk}^{i}(0)=0$, which we assume without loss of
generality.  We will prove\medskip

\noindent\textbf{Theorem 1.3.}  \textit{Suppose that $K(0)=|\nabla
K(0)|=0$,} $\det\text{Hess}\hspace{0.01in}K(0)\neq 0$ \textit{or}
$\text{Hess}\hspace{0.01in}K(0)$ \textit{has at least one negative
eigenvalue, and $K$, $a_{ij}$, $f\in C^{l}$, $l\geq 100$. Then
there exists a $C^{l-98}$ local solution of} (1.3).\medskip

\noindent\textbf{Remark.}  1) \textit{The methods carried out
below may be slightly modified to yield the same result for the
case when} $\text{Hess}\hspace{0.01in}K(0)$ \textit{has at least
one positive eigenvalue; and therefore ultimately include the case
of genuine second order vanishing, that is, when $K(0)=|\nabla
K(0)|=0$ and $|\nabla^{2}K(0)|\neq 0$.  It is conjectured that
local solutions exist whenever $K$ vanishes to finite order and
the $a_{ij}$ vanish to an order greater than that of $K$.}\par

2) \textit{Recently Q. Han} [9] \textit{together with the author,
have shown that local solutions exist for the isometric embedding
problem whenever $K$ vanishes to finite order and the zero set
$K^{-1}(0)$ consists of Lipschitz curves intersecting transversely
at the origin. Unfortunately the methods of} [9] \textit{breakdown
when the transversality assumption is removed. Therefore theorem}
1.3 (\textit{which allows tangential intersections}) \textit{and
the methods used to prove it, may be considered as a first step
towards the general conjecture.}\medskip

   Equation (1.3) is elliptic if $K>0$, hyperbolic if $K<0$, and
of mixed type if $K$ changes sign in a neighborhood of the origin.
Furthermore, the order to which $K$ vanishes determines how (1.3)
changes type in the following way.  If $K(0)=0$ and $\nabla
K(0)\neq 0$ [16], then (1.3) is a nonlinear perturbation of the
Tricomi equation:
\begin{equation*}
vz_{uu}+z_{vv}=0.
\end{equation*}
In our case, assuming that the origin is a nondegenerate critical
point for $K$, (1.3) is a nonlinear perturbation of Gallerstedt's
equation [3]:
\begin{equation*}
\pm v^{2}z_{uu}+z_{vv}=0.
\end{equation*}
Therefore, if sufficiently small linear perturbation terms are
added to the above two equations, then the first (second) partial
$v$-derivative of the $z_{uu}$ coefficient will not vanish for the
Tricomi (Gallerstedt) equation.  It is this fact, which allows one
to obtain appropriate estimates for the linearized equation of
(1.3) in both cases.  This observation, lemma 2.3 below, is the
key to our approach.\par
   From now on we only consider the case when
$\text{Hess}\hspace{0.01in}K(0)$ has at least one negative
eigenvalue, since the case of two positive eigenvalues may be
treated by the results in [15] when $K$ is nonnegative. Therefore,
we can assume without loss of generality that
\begin{equation*}
Kf(u,v,z,\nabla z)=-v^{2}+O(|u|^{2}+|v|^{3}+|z|^{2}+|\nabla
z|^{2}).
\end{equation*}
Let $\varepsilon$ be a small parameter and set
$u=\varepsilon^{4}x$, $v=\varepsilon^{2}y$,
$z=u^{2}/2-v^{4}/12+\varepsilon^{9}w$.  Then substituting into
(1.3) and cancelling $\varepsilon^{5}$ on both sides, equation
(1.3) becomes
\begin{equation}
-y^{2}w_{xx}+w_{yy}+\varepsilon\widetilde{F}(\varepsilon,x,y,w,
\nabla w,\nabla^{2} w)=0,
\end{equation}
where $\widetilde{F}(\varepsilon,x,y,p,q,r)$ is smooth with
respect to $\varepsilon$, $p$, $q$, and $r$.  Choose $x_{0}$,
$y_{0}>0$ and define the rectangle $X=\{(x,y) \mid |x|< x_{0},
|y|< y_{0}\}$.  Let $\psi \in C^{\infty}(X)$ be a cut-off function
such that
$$
\psi(x,y)=\begin{cases}
1 & \text{if $|x|\leq \frac{x_{0}}{2}$ and $|y|\leq \frac{y_{0}}{2}$},\\
0 & \text{if $|x|\geq \frac{3x_{0}}{4}$ or $|y|\geq
\frac{3y_{0}}{4}$},
\end{cases}
$$
and cut-off the nonlinear term by $F(\varepsilon,x,y,w,\nabla
w,\nabla^{2} w)=\psi\widetilde{F}$.  Then solving
\begin{equation}
\Phi(w)=-y^{2}w_{xx}+w_{yy}+\varepsilon F(\varepsilon,x,y,w,
\nabla w,\nabla^{2} w)=0\text{ }\text{ }\text{ }\text{ in }\text{
}\text{ }X,
\end{equation}
is equivalent to solving (1.3) locally at the origin.\par
   In the following sections, we shall study the linearization of
(1.5) about some function $w$.  The linearized equation is a small
perturbation of Gallerstedt's equation, which as mentioned above
admits certain estimates.  These estimates are sufficient for the
existence of weak solutions, however the perturbation terms cause
some difficulty in proving higher regularity.  To avoid this
problem, we will regularize the equation by appending a suitably
small fourth order operator.  In section $\S 2$ we shall prove the
existence of weak solutions for a boundary value problem
associated to this modified linearized equation. Regularity will
be obtained in section $\S 3$.  In section $\S 4$ we make the
appropriate estimates in preparation for the Nash-Moser iteration
procedure. Finally, in $\S 5$ we apply a modified version of the
Nash-Moser procedure and obtain a solution of (1.5).

\bigskip\bigskip
\noindent\textbf{2. Linear Existence Theory}
\setcounter{equation}{0} \setcounter{section}{2}
\bigskip\bigskip

\noindent In this section we will prove the existence of weak
solutions for a small perturbation of the linearized equation for
(1.5). Fix a constant $\Lambda>0$, and for all $i,j=1,2$ let
$b_{ij}$, $b_{i}$, $b\in C^{r}(\mathbb{R}^{2})$ be such
that:\bigskip

$i)$ the supports of $b_{ij}$, $b_{i}$, and $b$ are contained in
$X$, and\bigskip

$ii)$ $\sum |b_{ij}|_{C^{10}}+|b_{i}|_{C^{10}}+|b|_{C^{10}}\leq
\Lambda$.\bigskip

\noindent We will study the following generalization of the
linearization for (1.5),
\begin{equation}
L=\sum_{i,j}a_{ij}\partial_{x_{i}x_{j}}+\sum_{i}a_{i}\partial_{x_{i}}
+a,
\end{equation}
where $x_{1}=x$, $x_{2}=y$ and $a_{11}=-y^{2}+\varepsilon b_{11}$,
$a_{12}=\varepsilon b_{12}$, $a_{22}=1+\varepsilon b_{22}$,
$a_{1}=\varepsilon b_{1}$, $a_{2}=\varepsilon b_{2}$,
$a=\varepsilon b$.\par
   To simplify (2.1), we shall make a change of variables that will
eliminate the mixed second derivative term.  In constructing this
change of variables we will make use of the following lemma from
ordinary differential equations.\medskip

\noindent\textbf{Lemma 2.1 [1].}  \textit{Let $G(x,t)$ be a smooth
real valued function in the closed rectangle $|x-s|\leq T_{1}$,
$|t|\leq T_{2}$.  Let $M=\sup|G(x,t)|$ in this domain. Then the
initial-value problem $dx/dt=G(x,t)$, $x(0)=s$, has a unique
smooth solution defined on the interval $|t|\leq
\min(T_{2},T_{1}/M)$.}\medskip

   We now construct the desired change of variables.\medskip

\noindent\textbf{Lemma 2.2.}  \textit{For $\varepsilon$
sufficiently small, there exists a $C^{r}$ diffeomorphism
\begin{equation*}
\xi=\xi(x,y), \eta=y
\end{equation*}
of $X$ onto itself, such that in the new variables $(\xi,\eta)$
\begin{equation*}
L=\sum_{i,j}\overline{a}_{ij}\partial_{x_{i}x_{j}}+
\sum_{i}\overline{a}_{i}\partial_{x_{i}} +\overline{a},
\end{equation*}
where $x_{1}=\xi$, $x_{2}=\eta$,
$\overline{a}_{11}=-\eta^{2}+\varepsilon \overline{b}_{11}$,
$\overline{a}_{12}\equiv0$, $\overline{a}_{22}=1+\varepsilon
\overline{b}_{22}$, $\overline{a}_{1}=\varepsilon
\overline{b}_{1}$, $\overline{a}_{2}=\varepsilon
\overline{b}_{2}$, $\overline{a}=\varepsilon\overline{b}$, and
$\overline{b}_{ij}$, $\overline{b}_{i}$, $\overline{b}$
satisfy:}\bigskip

\textit{ $i)$ $\overline{b}_{ij},\overline{b}_{i},\overline{b} \in
C^{r-2}(\overline{X})$,}\bigskip

\textit{ $ii)$ $\overline{b}_{ij}$, $\overline{b}_{i}$, and
$\overline{b}$ vanish in a neighborhood of the lines $\xi=\pm
x_{0}$, and}\bigskip

\textit{ $iii)$ $\sum
|\overline{b}_{ij}|_{C^{8}(\overline{X})}+|\overline{b}_{i}|_{C^{8}(\overline{X})}
+|\overline{b}|_{C^{8}(\overline{X})}\leq \Lambda'$,}\bigskip

\noindent \textit{for some fixed $\Lambda'$.}\medskip

\noindent\textit{Proof.} Using the chain rule we find that
$\overline{a}_{12}=a_{12}\xi_{x}+a_{22}\xi_{y}$. Therefore, we
seek a smooth function $\xi(x,y)$ such that
\begin{equation}
a_{12}\xi_{x}+a_{22}\xi_{y}=0\text{ }\text{ in $X$, }\text{ }
\xi(x,0)=x, \text{ }\text{ } \xi(\pm x_{0},y)=\pm x_{0}.
\end{equation}
The boundary condition $\xi(\pm x_{0},y)=\pm x_{0}$ states that
the vertical sides of $\partial X$ will be mapped identically onto
themselves under the transformation $(\xi,\eta)$.  Moreover, the
horizontal portion of $\partial X$ will be mapped identically onto
itself since $\eta=y$.  Thus, $(\xi,\eta)$ will act as the
identity map on $\partial X$.\par
   Since $a_{12}=\varepsilon b_{12}$ and $a_{22}=1+\varepsilon
b_{22}$, by property $(ii)$ if $\varepsilon$ is sufficiently small
the line $y=0$ will be non-characteristic for (2.2).  Then by the
theory of first order partial differential equations, (2.2) is
reduced to the following system of first order ODE:
\begin{eqnarray*}
\dot{x}&=&\frac{a_{12}}{a_{22}},\text{ }\text{ }\text{
}x(0)=s,\text{ }-x_{0} \leq
s\leq x_{0},\\
\dot{y}&=&1,\text{ }\text{ }\text{ }\text{ }\text{ }\text{ }y(0)=0,\\
\dot{\xi}&=&0,\text{ }\text{ }\text{ }\text{ }\text{ }\text{
}\xi(0)=s,\text{ }\text{ }\text{ }\text{ }\xi(\pm x_{0},y)=\pm
x_{0},
\end{eqnarray*}
where $x=x(t)$, $y=y(t)$, $\xi(t)=\xi(x(t),y(t))$ and $\dot{x}$,
$\dot{y}$, $\dot{\xi}$ are derivatives with respect to $t$.\par
   We first show that the characteristic curves, given
parametrically by $(x,y)=(x(t),t)$, exist globally for $-y_{0}\leq
t\leq y_{0}$.  We apply lemma 2.1 with $T_{1}=2x_{0}$ and
$T_{2}=y_{0}$ to the initial-value problem
$\dot{x}=\frac{a_{12}}{a_{22}}$, $x(0)=s$.  By property $(ii)$ for
the $b_{ij}$
\begin{equation*}
M\leq \sup_{X}|\frac{a_{12}}{a_{22}}|=\varepsilon
\sup_{X}|\frac{b_{12}}{1+\varepsilon b_{22}}|\leq \varepsilon
C_{0},
\end{equation*}
so for $\varepsilon$ small, $M\leq \frac{2x_{0}}{y_{0}}$.  Thus
$\min(T_{2},T_{1}/M)=y_{0}$, and lemma 2.1 gives the desired
global existence.\par
   We observe that $\xi=s$ is constant along each characteristic.
In particular, since $\frac{a_{12}}{a_{22}}|_{(\pm x_{0},y)}=0$
the characteristics passing through $(\pm x_{0},0)$ are the
vertical lines $(\pm x_{0},t)$, so that $\xi(\pm x_{0},y)=\pm
x_{0}$ is satisfied.\par
   We now show that the map $\rho:X\rightarrow X$ given by
\begin{equation*}
(s,t)\mapsto(x(s,t),y(s,t))=(x(s,t),t)
\end{equation*}
is a diffeomorphism, from which we will conclude that $\xi=s(x,y)$
is a smooth function of $(x,y)$.  To show that $\rho$ is 1-1,
suppose that $\rho(s_{1},t_{1})=\rho(s_{2},t_{2})$. Then
$t_{1}=t_{2}$ and $x(s_{1},t_{1})=x(s_{2},t_{2})$, which implies
that $s_{1}=s_{2}$ by uniqueness for the initial-value problem for
ordinary differential equations.  To show that $\rho$ is onto,
take an arbitrary point $(x_{1},y_{1})\in X$, then we will show
that there exists $s\in [-x_{0},x_{0}]$ such that
$\rho(s,y_{1})=(x(s,y_{1}),y_{1})=(x_{1},y_{1})$.  Since the map
$x(s,\cdot):[-x_{0},x_{0}]\rightarrow [-x_{0},x_{0}]$ is
continuous and $x(\pm x_{0},\cdot)=\pm x_{0}$, the intermediate
value theorem guarantees that there is $s\in [-x_{0},x_{0}]$ with
$x(s,y_{1})=x_{1}$, showing that $\rho$ is onto.  Therefore,
$\rho$ has a well-defined inverse.\par
   To show that $\rho^{-1}$ is smooth it is sufficient, by the
inverse function theorem, to show that the Jacobian of $\rho$ does
not vanish at each point of $X$.  Since
\begin{equation*}
D\rho=\left(%
\begin{array}{cc}
  x_{s} & x_{t} \\
  0 & 1 \\
\end{array}%
\right),
\end{equation*}
this is equivalent to showing that $x_{s}$ does not vanish in $X$.
Differentiate the equation for $x$ with respect to $s$ to obtain,
$\frac{d}{dt}(x_{s})=(\frac{a_{12}}{a_{22}})_{x}x_{s}$,
$x_{s}(0)=1$.  Then by the mean value theorem
\begin{equation*}
|x_{s}(s,t)-1|=|x_{s}(s,t)-x_{s}(s,0)|\leq y_{0}
\sup_{X}|(\frac{a_{12}}{a_{22}})_{x}|\sup_{X}|x_{s}|
\end{equation*}
for all $(s,t)\in X$.  Thus, by property $(ii)$ for the $b_{ij}$,
\begin{equation*}
1-\varepsilon C_{1} y_{0} \sup_{X}|x_{s}|\leq x_{s}(s,t)\leq
\varepsilon C_{1} y_{0} \sup_{X}|x_{s}|+1
\end{equation*}
for all $(s,t)\in X$.  Hence for $\varepsilon$ sufficiently small,
$x_{s}(s,t)>0$ in $X$.  We have now shown that $\rho$ is a
diffeomorphism.  Moreover, by lemma 2.1 and the inverse function
theorem we have $\rho,\rho^{-1}\in C^{r}$.\par
   Lastly, we calculate $\overline{a}_{11}$, $\overline{a}_{22}$,
$\overline{a}_{1}$, $\overline{a}_{2}$, and show that they possess
the desired properties.  It will first be necessary to estimate
the derivatives of $\xi$.  By differentiating (2.2) with respect
to $x$, we obtain
\begin{equation*}
(\frac{a_{12}}{a_{22}})(\xi_{x})_{x}+(\xi_{x})_{y}
=-(\frac{a_{12}}{a_{22}})_{x}\xi_{x}, \text{ }\text{ }\text{
}\xi_{x}(x,0)=1.
\end{equation*}
As above, let $(x(t),y(t))$ be the parameterization of an
arbitrary characteristic, then $\xi_{x}(t)=\xi_{x}(x(t),y(t))$
satisfies $\dot{\xi}_{x}=-(\frac{a_{12}}{a_{22}})_{x}\xi_{x}$,
$\xi_{x}(0)=1$.  By the mean value theorem
\begin{equation*}
|\xi_{x}(t)-1|=|\xi_{x}(t)-\xi_{x}(0)|\leq
y_{0}\sup_{X}|(\frac{a_{12}}{a_{22}})_{x}|\text{
}\sup_{X}|\xi_{x}|.
\end{equation*}
By property $(ii)$ for the $b_{ij}$,
\begin{equation*}
1-\varepsilon C_{1} y_{0} \sup_{X}|\xi_{x}|\leq \xi_{x}(t)\leq
\varepsilon C_{1} y_{0} \sup_{X}|\xi_{x}|+1.
\end{equation*}
Since this holds for any characteristic, we obtain
\begin{equation*}
\sup_{X}|\xi_{x}|\leq \frac{1}{1-\varepsilon C_{1} y_{0}}:=C_{2}.
\end{equation*}
It follows from (2.2) that
\begin{equation*}
\sup_{X}|\xi_{y}|\leq C_{3},
\end{equation*}
where $C_{2}$, $C_{3}$ are independent of $\varepsilon$ and
$b_{ij}$.  In order to estimate $\xi_{xx}$, differentiate (2.2)
two times with respect to $x$:
\begin{equation*}
(\frac{a_{12}}{a_{22}})(\xi_{xx})_{x}+(\xi_{xx})_{y}
=-2(\frac{a_{12}}{a_{22}})_{x}\xi_{xx}-
(\frac{a_{12}}{a_{22}})_{xx}\xi_{x}, \text{ }\text{ }\text{
}\xi_{xx}(x,0)=0.
\end{equation*}
Then the same procedure as above yields
\begin{equation*}
\sup_{X}|\xi_{xx}|\leq \varepsilon C_{4} y_{0}
\sup_{X}|\xi_{xx}|+\varepsilon C_{5} y_{0},
\end{equation*}
implying that
\begin{equation*}
\sup_{X}|\xi_{xx}|\leq \frac{\varepsilon C_{5}
y_{0}}{1-\varepsilon C_{4} y_{0}}:=\varepsilon C_{6}.
\end{equation*}
Furthermore, using the above estimates we can differentiate (2.2)
to obtain
\begin{equation*}
\sup_{X}|\xi_{xy}|\leq\varepsilon C_{7},\text{ }\text{ }\text{ }
\sup_{X}|\xi_{yy}|\leq\varepsilon C_{8},
\end{equation*}
for some constants $C_{7}$, $C_{8}$ independent of $\varepsilon$
and $b_{ij}$.  This procedure may be continued to yield,
\begin{equation*}
|\partial^{\alpha}\xi|\leq \varepsilon C_{9},
\end{equation*}
for any multi-index $\alpha$ satisfying $2\leq|\alpha|\leq
10$.\par
   We now show that $\overline{a}_{11}$, $\overline{a}_{22}$,
$\overline{a}_{1}$, $\overline{a}_{2}$ satisfy properties $(i)$,
$(ii)$, $(iii)$ and have the desired form.  Calculation shows
that,
\begin{equation*}
\overline{a}_{11}=a_{11}\xi_{x}^{2}+2a_{12}\xi_{x}\xi_{y}+a_{22}\xi_{y}^{2},\text{
}\text{ }\text{
}\overline{a}_{1}=a_{11}\xi_{xx}+2a_{12}\xi_{xy}+a_{22}\xi_{yy}+a_{1}\xi_{x}
+a_{2}\xi_{y}.
\end{equation*}
Furthermore, according to the above estimates and the fact that
the $b_{ij}$ vanish in a neighborhood of $\partial X$, we may
write
\begin{equation*}
\xi_{x}=1+\varepsilon\chi,
\end{equation*}
where $\chi\in C^{r-1}(\overline{X})$ vanishes in a neighborhood
of the lines $x=\pm x_{0}$.  It follows that,
\begin{equation*}
\overline{a}_{11}=-\eta^{2}+\varepsilon\overline{b}_{11},\text{ }
\text{ }\text{ }\overline{a}_{1}=\varepsilon\overline{b}_{1},
\end{equation*}
where $\overline{b}_{11}$ and $\overline{b}_{1}$ satisfy
properties $(i)$, $(ii)$, $(iii)$.  Moreover, since
$\overline{a}_{22}=a_{22}$ and $\overline{a}_{2}=a_{2}$,
properties $(i)$, $(ii)$, $(iii)$ hold for these coefficients as
well.\qed
   For the remainder of this section and section 3, $(\xi,\eta)$ will be the
coordinates of the plane.  For simplicity of notation we put
$x=\xi$, $y=\eta$, and $a_{ij}=\overline{a}_{ij}$,
$a_{i}=\overline{a}_{i}$, $a=\overline{a}$,
$b_{ij}=\overline{b}_{ij}$, $b_{i}=\overline{b}_{i}$,
$b=\overline{b}$.\par
   In order to obtain a well-posed boundary value problem, we will
study a regularization of $L$ in the infinite strip
$\Omega=\{(x,y)\mid |x|<x_{0}\}$.  More precisely, define the
operator
\begin{equation*}
L_{\theta}^{'}=-\theta \partial_{xxyy}+L,
\end{equation*}
where $\theta >0$ is a small constant that will tend to zero in
the Nash-Moser iteration procedure.  Furthermore, we will need to
modify some of the coefficients of $L$ away from $X$ as follows.
First cut $b_{ij}$, $b_{i}$, and $b$ off near the lines $y=\pm
y_{0}$, so that by property $(ii)$ of lemma 2.2 these functions
vanish in a neighborhood of $\partial X$, and the coefficients
$a_{ij}$, $a_{i}$, and $a$ are now defined on all of $\Omega$.
Choose values $y_{1}$, $y_{2}$, and $y_{3}$ such that
$y_{0}<y_{1}<y_{2}<y_{3}$, and let $\delta>0$ be a small constant
that depends on $y_{2}-y_{1}$ and $y_{3}-y_{2}$. Then redefine the
coefficient $a$ in the domain $\Omega-X$ so that:\bigskip

$i)$ $a\in C^{r-2}(\overline{\Omega})$,\bigskip

$ii)$ $a\equiv 1$ if $|y|\geq y_{1}$,\bigskip

$iii)$ $a\geq 0$ for $|y|\geq y_{0}$,\bigskip

$iv)$ $\partial_{y}a\geq 0$ if $y\geq y_{0}$, and
$\partial_{y}a\leq 0$ if $y\leq -y_{0}$.\bigskip

\noindent Redefine $a_{11}$ in $\Omega-X$ and near
$\partial\Omega$ so that:\bigskip

$i)$ $a_{11}\in C^{r-2}(\overline{\Omega})$,\bigskip

$ii)$ $a_{11}=\begin{cases}
-y^{2} & \text{if }y_{0}\leq |y|\leq y_{1},\\
-(\frac{y_{1}+y_{2}}{2})^{2} & \text{if }|y|\geq y_{2},
\end{cases}$\bigskip

$iii)$ $\partial_{y}a_{11}<0$ if $y\geq y_{0}$, and
$\partial_{y}a_{11}>0$ if $y\leq -y_{0}$,\bigskip

$iv)$ $\sup_{\Omega}\partial_{yy}a_{11}\leq \delta$,\bigskip

$v)$ $a_{11}|_{\partial\Omega}\leq-\theta$,
$\partial_{x}^{\alpha}a_{11}|_{\partial\Omega}=0$, $\alpha\leq
r-2$, and
$\sup_{\Omega}|\partial_{x}^{\beta}a_{11}|\leq\varepsilon\Lambda^{'}$,
$1\leq\beta\leq 8$.\bigskip

\noindent Lastly, redefine $a_{2}$ in $\Omega-X$ so that:\bigskip

$i)$ $a_{2}\in C^{r-2}(\overline{\Omega})$,\bigskip

$ii)$ $a_{2}=\begin{cases}
0 & \text{if }y_{0}\leq |y|\leq y_{2},\\
-\delta y+\delta(\frac{y_{2}+y_{3}}{2}) & \text{if }y\geq y_{3},\\
-\delta y-\delta(\frac{y_{2}+y_{3}}{2}) & \text{if }y\leq -y_{3},
\end{cases}$\bigskip

$iii)$ $a_{2}\leq 0$ if $y\geq y_{2}$, and $a_{2}\geq 0$ if $y\leq
-y_{2}$,\bigskip

$iv)$ $\sup_{|y|\geq y_{2}}|\partial_{y}a_{2}|\leq\delta$.\bigskip

\noindent Denote the operator $L$ with coefficients modified as
above by $L'$, and define
\begin{equation*}
L_{\theta}=-\theta\partial_{xxyy}+L'.
\end{equation*}
Note that since we are studying a local problem, as stated in the
introduction, we may modify the coefficients of the linearization
away from a fixed neighborhood of the origin.  This will become
clear in the final section, where a modified version of the
Nash-Moser iteration scheme is used.\par
   Consider the following boundary value problems
\begin{equation}
L_{\theta}u=f \text{ }\text{ }\text{ in }\text{ $\Omega$, }\text{
}\text{ } u|_{\partial\Omega}=0,
\end{equation}

\begin{equation}
\text{ }L_{\theta}u=f \text{ }\text{ }\text{ }\text{ in }\text{
$\Omega$, }\text{ }\text{ } u_{x}|_{\partial\Omega}=0,
\end{equation}

\noindent and the corresponding adjoint problems
\begin{equation}
L_{\theta}^{*}v=g \text{ }\text{ }\text{ in }\text{ $\Omega$,
}\text{ } v|_{\partial\Omega}=0,
\end{equation}

\begin{equation}
\text{ }L_{\theta}^{*}v=g \text{ }\text{ }\text{ }\text{ in
}\text{ $\Omega$, }\text{ } v_{x}|_{\partial\Omega}=0,
\end{equation}

\noindent where $L_{\theta}^{*}$ is the formal adjoint of
$L_{\theta}$.  The main result of this section is to obtain weak
solutions for all four problems.\par
   We will make extensive use of the following function spaces.
For $m,n\in \mathbb{Z}_{\geq 0}$ let
\begin{eqnarray*}
C^{(m,n)}(\overline{\Omega})&=&\{u:\Omega\rightarrow
\mathbb{R}\mid
\partial_{x}^{\alpha}\partial_{y}^{\beta}u\in
C^{0}(\overline{\Omega}),\text{ $\alpha\leq m$, $\beta\leq n$}\},\\
\widetilde{C}^{(m,n)}(\overline{\Omega})&=&\{u\in
C^{(m,n)}(\overline{\Omega})\mid u|_{\partial\Omega}=0,\text{ $u$
has bounded support}\},\\
\widetilde{C}_{x}^{(m,n)}(\overline{\Omega})&=&\{u\in
C^{(m,n)}(\overline{\Omega})\mid u_{x}|_{\partial\Omega}=0,\text{
$u$ has bounded support}\}.
\end{eqnarray*}
Define the norm
$$
\parallel
u\parallel_{(m,n)}=(\sum_{\alpha\leq m,\text{ }\beta\leq
n}\parallel
\partial_{x}^{\alpha}\partial_{y}^{\beta}u\parallel_{L^{2}(\Omega)}^{2})^{1/2},
$$
\noindent and let $\widetilde{H}^{(m,n)}(\Omega)$ and
$\widetilde{H}_{x}^{(m,n)}(\Omega)$ be the respective closures of
$\widetilde{C}^{(m,n)}(\overline{\Omega})$ and
$\widetilde{C}_{x}^{(m,n)}(\overline{\Omega})$ in the norm
$\parallel\cdot\parallel_{(m,n)}$.  Furthermore, let
$H^{m}(\Omega)$ denote the Sobolev space of square integrable
derivatives up to and including order $m$, with norm
$\parallel\cdot\parallel_{m}$.  Denote the $L^{2}(\Omega)$ inner
product and norm by $(\cdot,\cdot)$ and $\parallel\cdot\parallel$
respectively, and define the negative norm
\begin{equation*}
\parallel u\parallel_{(-m,-n)}=\sup_{v\in
\widetilde{H}^{(m,n)}(\Omega)}\frac{|(u,v)|}{\text{ }\text{
}\text{ }\text{ }\text{ }\parallel v\parallel_{(m,n)}}.
\end{equation*}
\noindent  Let $\widetilde{H}^{(-m,-n)}(\Omega)$ be the closure of
$L^{2}(\Omega)$ in the norm $\parallel\cdot\parallel_{(-m,-n)}$,
then $\widetilde{H}^{(-m,-n)}(\Omega)$ is the dual space of
$\widetilde{H}^{(m,n)}(\Omega)$.  The dual space of
$\widetilde{H}_{x}^{(m,n)}(\Omega)$ is defined similarly.\par
   Let $f\in L^{2}(\Omega)$.  A function $u\in L^{2}(\Omega)$ is said to be
a weak solution of (2.3) (respectively (2.4)) if
\begin{equation*}
(u,L_\theta^{*}v)=(f,v), \text{ }\text{ for }\text{ all }\text{
}v\in \widetilde{C}^{\infty}(\overline{\Omega})\text{ }\text{ (for
}\text{ all }v\in \widetilde{C}_{x}^{\infty}(\overline{\Omega})).
\end{equation*}
We shall employ the energy integral method, developed by K. O.
Friedrichs and others, to prove the existence of weak solutions
for (2.3) and (2.4).  The first step is to establish an a priori
estimate.\medskip

\noindent\textbf{Lemma 2.3 (Basic Estimate).}  \textit{If
$\varepsilon$, $\theta$, and $\delta$ are sufficiently small, then
there exist constants $C_{1},C_{2}>0$ independent of
$\varepsilon$, $\theta$, $\delta$, and functions $A,B,C,D,E\in
C^{\infty}(\overline{\Omega})$ where $E>0$ and $E=O(|y|)$ as
$|y|\rightarrow\infty$, such that:}
\begin{equation*}
(Au+Bu_{x}+Cu_{y}+Du_{yy},L_{\theta}u)\geq
\end{equation*}
\begin{equation*}
C_{1}[\parallel u\parallel^{2}+\parallel E
u_{y}\parallel^{2}+\theta(\parallel u_{x}\parallel^{2}+\parallel
u_{xy}\parallel^{2}+\parallel u_{yy}\parallel^{2}+\theta\parallel
u_{xyy}\parallel^{2})],
\end{equation*}
\textit{for all $u\in C^{\infty}(\overline{\Omega})$ with bounded
support such that $u_{x}(-x_{0},y)=0$, and either $u(x_{0},y)=0$
or $u_{x}(x_{0},y)=0$.  Furthermore,}
\begin{equation*}
\parallel u\parallel+\parallel
u_{y}\parallel+\sqrt{\theta}(\parallel u_{x}\parallel+\parallel
u_{xy}\parallel+\parallel u_{yy}\parallel+\sqrt{\theta}\parallel
u_{xyy}\parallel)\leq C_{2}\parallel L_{\theta}u\parallel,
\end{equation*}
\textit{for all $u\in \widetilde{C}^{\infty}(\overline{\Omega})$
and for all $u\in \widetilde{C}^{\infty}_{x}(\overline{\Omega})$}.

\medskip

\noindent\textit{Proof.}  We first define the functions $A,B,C$
and $D$. Let $\mu$ be a positive constant such that
$\frac{1}{4}\mu+ a_{11}\geq 1$ throughout $\Omega$, and let
$\gamma\in C^{\infty}([-x_{0},x_{0}])$ be such that
$$
\gamma(x)=
\begin{cases}
1 & \text{if $-x_{0}\leq x\leq\frac{x_{0}}{2}$},\\
0 & \text{if $x=x_{0}$},
\end{cases}
$$

\noindent with $\gamma(x)>0$ except at $x=x_{0}$, and
$\gamma^{'}\leq 0$.  Define
\begin{equation*}
A=\frac{1}{2}\partial_{y}C-a_{11},
\end{equation*}
\begin{equation*}
B=-\theta\gamma,
\end{equation*}
\begin{equation*}
C=\begin{cases}
\mu\partial_{y}a_{11} & \text{if $|y|<y_{0}$},\\
-2\mu y & \text{if $|y|\geq y_{0}$},
\end{cases}
\end{equation*}
\begin{equation*}
D=\theta,
\end{equation*}

\noindent and note that $A,B,C,D\in
C^{\infty}(\overline{\Omega})$.\par
   We now prove the first estimate.  Let $u\in C^{\infty}(\overline{\Omega})$ satisfy
the given hypotheses.  Let $(n_{1},n_{2})$ denote the unit outward
normal to $\partial\Omega$.  Then integrate by parts to obtain:
\begin{equation*}
(Au+Bu_{x}+Cu_{y}+Du_{yy},L_{\theta}u)=
\end{equation*}
\begin{equation*}\int\int_{\Omega}I_{1}u_{xyy}^{2}+
I_{2}u_{yy}^{2}+2I_{3}u_{yy}u_{xy}+I_{4}u_{xy}^{2}+2I_{5}u_{xy}u_{xx}
\end{equation*}
\begin{equation*}
+2I_{6}u_{xy}u_{y}+I_{7}u_{x}^{2}+2I_{8}u_{x}u_{y}+I_{9}u_{y}^{2}+I_{10}u^{2}
\end{equation*}
\begin{equation*}
+
\int_{\partial\Omega}J_{1}u_{xy}^{2}+J_{2}u_{xy}u_{x}+J_{3}u_{x}^{2}+
J_{4}u_{y}^{2}+J_{5}u^{2};\text{ }\text{ }\text{ }\text{ }\text{
}\text{ }\text{ }\text{ }\text{ }\text{ }\text{ }\text{ }
\end{equation*}

\noindent where

\begin{equation*}
\text{ }\text{ } J_{1}=\frac{1}{2}\theta Bn_{1},\text{ }\text{
}\text{ }\text{ }J_{2}=\theta B_{y}n_{1},\text{ }\text{ }\text{
}\text{ }J_{3}=\frac{1}{2}Ba_{11}n_{1},
\end{equation*}
\begin{equation*}
2J_{4}=-\theta A_{x}n_{1}-\theta C_{xy}n_{1}
+(Da_{11})_{x}n_{1}-Da_{1}n_{1},
\end{equation*}
\begin{equation*}
\! 2J_{5}=-(Aa_{11})_{x}n_{1}+Aa_{1}n_{1}+Ban_{1}+\theta
A_{xyy}n_{1},
\end{equation*}
and the remaining $I_{1},\ldots,I_{10}$ will be given below as
each term is estimated.  First note that
$J_{2}|_{\partial\Omega}=J_{4}|_{\partial\Omega}\equiv 0$.
Furthermore $J_{1}=\cdots=J_{5}\equiv 0$ on the portion of the
boundary, $x=x_{0}$, since $\gamma(x_{0})=0$.  Whereas on the
other half of the boundary, $x=-x_{0}$, we have
$u_{x}(-x_{0},y)=0$ and $J_{5}=\frac{1}{2}Ban_{1}\geq 0$.  It
follows that the entire boundary integral is nonnegative.\par
   We now proceed to estimate the integral over $\Omega$, beginning
with $I_{1}$, $I_{5}$, and $I_{10}$, which are given by
\begin{equation*}
\text{ }\text{ }\text{ }\text{ }I_{1}=\theta D,\text{ }\text{
}\text{ }\text{ }\text{ }\text{ }I_{5}=-\frac{1}{2}\theta
B_{y},\text{ }\text{ }
\end{equation*}
\begin{equation*}
2I_{10}=(Aa_{11})_{xx}+(Aa_{22})_{yy}- (Aa_{1})_{x}-(Aa_{2})_{y}
\end{equation*}
\begin{equation*}
\text{ }\text{ }\text{ }\text{ }\text{ }\text{ }\text{ }\text{
}\text{ } +2Aa -(Ca)_{y}-(Ba)_{x}-\theta A_{xxyy}+(Da)_{yy}.
\end{equation*}
Since $B$ is a function of $x$ alone, $I_{5}\equiv 0$, and by
definition of $D$, $I_{1}=\theta^{2}$.  It will now be shown that
$I_{10}\geq M_{1}$ in $\Omega$, for some constant $M_{1}>0$
independent of $\varepsilon$ and $\theta$.  In order to accomplish
this we shall treat the regions $|y|\leq y_{0}$, $y_{0}\leq
|y|\leq y_{1}$,\linebreak $y_{1}\leq |y|\leq y_{2}$, and $|y|\geq
y_{2}$ separately.  Moreover, throughout this proof
$M_{i}$,\linebreak $i=1,2,\ldots$, will always denote positive
constants independent of $\varepsilon$ and $\theta$.  A
computation yields,
\begin{equation*}
I_{10}=-a_{22}\partial_{yy}a_{11}-
a_{11}a-\frac{1}{2}C\partial_{y}a-\frac{1}{2}(Aa_{2})_{y}
+O(\varepsilon+\theta).
\end{equation*}
In the region $|y|\leq y_{0}$ we have
$a,\partial_{y}a,a_{2},\partial_{y}a_{2}=O(\varepsilon)$,
$a_{22}=1+O(\varepsilon)$, and\linebreak $\partial_{yy}
a_{11}=-2+O(\varepsilon)$, so that here $I_{10}\geq M_{2}$.  If
$y_{0}\leq |y|\leq y_{1}$, the conditions placed on $a$ guarantee
that
\begin{equation*}
-a_{11}a-\frac{1}{2}C\partial_{y}a\geq 0;
\end{equation*}
furthermore $a_{22}$, $a_{11}$, and $a_{2}$ have the same
properties in this region as in the previous.  Hence, $I_{10}\geq
M_{3}$ when $y_{0}\leq |y|\leq y_{1}$.  If $y_{1}\leq |y|\leq
y_{2}$ then
\begin{equation*}
-a_{22}\partial_{yy}a_{11}=O(\delta),\text{ }\text{ }\text{
}\text{ }\text{ } \text{ }-a_{11}a\geq y_{1}^{2},\text{ }\text{
}\text{ }\text{ }\text{ }\text{ }\text{ }\text{ }\text{
}a_{2}=\partial_{y}a\equiv 0,
\end{equation*}
showing that $I_{10}\geq M_{4}$ in this region.  Lastly, when
$|y|\geq y_{2}$ we have $I_{10}\geq M_{5}$ since
\begin{equation*}
\partial_{yy}a_{11}=\partial_{y}a\equiv 0,\text{ }\text{ }
\text{ }\text{ }\text{ }
-a_{11}a=(\frac{y_{1}+y_{2}}{2})^{2},\text{ }\text{ }\text{
}\text{ }\text{ }\text{ } -\frac{1}{2}(Aa_{2})_{y}=O(\delta).
\end{equation*}
The desired conclusion now follows by combining the above
estimates.\par
   Next we show that
\begin{equation*}
\int\int_{\Omega}I_{2}u_{yy}^{2}+2I_{3}u_{yy}u_{xy}+I_{4}u_{xy}^{2}\geq
M_{6}\theta(\parallel u_{yy}\parallel^{2}+\parallel
u_{xy}\parallel^{2}),
\end{equation*}
where
\begin{equation*}
\text{ }I_{2}=-\frac{1}{2}\theta D_{xx}+Da_{22}, \text{ }\text{
}\text{ }\text{ }\text{ }I_{3}=-\frac{1}{2}\theta C_{x},
\end{equation*}
\begin{equation*}
I_{4}=-\frac{1}{2}\theta C_{y}-\frac{1}{2}\theta B_{x}-\theta
A+Da_{11}.
\end{equation*}
This will follow if $I_{2}\geq M_{7}\theta$, $I_{4}\geq
M_{8}\theta$, and $I_{2}I_{4}-I_{3}^{2}>0$.  A calculation shows
that
\begin{equation*}
I_{2}=\theta a_{22}=\theta(1+O(\varepsilon)),\text{ }\text{
}\text{ }\text{ }\text{ }\text{ }\text{
}I_{3}=O(\varepsilon\theta),
\end{equation*}
\begin{equation*}
I_{4}=2\theta(
a_{11}-\frac{1}{2}C_{y})+O(\varepsilon\theta)=2\theta(\mu+
a_{11}+O(\varepsilon)).
\end{equation*}
Therefore, since $\mu$ was chosen so that $\mu+a_{11}\geq 1$ in
$\Omega$, the desired conclusion follows if $\varepsilon$ is
sufficiently small.\par
   We now show that
\begin{equation*}
\int\int_{\Omega}I_{7}u_{x}^{2}+2I_{8}u_{x}u_{y}+I_{9}u_{y}^{2}\geq
M_{9}(\theta\parallel u_{x}\parallel^{2}+\parallel E
u_{y}\parallel^{2}),
\end{equation*}
where
\begin{equation*}
\text{ }\text{ }
2I_{7}=-2Aa_{11}-(Ba_{11})_{x}+2Ba_{1}+(Ca_{11})_{y}+\theta
B_{xyy}+\theta A_{yy}-(Da_{11})_{yy},
\end{equation*}
\begin{equation*}
\text{ }\text{ } 2I_{8}=-(Ba_{22})_{y}+Ba_{2}-(Ca_{11})_{x}
+Ca_{1}+\theta A_{xy}+(Da_{11})_{xy}-(Da_{1})_{y},\text{ }
\end{equation*}
\begin{equation*}
\text{ }\text{ } 2I_{9}=-2Aa_{22}-(Ca_{22})_{y}+ 2Ca_{2}+\theta
C_{xxy}+\theta A_{xx}\text{ }\text{ }\text{ }\text{ }\text{
}\text{ }\text{ }\text{ }\text{ }\text{ }\text{ }\text{ }\text{
}\text{ }\text{ }\text{ }\text{ }\text{ }\text{ }\text{ }\text{
}\text{ }\text{ }\text{ }\text{ }\text{ }\text{ }\text{ }\text{ }
\end{equation*}
\begin{equation*}
\text{ }\text{ }\text{ }\text{ }\text{
}-(Da_{11})_{xx}-(Da_{2})_{y}+(Da_{1})_{x}-2Da.\text{ }\text{ }
\text{ }\text{ }\text{ }\text{ }\text{ }\text{ }\text{ }\text{
}\text{ }\text{ }\text{ }\text{ }\text{ }\text{ }\text{ }\text{
}\text{ }\text{ }\text{ }\text{ }\text{ }\text{ }\text{ }\text{
}\text{ }\text{ }\text{ }
\end{equation*}
Again, this will follow if $I_{7}\geq M_{10}\theta$, $I_{9}\geq
M_{11}E^{2}$, and $I_{7}I_{9}-I_{8}^{2}>0$.  A calculation shows
that
\begin{equation*}
I_{7}=a_{11}^{2}+\frac{1}{2}C\partial_{y}a_{11}
+\theta(-\partial_{yy}a_{11}+\frac{1}{2}\gamma_{x}a_{11}+O(\varepsilon)),
\end{equation*}
\begin{equation*}
I_{8}=-\frac{1}{2}C_{x}a_{11}-\frac{1}{2}C\partial_{x}a_{11}
+\frac{1}{2}Ca_{1}+\frac{1}{2}Ba_{2}+O(\theta),
\end{equation*}
\begin{equation*}
I_{9}=(a_{11}-C_{y})a_{22}+Ca_{2}+O(\varepsilon+\theta)\text{ }
\text{ }\text{ }\text{ }\text{ }\text{ }\text{ }\text{ }\text{
}\text{ }\text{ }\text{ }\text{ }\text{ }\text{ }\text{ }\text{
}\text{ }
\end{equation*}
\begin{equation*}
\text{ }\text{ }\text{ }\text{ }\text{
}=(2\mu+a_{11}+O(\varepsilon))(1+O(\varepsilon))+Ca_{2}
+O(\varepsilon+\theta).
\end{equation*}
Then $I_{9}\geq M_{11}E^{2}$ immediately follows since
$Ca_{2}=O(\varepsilon)$ if $|y|\leq y_{0}$, $Ca_{2}\geq 0$ if
$|y|\geq y_{0}$, $Ca_{2}=O(|y|^{2})$ as $|y|\rightarrow\infty$,
and $2\mu+a_{11}\geq 1$.  To show that $I_{7}\geq M_{10}\theta$,
we consider the regions $|y|\leq y_{0}$ and $|y|\geq y_{0}$
separately.  If $|y|\leq y_{0}$ then
\begin{equation*}
C\partial_{y}a_{11}=\mu(\partial_{y}a_{11})^{2}\geq 0,\text{
}\text{ }\text{ }\text{ }
-\partial_{yy}a_{11}=2+O(\varepsilon),\text{ }\text{ }\text{
}\text{ }\text{ }\text{ } \gamma_{x}a_{11}\geq-O(\varepsilon),
\end{equation*}
so that here $I_{7}\geq 2\theta+O(\varepsilon\theta)$.
Furthermore, when $|y|\geq y_{0}$ we have $I_{7}\geq
y_{0}^{4}+O(\theta)$ since
\begin{equation*}
a_{11}^{2}\geq y_{0}^{4},\text{ }\text{ } \text{ }\text{ }\text{
}\text{ }\text{ }C\partial_{y}a_{11}\geq 0.
\end{equation*}
Finally, $I_{7}I_{9}-I_{8}^{2}>0$ follows from the next
calculation.  If $|y|\leq y_{0}$ then
\begin{equation*}
I_{7}I_{9}-I_{8}^{2}\geq(a_{11}^{2}
+\frac{\mu}{2}(\partial_{y}a_{11})^{2}+2\theta
+O(\varepsilon\theta))(1+O(\varepsilon+\theta))\text{ }\text{
}\text{ }\text{ }\text{ }\text{ }\text{ }\text{ }\text{ }\text{
}\text{ }
\end{equation*}
\begin{equation*}
\text{ }\text{ }\text{ }\text{
}-\frac{1}{4}O(\varepsilon^{2})a_{11}^{2}
-\frac{1}{4}O(\varepsilon^{2})(\partial_{y}
a_{11})^{2}-O(\varepsilon\theta+\theta^{2}),
\end{equation*}
whereas if $|y|\geq y_{0}$ then
\begin{equation*}
I_{7}I_{9}-I_{8}^{2}\geq(y_{0}^{4} +O(\theta))(1+O(\delta
y^{2}))-O(\theta^{2}y^{2}).
\end{equation*}\par
   Lastly, we deal with the term $2I_{6}u_{xy}u_{y}$.  Consider
the quadratic form:
\begin{equation*}
M_{6}\theta u_{xy}^{2}+2I_{6}u_{xy}u_{y}+M_{9}E^{2}u_{y}^{2},
\end{equation*}
where
\begin{equation*}
I_{6}=-\frac{1}{2}Ba_{22}.
\end{equation*}
Since
\begin{equation*}
(M_{2}\theta)(M_{3}E^{2})-I_{6}^{2}\geq
M_{11}\theta-M_{12}\theta^{2} (1+O(\varepsilon))
\end{equation*}
for some $M_{11}$, $M_{12}$, we obtain
\begin{equation*}
M_{6}\theta u_{xy}^{2}+2I_{6}u_{xy}u_{y}+M_{9}E^{2}u_{y}^{2}\geq
M_{13}(\theta u_{xy}^{2}+E^{2}u_{y}^{2}).
\end{equation*}
This completes the proof of the first estimate.\par
   To obtain the second estimate we need only observe that the above
arguments hold if $B\equiv 0$ and
$u\in\widetilde{C}^{\infty}(\overline{\Omega})$ or
$u\in\widetilde{C}^{\infty}_{x}(\overline{\Omega})$.  Then an
application of Cauchy's inequality ($ab\leq\lambda
a^{2}+\frac{1}{4\lambda}b^{2}$, $\lambda>0$) yields the desired
result.  The reason for including $B$ in the first estimate will
soon become clear.\qed
   Having established the basic estimate, our goal shall now be to
establish dual inequalities of the form:
\begin{equation*}
\parallel v\parallel\leq C_{1}\parallel
L_{\theta}^{*}v\parallel_{(-1,-2)} \text{ }\text{ for all }\text{
} v\in \widetilde{C}^{\infty}(\overline{\Omega}),
\end{equation*}
\begin{equation*}
\parallel v\parallel\leq C_{2}\parallel
L_{\theta}^{*}v\parallel_{(-1,-2)}\text{ }\text{ for all }\text{
}v\in \widetilde{C}_{x}^{\infty}(\overline{\Omega}).
\end{equation*}
The existence of weak solutions to problems (2.3) and (2.4) will
then easily follow from these two dual estimates, respectively. In
order to establish the dual estimates, we will need the following
lemma.  Let $P$ denote the differential operator
\begin{equation*}
P=D\partial_{y}^{2}+B\partial_{x}+C\partial_{y}+A,
\end{equation*}
where $A,B,C,$ and $D$ are defined in lemma 2.3.  Note that $P$ is
parabolic in $\Omega$, away from the portion of the boundary,
$x=x_{0}$.  This is the reason for including $B$ in the first
estimate of lemma 2.3.\medskip

\noindent\textbf{Lemma 2.4.}  \textit{For every $v\in
\widetilde{C}^{\infty}(\overline{\Omega})$ there exists a unique
solution $u\in C^{\infty}(\Omega)\cap H^{4}(\Omega)\subset
C^{\infty}(\Omega)\cap C^{2}(\overline{\Omega})$ of
\begin{equation*}
Pu=v\text{ }\textit{ in }\text{ }\Omega,\text{ }\text{
}u(-x_{0},y)=u_{x}(-x_{0},y)=0,\text{ }\text{ }u(x_{0},y)=0.
\end{equation*}
Furthermore, for every $v\in
\widetilde{C}_{x}^{\infty}(\overline{\Omega})$ there exists a
unique solution $u\in C^{\infty}(\Omega)\cap H^{4}(\Omega)\subset
C^{\infty}(\Omega)\cap C^{2}(\overline{\Omega})$ of
\begin{equation*}
Pu=v\text{ }\textit{ in }\text{ }\Omega,\text{ }\text{
}u_{x}(-x_{0},y)=0,\text{ }\text{ }u_{x}(x_{0},y)=0.
\end{equation*}}\medskip
\noindent\textit{Proof.}  Let $\tau>0$ be a small parameter, and
define the subdomains
\begin{equation*}
\Omega_{\tau}=\{(x,y)\mid -x_{0}<x<x_{0}-\tau\}.
\end{equation*}
Then $P$ is parabolic in $\overline{\Omega}_{\tau}$ for each
$\tau$.\par
   We now consider the case when $v\in
\widetilde{C}^{\infty}(\overline{\Omega})$.  The parabolicity of
$P$ guarantees the existence (see [13]) of a unique solution to
the Cauchy problem
\begin{equation*}
Pu=v \text{ }\text{ in }\text{ }\Omega,\text{ }\text{
}u(-x_{0},y)=0,
\end{equation*}
such that $u\in H^{\infty}(\Omega_{\tau})$ for every $\tau$.
Furthermore, $u_{x}(-x_{0},y)=0$ since
\begin{equation*}
Bu_{x}|_{(-x_{0},y)}=Pu|_{(-x_{0},y)}=v(-x_{0},y)=0.
\end{equation*}\par
   We shall now show that $u\in H^{4}(\Omega)$.  This will be
accomplished by estimating the $H^{4}(\Omega_{\tau})$ norm of $u$
in terms of the $H^{4}(\Omega)$ norm of $v$, independent of
$\tau$.  To facilitate the estimates, we first construct an
appropriate approximating sequence, $\{u^{k}\}_{k=1}^{\infty}$,
for $u$. Define functions $\nu_{k}\in C^{\infty}(\mathbb{R})$ by
\begin{equation}
\nu_{k}(y)=\begin{cases}
1 & \text{if $|y|\leq k$},\\
0 & \text{if $|y|\geq 3k$,}
\end{cases}
\end{equation}
such that $0\leq\nu_{k}\leq 1$, $\sup|\nu_{k}'|\leq\frac{1}{k}$,
and $|\nu_{k}|_{C^{4}(\overline{\Omega})}\leq M$ for some constant
$M$ independent of $k$.  Let $u^{k}=\nu_{k}u$, then:\bigskip

$i)$ $u^{k}\in C^{\infty}(\overline{\Omega}_{\tau})$ for all
$\tau$,\bigskip

$ii)$ $u^{k}$ has bounded support and
$u^{k}(-x_{0},y)=u^{k}_{x}(-x_{0},y)=0$,\bigskip

$iii)$ $\parallel u-u^{k}\parallel_{4,\Omega_{\tau}}\rightarrow 0$
as $k\rightarrow\infty$,\bigskip

$iv)$ $\parallel
Cu_{y}-Cu^{k}_{y}\parallel_{\Omega_{\tau}}\rightarrow 0$ as
$k\rightarrow\infty$,\bigskip

\noindent where $C$ was defined in lemma 2.3.  All of the above
properties are evident except for $(iv)$, and $(iv)$ follows from
the following calculation.  Let
\begin{equation*}
\Omega_{\tau}^{(k_{1},k_{2})}=\{(x,y)\in\Omega_{\tau}\mid
k_{1}\leq |y|\leq k_{2}\},
\end{equation*}
then
\begin{eqnarray*}
\parallel Cu_{y}-Cu^{k}_{y}\parallel_{\Omega_{\tau}}^{2}&\leq&
\parallel C(u_{y}-\nu_{k}u_{y})\parallel^{2}+
\parallel C\nu_{k}'u\parallel^{2}\\
&\leq&\int\int_{\Omega_{\tau}^{(k,\infty)}}C^{2}u_{y}^{2}+
\int\int_{\Omega_{\tau}^{(k,3k)}}(C\nu_{k}')^{2}u^{2}\\
&\leq&\int\int_{\Omega_{\tau}^{(k,\infty)}}C^{2}u_{y}^{2}+
\int\int_{\Omega_{\tau}^{(k,3k)}}(6\mu
k)^{2}(\frac{1}{k})^{2}u^{2},
\end{eqnarray*}
where $\mu$ was defined in the proof of lemma 2.3.  By solving for
$Cu_{y}$ in the equation $Pu=v$, we have
\begin{equation*}
Cu_{y}=v-Du_{yy}-Bu_{x}-Au\in L^{2}(\Omega_{\tau}).
\end{equation*}
Therefore
\begin{equation*}
\int\int_{\Omega_{\tau}^{(k,\infty)}}C^{2}u_{y}^{2}\rightarrow 0
\text{ }\text{ }\text{ as }\text{ }\text{ }k\rightarrow\infty.
\end{equation*}
Furthermore,
\begin{equation*}
\int\int_{\Omega_{\tau}^{(k,3k)}}36\mu^{2}u^{2}\rightarrow 0
\text{ }\text{ }\text{ as }\text{ }\text{ }k\rightarrow\infty
\end{equation*}
since $u\in L^{2}(\Omega_{\tau})$.  This proves $(iv)$.\par
   We now proceed to estimate the $H^{4}(\Omega_{\tau})$ norm of
$u$.  Let $\zeta=\zeta(y)\in C^{\infty}(\mathbb{R})$ be such that
$\zeta<0$, $\zeta(y)=-|y|^{-1/2}$ if $|y|\geq y_{1}$,
$\zeta^{'}(y)\geq 0$ if $y\geq 0$, and $\zeta^{'}(y)\leq 0$ if
$y\leq 0$.  Then set $\kappa=2\sup|\zeta a_{11}|$, and integrate
by parts to obtain
\begin{eqnarray*}
\int\int_{\Omega_{\tau}}(\kappa u^{k}_{yy}+\zeta u^{k})Pu^{k}&=&
\int\int_{\Omega_{\tau}}[\kappa
D](u^{k}_{yy})^{2}+[-D\zeta+\kappa(\frac{1}{2}B_{x}
-\frac{1}{2}C_{y}-A)](u^{k}_{y})^{2}\\
& &+[\frac{1}{2}\kappa A_{yy}+\frac{1}{2}(D\zeta)_{yy}-
\frac{1}{2}(B\zeta)_{x}-\frac{1}{2}(C\zeta)_{y}+\zeta
A](u^{k})^{2}\\
& &+\int_{\partial\Omega_{\tau}}[-\frac{1}{2}\kappa Bn_{1}]
(u^{k}_{y})^{2} +[\frac{1}{2}B\zeta n_{1}](u^{k})^{2}.
\end{eqnarray*}
The boundary integral is nonnegative since $u^{k}(-x_{0},y)=
u^{k}_{y}(-x_{0},y)=0$, and\linebreak $-\kappa
Bn_{1}|_{(x_{0}-\tau,y)}, B\zeta n_{1}|_{(x_{0}-\tau,y)}>0$. Also,
$\kappa D>0$,
\begin{equation*}
-D\zeta+\kappa(\frac{1}{2}B_{x} -\frac{1}{2}C_{y}-A)\geq
\kappa(2\mu+a_{11}+O(\varepsilon+\theta))\geq \kappa,
\end{equation*}
and
\begin{equation*}
\frac{1}{2}\kappa A_{yy}+\frac{1}{2}(D\zeta)_{yy}-
\frac{1}{2}(B\zeta)_{x}-\frac{1}{2}(C\zeta)_{y}+\zeta A
\end{equation*}
\begin{equation*}
=-\frac{1}{2}\kappa\partial_{yy}a_{11}-\frac{1}{2}C\zeta_{y}
-\zeta a_{11}+\frac{1}{2}(D\zeta)_{yy}-\frac{1}{2}
(B\zeta)_{x}+O(\varepsilon)
\end{equation*}
$$
\geq\begin{cases}
\kappa-\zeta a_{11}+O(\varepsilon+\theta) & \text{if $|y|\leq y_{1}$},\\
|y|^{-1/2}[\frac{1}{2}\mu+a_{11}+O(\theta)]+O(\kappa\delta) &
\text{if $y_{1}\leq |y|\leq y_{2}$,}\\
|y|^{-1/2}[\frac{1}{2}\mu+a_{11}+O(\theta)] & \text{if $|y|\geq
y_{2}$.}
\end{cases}
$$
Therefore if $\varepsilon$, $\theta$, and $\delta$ are
sufficiently small, we may apply the Schwartz inequality followed
by Cauchy's inequality to obtain
\begin{equation*}
\parallel\sqrt{-\zeta}u^{k}\parallel_{\Omega_{\tau}}+\parallel
u^{k}_{y}\parallel_{\Omega_{\tau}}+\parallel
u^{k}_{yy}\parallel_{\Omega_{\tau}}\leq M_{1}\parallel
Pu^{k}\parallel_{\Omega_{\tau}},
\end{equation*}
for some constant $M_{1}$ independent of $\tau$.  The properties
of $u^{k}$ guarantee that by letting $k\rightarrow\infty$, we
obtain
\begin{equation*}
\parallel\sqrt{-\zeta}u\parallel_{\Omega_{\tau}}+\parallel
u_{y}\parallel_{\Omega_{\tau}}+\parallel
u_{yy}\parallel_{\Omega_{\tau}}\leq M_{1}\parallel
Pu\parallel_{\Omega_{\tau}}=M_{1}\parallel
v\parallel_{\Omega_{\tau}}\leq M_{1}\parallel v\parallel.
\end{equation*}\par
   We now estimate $\partial_{x}^{\alpha}\partial_{y}^{\beta}u$
for $\alpha=1,\ldots,4$, and $\beta=0,1,2$.  Differentiate $Pu=v$
with respect to $x$:
\begin{equation}
D(u_{x})_{yy}+B(u_{x})_{x}+C(u_{x})_{y}+(A+B_{x})u_{x}
=v_{x}-C_{x}u_{y}-A_{x}u.
\end{equation}
Since $u_{x}(-x_{0},y)=0$ and $A_{x}$, $C_{x}$ vanish outside a
compact set, we can apply the same procedure as above to obtain
\begin{eqnarray*}
\parallel\sqrt{-\zeta}u_{x}\parallel_{\Omega_{\tau}}+\parallel
u_{xy}\parallel_{\Omega_{\tau}}+\parallel
u_{xyy}\parallel_{\Omega_{\tau}}&\leq& M_{1}\parallel
v_{x}-C_{x}u_{y}-A_{x}u\parallel_{\Omega_{\tau}}\\
&\leq& M_{2}(\parallel v_{x}\parallel_{\Omega_{\tau}}+
\parallel u_{y}\parallel_{\Omega_{\tau}}+\parallel u
\parallel_{\Omega_{\tau}})\\
&\leq& M_{3}(\parallel v\parallel+\parallel v_{x}\parallel).
\end{eqnarray*}
Differentiating (2.8) with respect to $x$ produces
\begin{equation*}
D(u_{xx})_{yy}+B(u_{xx})_{x}+C(u_{xx})_{y}+(A+2B_{x})u_{xx}
\end{equation*}
\begin{equation*}
=v_{xx}-\partial_{x}(C_{x}u_{y}+A_{x}u)
-C_{x}u_{xy}-(A_{x}+B_{xx})u_{x}:=v_{1}.
\end{equation*}
Again we apply the same method.  However, since $u_{xx}(-x_{0},y)
=B^{-1}v_{x}|_{(-x_{0},y)}$ from (2.8), we now have
\begin{eqnarray*}
\parallel\sqrt{-\zeta}u_{xx}\parallel_{\Omega_{\tau}}\!+\!\parallel
u_{xxy}\parallel_{\Omega_{\tau}}\!+\!\parallel
u_{xxyy}\parallel_{\Omega_{\tau}}\!&\leq&\! M_{1}\parallel
v_{1}\parallel_{\Omega_{\tau}}\!+M_{4}\\
\!&\leq&\! M_{5}(\parallel v\parallel\!+\!\parallel v_{x}\parallel
\!+\!\parallel v_{xx}\parallel)\!+\!M_{4},
\end{eqnarray*}
where $M_{4}=\kappa|B|^{-1}(\int_{x=-x_{0}}v_{xy}^{2}
+v_{x}^{2})^{1/2}$ which is independent of $\tau$. We can estimate
$\parallel\sqrt{-\zeta}\partial_{x}^{\alpha}u\parallel_{\Omega_{\tau}}$,
$\alpha=3,4$, and
$\parallel\partial_{x}^{\alpha}\partial_{y}^{\beta}u\parallel_{\Omega_{\tau}}$,
$\alpha=3,4$, $\beta=1,2$, in a similar manner.\par
   To estimate $u_{yyy}$, differentiate $Pu=v$ with respect to
$y$:
\begin{equation}
D(u_{y})_{yy}+B(u_{y})_{x}+C(u_{y})_{y}+(A+C_{y})u_{y}
=v_{y}-A_{y}u.
\end{equation}
Since $u_{y}(-x_{0},y)=0$, $C_{y}<0$, and $A_{y}$ vanishes outside
a compact set, the same method as above yields
\begin{eqnarray*}
\parallel\sqrt{-\zeta}u_{y}\parallel_{\Omega_{\tau}}+\parallel
u_{yy}\parallel_{\Omega_{\tau}}+\parallel
u_{yyy}\parallel_{\Omega_{\tau}}&\leq& M_{1}\parallel
v_{y}-A_{y}u\parallel_{\Omega_{\tau}}\\
&\leq& M_{6}(\parallel v\parallel+\parallel v_{y}\parallel).
\end{eqnarray*}
Furthermore, $\parallel u_{xyyy}\parallel_{\Omega_{\tau}}$ and
$\parallel u_{yyyy}\parallel_{\Omega_{\tau}}$ can be estimated by
differentiating (2.9) with respect to $x$ and $y$,
respectively.\par
   The combination of all the above estimates produces,
\begin{equation*}
\sum_{\alpha=0}^{4}\parallel\sqrt{-\zeta}\partial_{x}^{\alpha}u\parallel_{\Omega_{\tau}}
+\sum_{\alpha+\beta\leq 4\atop\beta\neq
0}\parallel\partial_{x}^{\alpha}\partial_{y}^{\beta}u\parallel_{\Omega_{\tau}}
\leq M_{7}\parallel v\parallel_{4}+M_{8},
\end{equation*}
where $M_{7}$ and $M_{8}$ are independent of $\tau$.  Then letting
$\tau \rightarrow 0$ we find that\linebreak
$\partial_{x}^{\alpha}\partial_{y}^{\beta}u\in L^{2}(\Omega)$,
$\alpha+\beta\leq 4$, $\beta\neq 0$, and that
$\sqrt{-\zeta}\partial_{x}^{\alpha}u\in L^{2}(\Omega)$,
$\alpha=0,\ldots,4$.  It follows that $u\in H^{4}(K)$ for every
compact $K\subset\Omega$, so that $u\in
C^{2}(\overline{\Omega})$.\par
   We now show that $\partial_{x}^{\alpha}u\in L^{2}(\Omega)$,
$\alpha=0,\ldots,4$.  Let $\varrho_{1},\varrho_{2}\in
C^{\infty}(\mathbb{R})$ be given by
$$
\varrho_{1}(x)=
\begin{cases}
-B+\theta & \text{if  $-x_{0}\leq x\leq\frac{-x_{0}}{2}$},\\
0 & \text{if  $0\leq x\leq x_{0}$.}
\end{cases}
\text{ }\text{ }\text{ }\text{ }\text{ }
\varrho_{2}(y)=
\begin{cases}
-y & \text{if  $|y|\leq y_{0}$,}\\
0 & \text{if  $|y|\geq T$,}
\end{cases}
$$
such that $\varrho_{2}(y)\leq 0$ if $y>0$ and $\varrho_{2}(y)\geq
0$ if $y<0$, where $T>0$ is large enough so that
$-1\leq\varrho_{2}^{'}\leq\varepsilon$.  Then define
$\overline{B}=B+\varrho_{1}$ and
$\overline{C}=C+\varrho_{2}-\varepsilon\mu\partial_{y}
b_{11}=-2\mu y+\varrho_{2}$, and set
\begin{equation*}
\overline{P}=\overline{B}\partial_{x}+\overline{C}\partial_{y}+A.
\end{equation*}\par
   If $w\in C_{c}^{\infty}(\overline{\Omega})$, then integrating
by parts yields
\begin{equation*}
(w,\overline{P}^{*}w)=\int\int_{\Omega}[-\frac{1}{2}
\overline{B}_{x}-\frac{1}{2}\overline{C}_{y}+A]w^{2}
+\int_{\partial\Omega}[-\frac{1}{2}\overline{B}n_{1}]w^{2}.
\end{equation*}
The boundary integral is nonnegative since
$\overline{B}(-x_{0},y)=\theta$ and $\overline{B}(x_{0},y)=0$.
Furthermore,
\begin{equation*}
-\frac{1}{2} \overline{B}_{x}-\frac{1}{2}\overline{C}_{y}+A
=-\varrho_{2}^{'}-a_{11}+O(\varepsilon+\theta)\geq M_{9}
\end{equation*}
for some constant $M_{9}>0$.  Thus
\begin{equation}
\parallel w\parallel\leq M_{10}\parallel\overline{P}^{*}w\parallel.
\end{equation}
Since
$v-Du_{yy}+\varrho_{1}u_{x}+(\varrho_{2}-\varepsilon\mu\partial_{y}
b_{11})u_{y}\in L^{2}(\Omega)$, (2.10) implies (see the proof of
theorem 2.1 below) the existence of a weak solution
$\widetilde{u}\in L^{2} (\Omega)$ of
\begin{equation*}
\overline{P}\widetilde{u}=v-Du_{yy}+\varrho_{1}u_{x}+(\varrho_{2}-\varepsilon\mu\partial_{y}
b_{11})u_{y},\text{ }\text{ }\text{ }\widetilde{u}(-x_{0},y)=0.
\end{equation*}\par
   We shall now show that $u\equiv\widetilde{u}$.  Since
$\overline{P}$ is a first order differential operator, we may
apply G. Peyser's extension [21] of Friedrichs' result [2] on the
identity of weak and strong solutions to obtain a sequence
$\{\widetilde{u}^{k}\}_{k=1}^{\infty}$, such that
$\widetilde{u}^{k}\in C^{\infty}(\overline{\Omega})$ has bounded
support, satisfies $\widetilde{u}^{k}(-x_{0},y)=0$, and
\begin{equation*}
\parallel\widetilde{u}-\widetilde{u}^{k}\parallel+
\parallel\overline{P}\widetilde{u}^{k}-(v-Du_{yy}+\varrho_{1}u_{x}+(\varrho_{2}-\varepsilon
\mu\partial_{y}b_{11})u_{y})\parallel\rightarrow 0\text{ }\text{
}\text{ as }\text{ }\text{ }k\rightarrow\infty.
\end{equation*}
Set $v^{k}=u-\widetilde{u}^{k}$.  Using the fact that
$|y|^{-1/4}v^{k}\rightarrow|y|^{-1/4}(u-\widetilde{u})\in
L^{2}(\Omega)$ and recalling the definition of $\overline{P}$, we
have
\begin{eqnarray*}
|(-|y|^{-1/4}v^{k},\overline{P}v^{k})|&\leq&
\parallel|y|^{-1/4}v^{k}\parallel\parallel\overline{P}v^{k}\parallel\\
&\leq& M_{11}\parallel
v-Du_{yy}+\varrho_{1}u_{x}+(\varrho_{2}-\varepsilon\mu\partial_{y}
b_{11})u_{y}-\overline{P}\widetilde{u}^{k}\parallel\rightarrow 0.
\end{eqnarray*}
Then the following calculation shows that $\parallel
u-\widetilde{u}^{k}\parallel_{L^{2}(K)}\rightarrow 0$ for every
compact $K\subset\Omega$:
\begin{eqnarray*}
(-|y|^{-1/4}v^{k},\overline{P}v^{k})&=& \lim_{t\rightarrow\infty}
\int\int_{\Omega^{(0,t)}}[\frac{1}{2}|y|^{-1/4}\overline{B}_{x}+
\frac{1}{2}(|y|^{-1/4}\overline{C})_{y}-|y|^{-1/4}A](v^{k})^{2}\\
& &+
\int_{\partial\Omega^{(0,t)}}[-\frac{1}{2}|y|^{-1/4}\overline{C}n_{2}
-\frac{1}{2}|y|^{-1/4}\overline{B}n_{1}](v^{k})^{2}\\
&\geq&\lim_{t\rightarrow\infty}
\int\int_{\Omega^{(0,t)}}[|y|^{-1/4}(\frac{1}{4}\mu+
a_{11}-\frac{1}{2}+O(\varepsilon+\theta))](v^{k})^{2}\\
&\geq& M_{12}\parallel|y|^{-1/8}v^{k}\parallel_{K}^{2}.
\end{eqnarray*}
Therefore $u\equiv\widetilde{u}$ in $L^{2}(\Omega)$.\par
   Differentiating the equation $Pu=v$ with respect to
$\partial_{x}^{\alpha}$, $\alpha=1,\ldots,4$, and applying the
above procedure, shows that $\partial_{x}^{\alpha}u\in
L^{2}(\Omega)$, $\alpha=1,\ldots,4$.  We now have that $u\in
H^{4}(\Omega)$.\par
   To complete the case when $v\in
\widetilde{C}^{\infty}(\overline{\Omega})$, we must show that
$u(x_{0},y)=0$.  Since $B(x_{0},y)=0$, from the equation $Pu=v$ we
find that
\begin{equation*}
(Du_{yy}+Cu_{y}+Au)|_{(x_{0},y)}=v(x_{0},y)=0.
\end{equation*}
Furthermore since $u\in H^{4}(\Omega)$, $u\rightarrow 0$ as
$|y|\rightarrow\infty$.  Therefore, applying the maximum principle
to the above equation, we have $u(x_{0},y)=0$.\par
   We now consider the case when $v\in\widetilde{C}_{x}
^{\infty}(\overline{\Omega})$.  Let $h(y)\in
H^{\infty}(\mathbb{R})$ be the unique solution of the ODE:
\begin{equation*}
D(-x_{0},y)h''+C(-x_{0},y)h'+A(-x_{0},y)h=v(-x_{0},y).
\end{equation*}
Then as before, the parabolicity of $P$ guarantees the existence
of a unique solution to the Cauchy problem
\begin{equation*}
Pu=v\text{ }\text{ in }\text{ }\Omega,\text{ }\text{
}u(-x_{0},y)=h(y),
\end{equation*}
such that $u\in H^{\infty}(\Omega_{\tau})$ for every $\tau$.
Furthermore, $u_{x}(-x_{0},y)=0$ since
\begin{equation*}
Bu_{x}|_{(-x_{0},y)}=v(-x_{0},y)-(Du_{yy}+Cu_{y}+Au)|_{(-x_{0},y)}=0.
\end{equation*}
Moreover, the same methods used above can be used here to show
that $u\in H^{4}(\Omega)$.  Lastly, to show that
$u_{x}(x_{0},y)=0$, differentiate $Pu=v$ with respect to $x$ and
use that $B(x_{0},y)=0$ to obtain
\begin{equation*}
(D(u_{x})_{yy}+C(u_{x})_{y}+(A+B_{x})u_{x})|_{(x_{0},y)}=v_{x}(x_{0},y)
-(C_{x}u_{y}+A_{x}u)|_{(x_{0},y)}=0.
\end{equation*}
Since $u_{x}\rightarrow 0$ as $|y|\rightarrow\infty$, by the
maximum principle $u_{x}(x_{0},y)=0$.\qed
   With lemma 2.4 we are now in a position to establish the dual
inequalities.\medskip

\noindent\textbf{Proposition 2.1.}  \textit{There exist constants
$M_{1},M_{2}$ such that:}
\begin{equation*}
\parallel v\parallel\leq M_{1}\parallel
L_{\theta}^{*}v\parallel_{(-1,-2)}\text{ }\textit{ for all }\text{
}\text{ } v\in \widetilde{C}^{\infty}(\overline{\Omega}),
\end{equation*}
\begin{equation*}
\parallel v\parallel\leq M_{2}\parallel
L_{\theta}^{*}v\parallel_{(-1,-2)}\text{ }\textit{ for all }\text{
}\text{ } v\in \widetilde{C}_{x}^{\infty}(\overline{\Omega}).
\end{equation*}\medskip
\noindent\textit{Proof.}  We first consider the case when
$v\in\widetilde{C}^{\infty}(\overline{\Omega})$.  Let $u\in
C^{\infty}(\Omega)\cap H^{4}(\Omega)$ be the unique solution of
\begin{equation*}
Pu=v\text{ }\text{ in }\text{ }\Omega,\text{ }\text{
}u(-x_{0},y)=u_{x}(-x_{0},y)=0,\text{ }\text{ }u(x_{0},y)=0,
\end{equation*}
given by lemma 2.4.\par
   We now show that
\begin{equation*}
(Au+Bu_{x}+Cu_{y}+Du_{yy},L_{\theta}u)\geq
\end{equation*}
\begin{equation*}
C_{1}[\parallel u\parallel^{2}+\parallel Eu_{y}\parallel^{2}
+\theta(\parallel u_{x}\parallel^{2}+\parallel u_{xy}\parallel^{2}
+\parallel u_{yy}\parallel^{2}+\theta\parallel
u_{xyy}\parallel^{2})],
\end{equation*}\par
\noindent where $A,B,C,D,E$, and $C_{1}$ were given in lemma 2.3.
Let $\nu_{k}$ be given by (2.7) and define the sequence
$\{u^{k}\}_{k=1}^{\infty}$, where $u^{k}=\nu_{k}u$.  Then as in
the proof of lemma 2.4 we have:\bigskip

$i)$ $u^{k}\in C^{\infty}(\Omega)\cap H^{4}(\Omega)$,\bigskip

$ii)$ $u^{k}$ has bounded support and $u^{k}_{x}(-x_{0},y)=0$,
$u^{k}(x_{0},y)=0$,\bigskip

$iii)$ $\parallel u-u^{k}\parallel_{4}\rightarrow 0$ as
$k\rightarrow\infty$,\bigskip

$iv)$ $\parallel Eu_{y}-Eu^{k}_{y}\parallel\rightarrow 0$ as
$k\rightarrow\infty$.\bigskip

\noindent Let $\{u_{k}\}_{k=1}^{\infty}$ a $C^{\infty}$
approximation of $\{u^{k}\}_{k=1}^{\infty}$ such that:\bigskip

$i)$ $u_{k}\in C^{\infty}(\overline{\Omega})$,\bigskip

$ii)$ $u_{k}$ has bounded support and $(u_{k})_{x}(-x_{0},y)=0$,
$u_{k}(x_{0},y)=0$,\bigskip

$iii)$ $\parallel u^{k}-u_{k}\parallel_{4}\rightarrow 0$ as
$k\rightarrow\infty$,\bigskip

$iv)$ $\parallel Eu^{k}_{y}-E(u_{k})_{y}\parallel\rightarrow 0$ as
$k\rightarrow\infty$.\bigskip

\noindent  Then applying lemma 2.3 we have,
\begin{equation*}
(Au+Bu_{x}+Cu_{y}+Du_{yy},L_{\theta}u)\text{ }\text{ } \text{
}\text{ }\text{ }\text{ }\text{ }\text{ }\text{ } \text{ }\text{
}\text{ }\text{ }\text{ }\text{ }\text{ }\text{ }\text{ }\text{
}\text{ }\text{ }\text{ }\text{ }\text{ }\text{ }\text{ }\text{
}\text{ }\text{ }\text{ }\text{ }\text{ }\text{ }\text{ }\text{
}\text{ }\text{ }\text{ }\text{ }\text{ }\text{ }\text{ }\text{
}\text{ }\text{ }\text{ }\text{ }\text{ }\text{ }\text{ }\text{
}\text{ }\text{ }\text{ }\text{ }\text{ }\text{ }\text{ }\text{
}\text{ }\text{ }\text{ }\text{ }
\end{equation*}
\begin{equation*}
=\lim_{k\rightarrow\infty}(Au_{k}+B(u_{k})_{x}
+C(u_{k})_{y}+D(u_{k})_{yy},L_{\theta}u_{k})\text{ }\text{ }
\text{ }\text{ }\text{ }\text{ }\text{ }\text{ }\text{ }\text{ }
\text{ }\text{ }\text{ }\text{ }\text{ }\text{ }\text{ }\text{ }
\text{ }\text{ }\text{ }\text{ }\text{ }\text{ }\text{ }\text{
}\text{ }\text{ }\text{ }\text{ }\text{ }\text{ }\text{ }\text{
}\text{ }\text{ }\text{ }\text{ }\text{ }\text{ }\text{ }\text{
}\text{ }\text{ }\text{ }\text{ }\text{ }\text{ }\text{ }
\end{equation*}
\begin{equation*} \geq\lim_{k\rightarrow\infty}\!C_{1}[\parallel
u_{k}\!\parallel^{2} \!+\!\parallel
E(u_{k})_{y}\!\parallel^{2}\!+\theta(\parallel
(u_{k})_{x}\!\parallel^{2}\!+\!\parallel
(u_{k})_{xy}\!\parallel^{2}\! +\!\parallel
(u_{k})_{yy}\!\parallel^{2}\!+\theta\parallel
(u_{k})_{xyy}\!\parallel^{2})]
\end{equation*}
\begin{equation*}
=C_{1}[\parallel u\parallel^{2} +\parallel Eu_{y}\parallel^{2}
+\theta(\parallel u_{x}\parallel^{2}+\parallel u_{xy}\parallel^{2}
+\parallel u_{yy}\parallel^{2}+\theta\parallel
u_{xyy}\parallel^{2})].\text{ }\text{ }\text{ }\text{ }\text{ }
\text{ }\text{ }\text{ }\text{ }\text{ }\text{ }\text{ }\text{ }
\text{ }\text{ }\text{ }\text{ }\text{ }
\end{equation*}\par
   By the above estimate and definition of the negative norms, it
follows that
\begin{eqnarray*}
\parallel L_{\theta}^{*}v\parallel_{(-1,-2)}\parallel u\parallel
_{(1,2)}&\geq& (L_{\theta}^{*}v,u)\\ &=&(v,L_{\theta}u)\\
&=&(Au+Bu_{x}+Cu_{y}+Du_{yy},L_{\theta}u)\\
&\geq& C_{1}[\parallel u\parallel^{2} +\parallel
Eu_{y}\parallel^{2}\\
& & +\theta(\parallel u_{x}\parallel^{2}+\parallel
u_{xy}\parallel^{2} +\parallel u_{yy}\parallel^{2}+\theta\parallel
u_{xyy}\parallel^{2})].
\end{eqnarray*}
Furthermore, using Cauchy's inequality and the equation $Pu=v$, we
obtain
\begin{eqnarray*}
\parallel L_{\theta}^{*}v\parallel_{(-1,-2)}\!\!&\geq&\!\!
C_{1}'[\parallel u\parallel+\parallel Eu_{y}\parallel
+\sqrt{\theta}(\parallel u_{x}\parallel+\parallel
u_{xy}\parallel+\parallel u_{yy}\parallel+\sqrt{\theta}\parallel
u_{xyy}\parallel)]\\
&\geq& \!M_{1}^{-1}\parallel v\parallel,
\end{eqnarray*}
for some constants $C_{1}^{'},M_{1}>0$.  Moreover, similar
arguments may be used to treat the case when $v\in
\widetilde{C}_{x}^{\infty}(\overline{\Omega})$.\qed
   The existence of weak solutions to problems (2.3) and (2.4)
immediately follows from proposition 2.1 by a standard functional
analytic argument.  We include the proof here for
convenience.\medskip

\noindent\textbf{Theorem 2.1.}  \textit{For each $f\in
L^{2}(\Omega)$ there exists a weak solution $u\in
\widetilde{H}^{(1,2)}(\Omega)$,
$\widetilde{H}_{x}^{(1,2)}(\Omega)$ of} (2.3), (2.4)
\textit{respectively.}\medskip

\noindent\textit{Proof.}  We shall first treat problem (2.3).  Let
$W= L_{\theta}^{*}(\widetilde{C}^{\infty}(\overline{\Omega}))$ and
define the linear functional $F:W\rightarrow\mathbb{R}$ by
\begin{equation*}
F(L_{\theta}^{*}v)=(f,v).
\end{equation*}
Using proposition 2.1, the following calculation will show that
$F$ is bounded as a linear functional on the subspace $W$ of
$\widetilde{H}^{(-1,-2)}(\Omega)$,
\begin{eqnarray*}
|F(L_{\theta}^{*}v)|&=&|(f,v)|\\
&\leq&\parallel f\parallel\parallel v\parallel\\
&\leq&M_{1}\parallel f\parallel\parallel L_{\theta}^{*}v
\parallel_{(-1,-2)}.
\end{eqnarray*}
Use the Hahn-Banach theorem to extend $F$ from $W$ onto the whole
space $\widetilde{H}^{(-1,-2)}(\Omega)$.  It follows from the
Riesz representation theorem that there exists $u\in
\widetilde{H}^{(1,2)}(\Omega)$ such that
\begin{equation*}
F(w)=(u,w)\text{ }\text{ for all }\text{
}w\in\widetilde{H}^{(-1,-2)}(\Omega).
\end{equation*}
Thus, restricting $w$ to $W$ we have
\begin{equation*}
(u,L_{\theta}^{*}v)=F(L_{\theta}^{*}v)=(f,v)\text{ }\text{ for all
}\text{ } v\in\widetilde{C}^{\infty}(\overline{\Omega}).
\end{equation*}
The case of problem (2.4) may be treated in a similar manner.\qed
   We now prove the existence of weak solutions for the adjoint
problems (2.5) and (2.6).  The existence of solutions for these
problems will be needed in the next section, where they will aid
in proving higher regularity for solutions of (2.3).\par
   The formal adjoint of $L_{\theta}$ is given by
\begin{equation*}
L_{\theta}^{*}=-\theta\partial_{xxyy}+a_{11}\partial_{xx}
+a_{22}\partial_{yy}
\end{equation*}
\begin{equation*}
+(2\partial_{x}a_{11}-a_{1})\partial_{x}
+(2\partial_{y}a_{22}-a_{2})\partial_{y}+(\partial_{xx}a_{11}+
\partial_{yy}a_{22}-\partial_{x}a_{1}-\partial_{y}a_{2}+a).
\end{equation*}
All the coefficients of $L_{\theta}^{*}$, denoted
$a_{ij}^{*},a_{i}^{*},a^{*}$, have the same properties as the
coefficients of $L_{\theta}$, except
$a_{2}^{*}=2\partial_{y}a_{22}-a_{2}$.  This difference will not
allow us to directly apply the above procedure to obtain weak
solutions for (2.5) and (2.6).  However, if
\begin{equation*}
h(x,y)=e^{2\int_{0}^{y}\frac{a_{2}(x,t)}{a_{22}(x,t)}dt},
\end{equation*}
then by setting $v=hw$, the equation $L_{\theta}^{*}v=g$ becomes
$\overline{L}_{\theta}^{*}w=g/h$, where
\begin{equation*}
\overline{L}_{\theta}^{*}=-\theta\partial_{xxyy}-2\theta\frac{h_{y}}{h}
\partial_{xxy}-2\theta\frac{h_{x}}{h}\partial_{xyy}
\end{equation*}
\begin{equation*}
+(a^{*}_{11}-\theta
\frac{h_{yy}}{h})\partial_{xx}-4\theta\frac{h_{xy}}{h}\partial_{xy}
+(a^{*}_{22}-\theta\frac{h_{xx}}{h})\partial_{yy}
\end{equation*}
\begin{equation*}
+(a^{*}_{2}+2a^{*}_{22}\frac{h_{y}}{h}-2\theta\frac{h_{xxy}}{h})\partial_{y}
+(a^{*}_{1}+2a^{*}_{11}\frac{h_{x}}{h}-2\theta\frac{h_{xyy}}{h})\partial_{x}
\end{equation*}
\begin{equation*}
+(a^{*}_{11}\frac{h_{xx}}{h}+a^{*}_{22}\frac{h_{yy}}{h}+a^{*}_{1}\frac{h_{x}}{h}
+a^{*}_{2}\frac{h_{y}}{h}+a^{*}-\theta\frac{h_{xxyy}}{h}).
\end{equation*}
The special choice of $h$ guarantees that the coefficient of
$\partial_{y}$ in $\overline{L}_{\theta}^{*}$ is $3a_{2}+O(
\varepsilon+\theta)$, so that all the coefficients of
$\overline{L}_{\theta}^{*}$ have the same properties as the
coefficients of $\overline{L}_{\theta}$, where
$\overline{L}_{\theta}w=f/h$ is the equation obtained from
$L_{\theta}u=f$ by setting $u=hw$.  Therefore if $g\in
L^{2}(\Omega)$, the problems
\begin{equation*}
\overline{L}_{\theta}^{*}w=g/h\text{ }\text{ }\text{ in }\text{
}\text{ }\Omega,\text{ }\text{ }\text{ } w|_{\partial\Omega}=0,
\end{equation*}
\begin{equation*}
\overline{L}_{\theta}^{*}w=g/h\text{ }\text{ }\text{ in }\text{
}\text{ }\Omega,\text{ }\text{ }\text{ }
w_{x}|_{\partial\Omega}=0,
\end{equation*}
have weak solutions of the form $w=v/h$, where $v\in\widetilde{H}
^{(1,2)}(\Omega)$, $\widetilde{H}_{x}^{(1,2)}(\Omega)$
respectively.  We then obtain\medskip

\noindent\textbf{Corollary 2.1.}  \textit{For each $g\in
L^{2}(\Omega)$ there exists a weak solution
$v\in\widetilde{H}^{(1,2)}(\Omega)$,}
$\widetilde{H}_{x}^{(1,2)}(\Omega)$ \textit{of} (2.5), (2.6)
\textit{respectively}.

\bigskip\bigskip
\noindent\textbf{3.  Linear Regularity}
\setcounter{equation}{0}
\setcounter{section}{3}
\bigskip\bigskip

\noindent The purpose of this section is to establish the
regularity, in $X$, of weak solutions to problem (2.3) for a
particular choice of the right-hand side, $f$.  This shall be
accomplished by establishing the uniqueness of weak solutions to
problems (2.3) and (2.4) in $L^{2}(\Omega)$, and then applying a
boot-strap argument.\par
   In order to obtain the uniqueness of weak solutions, we will
utilize the notion of a strong solution, in particular, for first
order systems.  The definition of a strong solution will be given
below.  We first introduce the notation and terminology that will
be used for first order systems.  Consider a boundary value
problem
\begin{equation}
SU=A_{1}U_{x}+A_{2}U_{y}+A_{3}U=F \text{ }\text{ }\text{ in
}\text{ }\text{ } \Omega, \text{ }\text{ } U|_{\partial\Omega}\in
N,
\end{equation}
where $A_{1},A_{2},A_{3}$ are $n\times n$ matrices, $U$ and $F$
are $n$-vectors, and $N$ is a linear subspace of the space of
$n$-vector valued functions restricted to $\partial\Omega$.  The
corresponding adjoint problem is given by
\begin{equation*}
S^{*}V=-A^{*}_{1}V_{x}-A^{*}_{2}V_{y}+(A^{*}_{3}-\partial_{x}A^{*}_{1}
-\partial_{y}A^{*}_{2})V=G \text{ }\text{ }\text{ in }\text{
}\text{ }\Omega,\text{ }\text{ } V|_{\partial\Omega}\in N^{*},
\end{equation*}
where $A^{*}_{i}$ denotes the transpose of $A_{i}$, and $N^{*}$ is
the orthogonal complement of $\bigtriangleup N$, where
$\bigtriangleup$ is the matrix defined on $\partial\Omega$ by
$A_{1}n_{1}+A_{2}n_{2}$, and $(n_{1},n_{2})$ is the unit outward
normal to $\partial\Omega$.\par
   Let $F\in L^{2}(\Omega)$.  The notion of a weak solution to problem
(3.1) is similar to the definition given in section $\S 2$ for
single equations.  That is, $U\in L^{2}(\Omega)$ is said to be a
weak solution of (3.1) whenever
\begin{equation*}
(S^{*}V,U)=(V,F),
\end{equation*}
for every $V\in C^{\infty}(\overline{\Omega})$ with bounded
support and such that $V|_{\partial\Omega}\in N^{*}$. We now give
the definition of a strong solution.\medskip

\noindent\textbf{Definition 3.1.}  \textit{$U\in L^{2}(\Omega)$ is
a strong solution of} \text{(3.1)}\textit{ if there exists a
sequence $\{U_{k}\}^{\infty}_{k=1}$, such that $U_{k}\in
C^{\infty}(\overline{\Omega})$ with bounded support,
$U_{k}|_{\partial\Omega}\in N$, and}
\begin{equation*}
\parallel U_{k}-U\parallel\rightarrow 0,\text{ }\text{ }
\parallel SU_{k}-F\parallel\rightarrow 0,\text{ }\text{ }\textit{ as }\text{ }\text{ }
k\rightarrow\infty.
\end{equation*}
\noindent Clearly, a strong solution is a weak solution. Moreover,
using techniques developed by Friedrichs [2] and Lax-Phillips
[14], G. Peyser [21] has obtained the following converse
statement.\medskip

\noindent\textbf{Theorem 3.1 (Identity of Weak and Strong
Solutions).} \textit{Let the following conditions on the operator
$S$ and the boundary space $N$ be satisfied:}

$i)$ \textit{the matrix $\bigtriangleup$ is of constant rank in a
neighborhood of the boundary},

$ii)$ $N$ \textit{is of constant dimension at each point of the
boundary},

$iii)$ $N$ \textit{contains the nullspace of} $\bigtriangleup$.

\noindent\textit{Then a weak solution $U\in L^{2}(\Omega)$ of}
(3.1) \textit{is also a strong solution.}\medskip

\noindent Note that for our particular domain,
$\bigtriangleup=A_{1}n_{1}$, so that condition $(i)$ is equivalent
to $A_{1}$ having constant rank in a neighborhood of
$\partial\Omega$.\par
   With the aim of applying theorem 3.1, we shall transform
problems (2.3), (2.4), (2.5), and (2.6) into the setting of first
order systems.  Let $f,g\in L^{2}(\Omega)$ be the right-hand sides
of (2.3), (2.4) and (2.5), (2.6) respectively, and define $A_{1}$,
$\widetilde{A}_{1}$, $A_{2}$, $\widetilde{A}_{2}$, $A_{3}$,
$\widetilde{A}_{3}$, $F$, and $G$ by
\begin{equation*}
A_{1}=\widetilde{A}_{1}=\left(\begin{array}{ccccc}
-\theta & 0 & a_{11} & 0 & 0 \\
0 & 0 & 0 & 0 & 1 \\
0 & 0 & 0 & 0 & 0 \\
0 & 0 & 0 & 0 & 0 \\
0 & 1 & 0 & 0 & 0
\end{array}\right),
\text{  }\text{  }
A_{2}=\widetilde{A}_{2}=\left(\begin{array}{ccccc}
0 & 0 & 0 & a_{22} & 0 \\
0 & 0 & 0 & 0 & 0 \\
0 & 0 & 0 & 0 & 1 \\
0 & 0 & 0 & 1 & 0 \\
0 & 0 & 0 & 0 & 0
\end{array}\right),
\end{equation*}
\begin{equation*}
A_{3}=\left(\begin{array}{ccccc}
0 & 0 & a_{1} & a_{2} & a \\
0 & 0 & -1 & 0 & 0 \\
0 & 0 & 0 & -1 & 0 \\
0 & -1 & 0 & 0 & 0 \\
-1 & 0 & 0 & 0 & 0
\end{array}\right),
\text{  }\text{  }
\widetilde{A}_{3}=\left(\begin{array}{ccccc}
0 & 0 & a^{*}_{1} & a^{*}_{2} & a^{*} \\
0 & 0 & -1 & 0 & 0 \\
0 & 0 & 0 & -1 & 0 \\
0 & -1 & 0 & 0 & 0 \\
-1 & 0 & 0 & 0 & 0
\end{array}\right),
\end{equation*}
\begin{equation*}
F=\left(\begin{array}{c} f \\ 0 \\ 0 \\ 0 \\ 0
\end{array}\right),
\text{  }\text{  }\text{  }\text{  }
G=\left(\begin{array}{c} g \\
0 \\ 0 \\ 0 \\ 0
\end{array}\right).
\end{equation*}
Define boundary spaces $N_{1}$ and $N_{2}$ by
\begin{equation*}
N_{1}=\{(u_{1},\ldots,u_{5})|_{\partial\Omega}\mid u_{5}
|_{\partial\Omega}=0\},
\end{equation*}
\begin{equation*}
N_{2}=\{(u_{1},\ldots,u_{5})|_{\partial\Omega}\mid (-\theta
u_{1}+a_{11}u_{3})|_{\partial\Omega}=0\}.
\end{equation*}
Furthermore, define boundary value problems
\begin{equation}
S_{\theta}U=A_{1}U_{x}+A_{2}U_{y}+A_{3}U=F\text{ }\text{ }\text{
in }\text{ }\text{ } \Omega,\text{ }\text{ }U|_{\partial\Omega}\in
N_{1},
\end{equation}
\begin{equation}
S_{\theta}U=F\text{ }\text{ }\text{ in } \text{ }\text{
}\Omega,\text{ }\text{ }U|_{\partial\Omega}\in N_{2},
\end{equation}
\begin{equation}
\widetilde{S}_{\theta}V=\widetilde{A}_{1}V_{x}+\widetilde{A}_{2}V_{y}
+\widetilde{A}_{3}V=G\text{ }\text{ }\text{ in }\text{ }\text{ }
\Omega,\text{ }\text{ }V|_{\partial\Omega}\in N_{1},
\end{equation}
\begin{equation}
\widetilde{S}_{\theta}V=G\text{ }\text{ }\text{ in }\text{ }\text{
} \Omega,\text{ }\text{ }V|_{\partial\Omega}\in N_{2}.
\end{equation}
We now show that the weak solutions of (2.3), (2.4), (2.5), and
(2.6) given by theorem 2.1 and corollary 2.1 are also weak
solutions of (3.2), (3.3), (3.4), and (3.5) respectively.\medskip

\noindent\textbf{Lemma 3.1.}  \textit{Let
$u\in\widetilde{H}^{(1,2)}(\Omega)$,
$\widetilde{H}^{(1,2)}_{x}(\Omega)$ be a weak solution of} (2.3),
(2.4) \textit{respectively}, \textit{then
$U=(u_{xyy},u_{yy},u_{x},u_{y},u)\in L^{2}(\Omega)$ is a weak
solution of} (3.2), (3.3) \textit{respectively}.
\textit{Similarly, if $v\in\widetilde{H}^{(1,2)}(\Omega)$,
$\widetilde{H}^{(1,2)}_{x}(\Omega)$ is a weak solution of} (2.5),
(2.6) \textit{respectively}, \textit{then
$V=(v_{xyy},v_{yy},v_{x},v_{y},v)\in L^{2}(\Omega)$ is a weak
solution of} (3.4), (3.5) \textit{respectively}.\medskip

\noindent\textit{Proof.}  Let $u\in\widetilde{H}^{(1,2)}(\Omega)$
be a weak solution of problem (2.3).  We will show that
\begin{equation}
\int\int_{\Omega}U^{*}S^{*}_{\theta}V=\int\int_{\Omega}F^{*}V
\end{equation}
for all $V\in C^{\infty}(\overline{\Omega})$ with bounded support
such that $V|_{\partial\Omega}\in N^{*}_{1}$, where
\begin{equation*}
N^{*}_{1}=\{(v_{1},\ldots,v_{5})|_{\partial\Omega}\mid
v_{1}|_{\partial\Omega}=v_{5}|_{\partial\Omega}=0\}.
\end{equation*}
A calculation shows that
\begin{eqnarray}
\int\int_{\Omega}U^{*}S^{*}_{\theta}V & =
&\int\int_{\Omega}(\theta
u_{xyy}-a_{11}u_{x})\partial_{x}v_{1}-a_{22}u_{y}\partial_{y}v_{1}
-(u\partial_{x}v_{2}+u_{x}v_{2})\nonumber\\
& &+[(a_{1}-\partial_{x}a_{11})u_{x}+(a_{2}-\partial_{y}a_{22})u_{y}+au]v_{1}\\
& & -(u\partial_{y}v_{3}+u_{y}v_{3})-(u_{y}\partial_{y}v_{4}
+u_{yy}v_{4})-(u_{yy}\partial_{x}v_{5}+u_{xyy}v_{5}).\nonumber
\end{eqnarray}
Since $V|_{\partial\Omega}\in N^{*}_{1}$ and
$u\in\widetilde{H}^{(1,2)}(\Omega)$ is a weak solution of (2.3),
we can integrate by parts to obtain
\begin{equation*}
\int\int_{\Omega}U^{*}S^{*}_{\theta}V=\int\int_{\Omega}uL^{*}_{\theta}v_{1}=
\int\int_{\Omega}fv_{1}=\int\int_{\Omega}F^{*}V,
\end{equation*}
showing that $U$ is a weak solution of (3.2).\par
   Let $u\in\widetilde{H}_{x}^{(1,2)}(\Omega)$ be a weak solution
of (2.4).  We now show that (3.6) holds for all $V\in
C^{\infty}(\overline{\Omega})$ with bounded support such that
$V|_{\partial\Omega}\in N^{*}_{2}$, where
\begin{equation*}
N^{*}_{2}=\{(v_{1},\ldots,v_{5})|_{\partial\Omega}\mid
v_{2}|_{\partial\Omega}=v_{5}|_{\partial\Omega}=0\}.
\end{equation*}
From (3.7) it follows that,
\begin{eqnarray}
\int\int_{\Omega}U^{*}S^{*}_{\theta}V & =
&\int\int_{\Omega}(\theta
u_{xyy}-a_{11}u_{x})\partial_{x}v_{1}-a_{22}u_{y}\partial_{y}v_{1}\\
& &+[(a_{1}-\partial_{x}a_{11})u_{x}+(a_{2}-
\partial_{y}a_{22})u_{y}+au]v_{1}.\nonumber
\end{eqnarray}
In order to integrate by parts we construct an approximating
sequence $\{v_{1}^{k}\}_{k=1}^{\infty}$ for $v_{1}$, such that
$v_{1}^{k}\in\widetilde{C}_{x}^{\infty}(\overline{\Omega})$ and
\begin{equation*}
\parallel v_{1}^{k}-v_{1}\parallel+\parallel\partial_{x}
v_{1}^{k}-v_{1}\parallel\rightarrow 0\text{ }\text{ }\text{
}\text{ as } \text{ }\text{ }\text{ }k\rightarrow\infty.
\end{equation*}
Take a sequence
$\{v_{k}\}_{k=1}^{\infty}\subset\widetilde{C}^{\infty}(\overline{\Omega})$
with the property that $\parallel
v_{k}-\partial_{x}v_{1}\parallel\rightarrow 0$ as
$k\rightarrow\infty$, and define
\begin{equation*}
v^{k}_{1}=\int_{-x_{0}}^{x}v_{k}(t,y)dt+v_{1}(-x_{0},y).
\end{equation*}
Then since
\begin{eqnarray*}
(v_{1}^{k}-v_{1})^{2}&=&(\int_{-x_{0}}^{x}\partial_{t}
(v_{1}^{k}(t,y)-v_{1}(t,y))dt)^{2}\\
&\leq&2x_{0}\int_{-x_{0}}^{x}(\partial_{t}
v_{1}^{k}(t,y)-\partial_{t}v_{1}(t,y))^{2}dt,
\end{eqnarray*}
we have
\begin{eqnarray*}
\int\int_{\Omega}(v_{1}^{k}-v_{1})^{2}&\leq&
4x_{0}^{2}\int\int_{\Omega}(\partial_{x}v_{1}^{k}-\partial_{x}v_{1})^{2}\\
&=&4x_{0}^{2}\int\int_{\Omega}(v_{k}-\partial_{x}v_{1})^{2}\rightarrow
0,
\end{eqnarray*}
so that $v_{1}^{k}$ satisfies the desired properties. Therefore,
recalling that $a_{1}|_{\partial\Omega}=$\linebreak
$\partial_{x}a_{11}|_{\partial\Omega}=0$ by $(ii)$ of lemma 2.2,
and using the fact that $u$ is a weak solution of (2.4), we can
integrate by parts in (3.8) to obtain
\begin{eqnarray*}
\int\int_{\Omega}U^{*}S^{*}_{\theta}V & = &
\lim_{k\rightarrow\infty}\int\int_{\Omega}
uL^{*}_{\theta}v_{1}^{k}\\
& = &\lim_{k\rightarrow\infty}\int\int_{\Omega}
fv_{1}^{k}=\int\int_{\Omega}fv_{1}=\int\int_{\Omega}F^{*}V,
\end{eqnarray*}
showing that $U$ is a weak solution of (3.3).  Similar arguments
show that if\linebreak $v\in\widetilde{H}^{(1,2)}(\Omega)$,
$\widetilde{H}^{(1,2)}_{x}(\Omega)$ is a weak solution of (2.5),
(2.6) respectively, then\linebreak
$V=(v_{xyy},v_{yy},v_{x},v_{y},v)\in L^{2}(\Omega)$ is a weak
solution of (3.4), (3.5) respectively.\qed
   Now that the weak solutions of the previous section have
been translated into the setting of first order systems, theorem
3.1 is applicable.  As a result, we obtain\medskip

\noindent\textbf{Proposition 3.1.}  \textit{The weak solutions of
problems} (2.3) \textit{and} (2.4), \textit{given by theorem} 2.1,
\textit{are unique in $L^{2}(\Omega)$.}\medskip

\noindent\textit{Proof.}  Let $u\in\widetilde{H}^{(1,2)}(\Omega)$
be a weak solution of problem (2.3) with $f=0$, then
\begin{equation}
(L^{*}_{\theta}w,u)=0 \text{ }\text{ }\text{ for all }\text{
}\text{ } w\in\widetilde{C}^{\infty}(\overline{\Omega}).
\end{equation}
We will show that $u=0$ in $L^{2}(\Omega)$.\par
   Let $v\in\widetilde{H}^{(1,2)}(\Omega)$ be the weak solution of
(2.5) with $g=u$.  Then by lemma 3.1
$V=(v_{xyy},v_{yy},v_{x},v_{y},v)$ is a weak solution of (3.4). We
now show that the conditions of theorem 3.1 are satisfied for
problem (3.4).  Condition $(ii)$ is immediately satisfied, and
since $a_{11}^{*}\leq-\theta$ in a neighborhood of
$\partial\Omega$, condition $(i)$ is satisfied with
$\triangle=\pm\widetilde{A}_{1}$ having the constant rank of 3.
Furthermore, the nullspace of $\triangle$ is given by
\begin{equation*}
\{(v_{1},\ldots,v_{5})|_{\partial\Omega}\mid(-\theta
v_{1}+a_{11}v_{3})|_{\partial\Omega}=v_{2}|_{\partial\Omega}=v_{5}
|_{\partial\Omega}=0\},
\end{equation*}
which is contained in $N_{1}$ so that condition $(iii)$ is
satisfied.  Therefore, we can apply theorem 3.1 to obtain an
approximating sequence $\{V_{k}\}_{k=1}^{\infty}$ for $V$, such
that $V_{k}\in C^{\infty}(\overline{\Omega})$ with bounded
support, $V_{k}|_{\partial\Omega}\in N_{1}$, and
\begin{equation}
\parallel V_{k}-V\parallel\rightarrow 0,\text{  }\text{ }\text{ }
\parallel \widetilde{S}_{\theta}V_{k}-G\parallel\rightarrow
0\text{ }\text{ }\text{  as }\text{ }\text{ } k\rightarrow\infty.
\end{equation}\par
   From (3.10) it follows that
\begin{equation*}
\parallel
v_{k}^{1}-v_{xyy}\parallel\rightarrow 0,\text{ }\text{ }\text{ }
\parallel v_{k}^{2}-v_{yy}\parallel\rightarrow 0, \text{ }\text{ }\text{ }\parallel
v_{k}^{3}-v_{x}\parallel\rightarrow 0,
\end{equation*}
\begin{equation*}
\parallel
v_{k}^{4}-v_{y}\parallel\rightarrow 0,\text{ }\text{ }\text{ }
\parallel v_{k}^{5}-v\parallel\rightarrow 0,
\end{equation*}
and
\begin{equation*}
\parallel(-\theta\partial_{x}v_{k}^{1}+a_{11}^{*}\partial_{x}
v_{k}^{3}+a_{22}^{*}\partial_{y}v_{k}^{4}+a_{1}^{*}v_{k}^{3}
+a_{2}^{*}v_{k}^{4}+a^{*}v_{k}^{5})-u\parallel\rightarrow 0.
\end{equation*}
Hence,
\begin{equation*}
(u,u)=\lim_{k\rightarrow\infty}\int\!\!\int_{\Omega}[
-\theta\partial_{x}v_{k}^{1}+a_{11}^{*}\partial_{x}v_{k}^{3}
+a_{22}^{*}\partial_{y}v_{k}^{4}+a_{1}^{*}v_{k}^{3}+a_{2}^{*}v_{k}^{4}
+a^{*}v_{k}^{5}]u\text{ }\text{ }\text{ }\text{ }\text{ }\text{ }
\text{ }\text{ }\text{ }\text{ }\text{ }\text{ }\text{ }\text{ }
\text{ }\text{ }\text{ }\text{ }\text{ }\text{ }
\end{equation*}
\begin{equation*}
\text{ }\text{ }\text{ }
=\lim_{k\rightarrow\infty}\int\!\!\int_{\Omega}(\theta
v_{k}^{1}\!-\!a_{11}^{*}v_{k}^{3})u_{x}\!-\!a_{22}^{*}v_{k}^{4}u_{y}\!+\!
[(a_{1}^{*}\!-\!\partial_{x}a_{11}^{*})v_{k}^{3}\!+\!(a_{2}^{*}\!-\!\partial_{y}a_{22}^{*})
v_{k}^{4}\!+\!a^{*}v_{k}^{5}]u
\end{equation*}
\begin{equation*}
\text{ }\text{ }\text{ }\text{ }\text{ }\text{ }
=\int\!\!\int_{\Omega}(\theta
v_{xyy}-a_{11}^{*}v_{x})u_{x}-a_{22}^{*}v_{y}u_{y}+
[(a_{1}^{*}-\partial_{x}a_{11}^{*})v_{x}+(a_{2}^{*}-\partial_{y}a_{22}^{*})
v_{y}+a^{*}v]u.
\end{equation*}
Let $\{v_{n}\}_{n=1}^{\infty}\subset
\widetilde{C}^{\infty}(\overline{\Omega})$ be an approximating
sequence for $v$ in $\widetilde{H}^{(1,2)}(\Omega)$.  Then
integrating by parts and using (3.9), we obtain
\begin{equation*}
(u,u)=\lim_{n\rightarrow\infty}\int\int_{\Omega}(L_{\theta}^{*}v_{n})u=0.
\end{equation*}
Similar arguments hold for problem (2.4).\qed
   Having established the uniqueness of weak solutions, we are now
ready to apply a boot-strap procedure to obtain higher regularity
for problem (2.3) in the $x$-direction.\medskip

\noindent\textbf{Theorem 3.2.}  \textit{Let $u$ and $f$ be as in
problem} (2.3).  \textit{Let $s\leq r-4$ and $f\in H^{s}(\Omega)$
be such that $\partial_{x}^{\alpha}f|_{\partial\Omega}=0$ for
$\alpha\leq s-1$. If $\varepsilon=\varepsilon(s)$ is sufficiently
small, then for all $\alpha\leq s$,
$\partial_{x}^{\alpha}u\in\widetilde{H}^{(1,2)}(\Omega)$ when
$\alpha$ is even, and $\partial_{x}^{\alpha}u\in
\widetilde{H}_{x}^{(1,2)}(\Omega)$ when $\alpha$ is odd.}\medskip

\noindent\textit{Proof.}  The case $s=0$ is given by theorem 2.1.
Consider the case $s=1$.  Let $w=u_{x}$ and formally differentiate
the equation $L_{\theta}u=f$ with respect to $x$:
\begin{eqnarray*}
L_{1}w&:=&-\theta
w_{xxyy}+a_{11}w_{xx}+a_{22}w_{yy}+(a_{1}+\partial_{x}a_{11})w_{x}
+a_{2}w_{y}+(a+\partial_{x}a_{1})w\\
&=&f_{x}-u_{yy}\partial_{x}a_{22}-u_{y}\partial_{x}a_{2}-u\partial_{x}a
:=f_{1}.
\end{eqnarray*}
Observe that since
$\partial_{x}a_{11},\partial_{x}a_{1}=O(\varepsilon)$ and both
vanish outside $X$, the operator $L_{1}$ has the same existence
and uniqueness properties as $L_{\theta}$.  Furthermore, by
restricting $L_{\theta}u=f$ to the boundary of $\Omega$ and using
$u|_{\partial\Omega}=a_{1}|_{\partial\Omega}=0$, we obtain the
following ODE
\begin{equation}
(-\theta u_{xxyy}+a_{11}u_{xx})|_{\partial\Omega}=0,
\end{equation}
for which the only solution in $L^{2}(\partial\Omega)$ is
$u_{xx}|_{\partial\Omega}=0$.  Therefore, in the regular case
$w=u_{x}$ satisfies problem (2.4) with $L_{\theta}$ and $f$
replaced by $L_{1}$ and $f_{1}$.\par
   Let $u\in\widetilde{H}^{(1,2)}(\Omega)$ be the weak solution of
problem (2.3).  We now show that $u_{x}\in L^{2}(\Omega)$ is a
weak solution of (2.4) with $L_{\theta}$ and $f$ replaced by
$L_{1}$ and $f_{1}\in L^{2}(\Omega)$; we denote this problem by
$\text{(2.4)}_{1}$.  Let
$v\in\widetilde{C}^{\infty}_{x}(\overline{\Omega})$, then
\begin{eqnarray*}
(u_{x},L_{1}^{*}v)\!\!\!\!&=&\!\!\!\!-(u,(L_{1}^{*}v)_{x})=-(u,L^{*}(v_{x}))
+(u,L^{*}(v_{x})-(L_{1}^{*}v)_{x})\\
&=&\!\!\!\!-(f,v_{x})+(u,-v_{yy}\partial_{x}a_{22}\!+v_{y}[
\partial_{x}a_{2}-2\partial_{xy}a_{22}]
+v[-\partial_{x}a\!-
\partial_{xyy}a_{22}+\partial_{xy}a_{2}])\\
&=&\!\!\!\!(f_{x},v)+(-u_{yy}\partial_{x}a_{22}-u_{y}\partial_{x}a_{2}
-u\partial_{x}a,v)=(f_{1},v).
\end{eqnarray*}
Therefore $u_{x}$ is a weak solution of $\text{(2.4)}_{1}$, and by
the uniqueness result, proposition 3.1, $u_{x}$ must coincide with
the solution in $\widetilde{H}_{x}^{(1,2)}(\Omega)$ given by
theorem 2.1.  Hence
$u_{x}\in\widetilde{H}_{x}^{(1,2)}(\Omega)$.\par
   We now consider the case $s=2$.  Let $w=u_{xx}$ and formally
differentiate the equation $L_{1}u_{x}=f_{1}$ with respect to $x$:
\begin{eqnarray*}
L_{2}w&:=&-\theta
w_{xxyy}+a_{11}w_{xx}+a_{22}w_{yy}\\
& &+(a_{1}+2\partial_{x}a_{11})w_{x}
+a_{2}w_{y}+(a+2\partial_{x}a_{1}+\partial_{xx}a_{11})w\\
&=&\partial_{x}f_{1}-u_{xyy}\partial_{x}a_{22}-u_{xy}\partial_{x}a_{2}-u_{x}
(\partial_{x}a+\partial_{xx}a_{1}) :=f_{2}.
\end{eqnarray*}
Again, since
$\partial_{x}a_{11},\partial_{xx}a_{11},\partial_{x}a_{1}=O(\varepsilon)$
and all three vanish outside $X$, the operator $L_{2}$ has the
same existence and uniqueness properties as $L_{\theta}$, provided
that $\varepsilon$ is sufficiently small.  Also, when $u$ is
regular $u_{xx}|_{\partial\Omega}=0$ from (3.11).  Thus in the
regular case $w=u_{xx}$ satisfies (2.3) with $L_{\theta}$ and $f$
replaced by $L_{2}$ and $f_{2}\in L^{2}(\Omega)$; we denote this
problem by $\text{(2.3)}_{2}$.\par
   Let $u\in\widetilde{H}^{(1,2)}(\Omega)$ be the weak solution of
(2.3), then we know that
$u_{x}\in\widetilde{H}_{x}^{(1,2)}(\Omega)$.  We now show that
$u_{xx}\in L^{2}(\Omega)$ is a weak solution of
$\text{(2.3)}_{2}$. Note that $L_{\theta}u=f$ in $L^{2}(\Omega)$
and let $v\in\widetilde{C}^{\infty}(\overline{\Omega})$, then a
calculation produces
\begin{eqnarray*}
(u_{xx},L_{2}^{*}v)&=&(u_{xxyy},-\theta
v_{xx})+(u_{xx},(a_{11}v)_{xx})+(u_{yy},(a_{22}v)_{xx})+(u_{y},(a_{2}v)_{xx})\\
& &+(u_{x},[(a_{1}+2\partial_{x}a_{11})v]_{xx})
+(u,[(a+2\partial_{x}a_{1}+\partial_{xx}a_{11})v]_{xx})\\
&=&(L_{\theta}u,v_{xx})+(f_{2}-f_{xx},v)=(f,v_{xx})+(f_{2}-f_{xx},v)
=(f_{2},v).
\end{eqnarray*}
By the uniqueness of weak solutions for problem
$\text{(2.3)}_{2}$, $u_{xx}$ must coincide with the solution in
$\widetilde{H}^{(1,2)}(\Omega)$.  Thus
$u_{xx}\in\widetilde{H}^{(1,2)}(\Omega)$.\par
   To obtain the regularity of higher order derivatives, we observe that
the above procedure applied to $L_{\theta}u=f$ holds for
$L_{2}u_{xx}=f_{2}$, since for $\alpha\geq 1$
\begin{equation*}
\partial_{x}^{\alpha}a_{11}|_{\partial\Omega}=
\partial_{x}^{\alpha}a_{22}|_{\partial\Omega}=
\partial_{x}^{\alpha}a_{i}|_{\partial\Omega}=
\partial_{x}^{\alpha}a|_{\partial\Omega}=0,
\end{equation*}
so that $f_{2}|_{\partial\Omega}=0$.  Therefore $u_{xxx}\in
\widetilde{H}_{x}^{(1,2)}(\Omega)$ and
$u_{xxxx}\in\widetilde{H}^{(1,2)}(\Omega)$.  Furthermore, we can
continue this process until $f$ and the coefficients of
$L_{\theta}$ run out of derivatives, as long as $\varepsilon$ is
chosen sufficiently small depending on the size of $s$.\qed
   We now prove regularity in the $y$-direction for the weak
solution of problem (2.3).  The following standard lemma
concerning difference quotients will be needed.\medskip

\noindent\textbf{Lemma 3.2.}  \textit{Let $w\in L^{2}(\Omega)$
have bounded support, and define}
\begin{equation*}
w^{h}=\frac{1}{h}(w(x,y+h)-w(x,y)).
\end{equation*}
\textit{If $\parallel w^{h}\parallel\leq M$ where $M$ is
independent of $h$, then $w\in H^{(0,1)}(\Gamma)$ for any compact
$\Gamma\subset\Omega$. Furthermore, if $w\in H^{(0,1)}(\Omega)$
then $\parallel w^{h}\parallel\leq M\parallel
w_{y}\parallel$.}\medskip

\noindent\textbf{Theorem 3.3.}  \textit{Let the hypotheses of
theorem} 3.2 \textit{hold, then $u\in H^{s}(X)$.}\medskip

\noindent\textit{Proof.}  From theorem 3.2 we know that
$\partial_{x}^{\alpha}u\in H^{(1,2)}(\Omega)$ for $0\leq\alpha\leq
s$.  Therefore the following equality holds in $L^{2}(\Omega)$,
\begin{equation}
\widetilde{L}u_{yy}:=-\theta u_{xxyy}+a_{22}u_{yy}=
f-a_{11}u_{xx}-a_{1}u_{x}-a_{2}u_{y}-au:=\widetilde{f}.
\end{equation}
Since $|a_{2}|=O(|y|)$ as $|y|\rightarrow\infty$, we do not
necessarily know that $\widetilde{f}\in H^{(0,1)}(\Omega)$;
however, we do have $\widetilde{f}\in H^{(0,1)}(\Gamma)$ for any
compact $\Gamma\subset\Omega$.  Fix a constant $k>y_{0}$ and set
$w=\nu_{k}u_{yy}$, where $\nu_{k}$ is given by (2.7).  Then
\begin{equation}
\widetilde{L}w^{h}=(\nu_{k}\widetilde{f})^{h}-\nu_{k}(y+h)u_{yy}(x,y+h)a_{22}^{h}.
\end{equation}
Since $u\in\widetilde{H}^{(1,2)}(\Omega)$, by multiplying (3.13)
on both sides by $w^{h}$ and integrating by parts, we obtain
\begin{equation*}
\parallel w^{h}\parallel+\parallel w_{x}^{h}\parallel
\leq M_{1}(\parallel (\nu_{k}\widetilde{f})^{h}\parallel+1),
\end{equation*}
for some $M_{1}$ independent of $h$.  By lemma 3.2
\begin{equation*}
\parallel w^{h}\parallel+\parallel w_{x}^{h}\parallel
\leq M_{2}(\parallel\nu_{k}\widetilde{f}\parallel_{(0,1)}+1),
\end{equation*}
independent of $h$.  Therefore $w_{y}$, $w_{xy}\in L^{2}(X)$,
which implies that $\partial_{y}^{3}u$,\linebreak
$\partial_{x}\partial_{y}^{3}u\in L^{2}(X)$. Furthermore, by
differentiating $L_{\theta}u=f$ with respect to $x$,\linebreak
$\alpha=1,\ldots,s-3$ times, the same procedure yields
$\partial_{x}^{\alpha}\partial_{y}^{3}u\in L^{2}(X)$.\par
   Proceeding by induction on $l$, assume that $\partial_{x}
^{\alpha}\partial_{y}^{\beta}u\in L^{2}(X)$, $\alpha\leq s-\beta$,
$\beta\leq l$, and $3\leq l<s$.  Differentiate (3.12) with respect
to $y$, $l-2$ times:
\begin{equation}
\widetilde{L}\partial_{y}^{l}u=\partial_{y}^{l-2}\widetilde{f}
-\sum_{i=0}^{l-3}\partial_{y}^{i}(\partial_{y}a_{22}
\partial_{y}^{l-3-i}u_{yy}).
\end{equation}
Note that this equation holds in $L^{2}(\Omega)$, and that the
right-hand side is in $H^{(0,1)}(\Gamma)$ for any compact
$\Gamma\subset\Omega$.  Applying the method above yields
$\partial_{y}^{l+1}u,\partial_{x}\partial_{y}^{l+1}u\in L^{2}(X)$.
Moreover, differentiating (3.14) with respect to $x$,
$\alpha=1,\ldots,s-(l+1)$ times, and applying the same procedure,
yields $\partial_{x}^{\alpha}\partial_{y}^{l+1}u\in L^{2}(X)$. The
desired conclusion now follows by induction.\qed

\bigskip\medskip
\noindent\textbf{4.  The Moser Estimate}
\setcounter{equation}{0}
\setcounter{section}{4}
\bigskip\bigskip

\noindent Having established the existence of regular solutions to
a small perturbation of the linearized equation for (1.5), we
intend to apply a Nash-Moser type iteration procedure in the
following section, to obtain a smooth solution of (1.5) in a
subdomain of $X$ which contains the origin.  In the current
section, we shall make preparations for the Nash-Moser procedure
by establishing a certain a priori estimate. This estimate,
referred to as the Moser estimate, will establish the dependence
of the solution $u$ of (2.3), on the coefficients of $L_{\theta}$
as well as on the right-hand side, $f$.  The Moser estimate that
we seek has the form
\begin{equation}
\parallel u\parallel_{H^{s}(X)}\leq C_{s}(\parallel f\parallel
_{H^{s}(X)}+\Lambda_{s+s_{0}}\parallel f\parallel_{H^{2}(X)}),
\end{equation}
where
\begin{equation*}
\Lambda_{s+s_{0}}=\sum\parallel
a_{ij}\parallel_{H^{s+s_{0}}(X)}+\parallel
a_{i}\parallel_{H^{s+s_{0}}(X)}+\parallel
a\parallel_{H^{s+s_{0}}(X)}
\end{equation*}
for some $s_{0}>0$, and $C_{s}$ is a constant independent of
$\varepsilon$ and $\theta$.\par
   Estimate (4.1) will first be established in the coordinates
$(\xi,\eta)$, which we have been denoting by $(x,y)$ for
convenience, and later converted into the original coordinates
$(x,y)$ of the introduction.  We will need the following
preliminary lemmas.  The first is a modification of lemma 2.3, and
the second contains standard consequences of the interpolation
inequalities for Sobolev spaces.\medskip

\noindent\textbf{Lemma 4.1.}  \textit{Let
$w\in\widetilde{H}^{(2,2)}({\Omega})$} (\textit{or
$\widetilde{H}_{x}^{(2,2)}({\Omega})$}) \textit{be such that
$yw\in L^{2}(\Omega)$, and let
$p_{1}=\varepsilon\theta\widetilde{p_{1}}$,
$p_{2}=\varepsilon\widetilde{p_{2}}$,
$p_{3}=\varepsilon\widetilde{p_{3}}$, where $\widetilde{p_{i}}\in
C^{\infty}_{c}(X)$, $i=1,2,3$.  Then for $\varepsilon$ and
$\theta$ sufficiently small, there exists a constant $M$,
independent of $\varepsilon$ and $\theta$, such that}
\begin{equation*}
\parallel w\parallel+\parallel w_{y}\parallel\leq M\parallel
p_{1}w_{xyy}+p_{2}w_{x}+p_{3}w+L_{\theta}w\parallel.
\end{equation*}
\medskip

\noindent\textit{Proof.}  Assume temporarily that $w\in
\widetilde{C}^{\infty}(\overline{\Omega})$ (or
$\widetilde{C}_{x}^{\infty}(\overline{\Omega})$).  The properties
of $p_{2}$ and $p_{3}$ guarantee that lemma 2.3 holds for the
operator $p_{2}\partial_{x}+p_{3}+L_{\theta}$.  Therefore
\begin{equation}
(Aw+Cw_{y}+Dw_{yy},p_{2}w_{x}+p_{3}w+L_{\theta}w)\geq
\end{equation}
\begin{equation*}
C_{1}[\parallel w\parallel^{2}+\parallel w_{y}\parallel^{2}
+\theta(\parallel w_{x}\parallel^{2}+\parallel
w_{xy}\parallel^{2}+\parallel w_{yy}\parallel^{2})],
\end{equation*}
where $A,C,D$, and $C_{1}$ were given in lemma 2.3.  Furthermore,
integrating by parts yields
\begin{eqnarray}
(Aw+Cw_{y}+Dw_{yy},p_{1}w_{xyy})&=&\int\int_{\Omega}[-\frac{1}{2}
(Dp_{1})_{x}]w_{yy}^{2}+[-Cp_{1}]w_{xy}w_{yy}\nonumber
\\& &+[\frac{1}{2}(Cp_{1})_{xy}
+\frac{1}{2}(Ap_{1})_{x}]w_{y}^{2}\\
&
&+[(Ap_{1})_{y}]w_{x}w_{y}+[-\frac{1}{2}(Ap_{1})_{xyy}]w^{2}\nonumber
\end{eqnarray}
All the boundary integrals vanish since $p_{1}\in
C_{c}^{\infty}(X)$.  Moreover, the properties of $p_{1}$ guarantee
that by choosing $\varepsilon$ and $\theta$ sufficiently small, we
obtain the following by adding (4.2) and (4.3),
\begin{equation*}
(Aw+Cw_{y}+Dw_{yy},p_{1}w_{xyy}+p_{2}w_{x}+p_{3}w+L_{\theta}w)\geq
\end{equation*}
\begin{equation*}
C_{1}[\parallel w\parallel^{2}+\parallel w_{y}\parallel^{2}
+\theta(\parallel w_{x}\parallel^{2}+\parallel
w_{xy}\parallel^{2}+\parallel w_{yy}\parallel^{2})].
\end{equation*}
Then an application of Cauchy's inequality, and the use of an
approximating sequence $\{w_{k}\}_{k=1}^{\infty}$, as was
constructed in proposition 2.1, removes the assumption
that\linebreak $w\in \widetilde{C}^{\infty}(\overline{\Omega})$
(or $\widetilde{C}_{x}^{\infty}(\overline{\Omega})$) and completes
the proof.\qed

\noindent\textbf{Lemma 4.2 [24].}  \textit{Let $u,v\in H^{s}(X)$.}

$i)$\textit{  If $0\leq i\leq j \leq s$, then there exists a
constant $\mathcal{M}_{i,j,s}$ such that}
\begin{equation*}
\parallel u\parallel_{H^{j}(X)}\leq \mathcal{M}_{i,j,s}\parallel u\parallel
^{\frac{s-j}{s-i}}_{H^{i}(X)}\parallel
u\parallel^{\frac{j-i}{s-i}}_{H^{s}(X)}.
\end{equation*}

$ii)$\textit{  If $\alpha$ and $\beta$ are multi-indices such that
$|\alpha|+|\beta|=s$, then there exists a constant
$\mathcal{M}_{s}$ such that}
\begin{equation*}
\parallel\partial^{\alpha}u\partial^{\beta}v\parallel_{L^{2}(X)}\leq
\mathcal{M}_{s}(|u|_{L^{\infty}(X)}\parallel
v\parallel_{H^{s}(X)}+
\parallel u\parallel_{H^{s}(X)}|v|_{L^{\infty}(X)}).
\end{equation*}

$iii)$\textit{  Let $\Gamma\subset\mathbb{R}^{N}$ be compact and
contain the origin, and let $G\in C^{\infty}(\Gamma)$.  If $u\in
H^{s+2}(X,\Gamma)$ and $\parallel u\parallel _{H^{2}(X)}\leq
\mathcal{C}$ for some fixed $\mathcal{C}$, then there exists a
constant $\mathcal{M}_{s}$ such that}
\begin{equation*}
\parallel G\circ u\parallel_{H^{s}(X)}\leq\mathrm{Vol}(X)|G(0)|
+\mathcal{M}_{s}\parallel u\parallel_{H^{s+2}(X)}.
\end{equation*}\smallskip

   Estimate (4.1) will be established by induction on $s$, and we
begin by estimating the $x$-derivatives.  Let
$\parallel\cdot\parallel_{s,X}$ denote $\parallel\cdot\parallel
_{H^{s}(X)}$, and $|\cdot|_{\infty}$ denote
$|\cdot|_{L^{\infty}(X)}$.\medskip

\noindent\textbf{Proposition 4.1.}  \textit{Let $u$ and $f$ be as
in theorem} 3.2.  \textit{If $\varepsilon=\varepsilon(s)$ is
sufficiently small, then}
\begin{equation*}
\parallel\partial^{s}_{x}u\parallel+\parallel\partial^{s}_{x}u_{y}\parallel
\leq C_{s}(\parallel f\parallel_{s}+\parallel
u\parallel_{s-1,X}+\Lambda_{s+2}\parallel f\parallel_{2,X}),
\end{equation*}
\textit{for $s\leq r-6$, where $C_{s}$ is independent of
$\varepsilon$ and $\theta$, and}
\begin{equation*}
\Lambda_{s+2}=\sum\parallel a_{ij}\parallel_{s+2,X}+\parallel
a_{i}\parallel_{s+2,X}+\parallel a\parallel_{s+2,X}.
\end{equation*}

\noindent\textit{Proof.}   We proceed by induction on $s$.  The
case $s=0$ follows from lemma 2.3.  Differentiate $L_{\theta}u=f$
$s$-times with respect to $x$ and put $w=\partial_{x}^{s}u$, then
\begin{equation}
-\theta
w_{xxyy}+a_{11}w_{xx}+a_{22}w_{yy}+(a_{1}+s\partial_{x}a_{11})w_{x}
+a_{2}w_{y}+a_{s}w\text{ }\text{ }\text{ }\text{ }\text{ }\text{ }
\text{ }\text{ }\text{ }\text{ }\text{ }\text{ }\text{ }\text{
}\text{ }\text{ }\text{ }
\end{equation}
\begin{equation*}
=\partial_{x}^{s}f-\sum_{i=0}^{s-1}\partial_{x}^{i}(\partial_{x}
a_{22}\partial_{x}^{s-1-i}u_{yy})-\sum_{i=0}^{s-1}\partial_{x}^{i}(\partial_{x}
a_{2}\partial_{x}^{s-1-i}u_{y})
-\sum_{i=0}^{s-1}\partial_{x}^{i}(\partial_{x}
a_{s-1-i}\partial_{x}^{s-1-i}u),
\end{equation*}
\begin{equation*}
:=f_{s}\text{ }\text{ }\text{ }\text{ }\text{ }\text{ }\text{ }
\text{ }\text{ }\text{ }\text{ }\text{ }\text{ }\text{ }\text{ }
\text{ }\text{ }\text{ }\text{ }\text{ }\text{ }\text{ }\text{
}\text{ }\text{ }\text{ }\text{ }\text{ }\text{ }\text{ }\text{
}\text{ }\text{ }\text{ }\text{ }\text{ }\text{ }\text{ }\text{
}\text{ }\text{ }\text{ }\text{ }\text{ }\text{ }\text{ }\text{
}\text{ }\text{ }\text{ }\text{ }\text{ }\text{ }\text{ }\text{
}\text{ }\text{ }\text{ }\text{ }\text{ }\text{ }\text{ }\text{
}\text{ }\text{ }\text{ }\text{ }\text{ }\text{ }\text{ }\text{
}\text{ }\text{ }\text{ }\text{ }\text{ }\text{ }\text{ }\text{
}\text{ }\text{ }\text{ }\text{ }\text{ }\text{ }\text{ }\text{
}\text{ }\text{ }\text{ }\text{ }\text{ }\text{ }\text{ }\text{
}\text{ }\text{ }\text{ }\text{ }\text{ }\text{ }\text{ }
\end{equation*}
where $a_{s}=a+s\partial_{x}
a_{1}+\frac{s(s-1)}{2}\partial_{x}^{2}a_{11}$.  A calculation
shows that
\begin{equation*}
\sum_{i=0}^{s-1}\partial_{x}^{i}(\partial_{x}
a_{22}\partial_{x}^{s-1-i}u_{yy})=s\partial_{x}a_{22}
\partial_{x}^{s-1}u_{yy}
+\frac{s(s-1)}{2}\partial_{x}^{2}a_{22}
\partial_{x}^{s-2}u_{yy}
\end{equation*}
\begin{equation*}
\text{ }\text{ }\text{ }\text{ }\text{ }\text{ }\text{ } \text{
}\text{ }\text{ }\text{ }\text{ }\text{ }\text{ }\text{ }\text{
}\text{ }\text{ }\text{ }\text{ }\text{ }\text{ }\text{ }\text{
}\text{ }\text{ }\text{ } +\sum_{i=2}^{s-1}\sum_{j=2}^{i}
\left(\begin{array}{c}i\\j\end{array}\right)
\partial_{x}^{j+1}a_{22}\partial_{x}^{s-1-j}u_{yy}.
\end{equation*}
Note that the term $\partial_{x}^{s-1}u_{yy}$ contains too many
derivatives. However, since \linebreak $a_{22}=1+O(\varepsilon)$,
we can solve for $\partial_{x}^{s-1}u_{yy}$ in (4.4) with $s$
replaced by $s-1$ to obtain a more manageable expression:
\begin{equation*}
\partial_{x}^{s-1}u_{yy}=\frac{1}{a_{22}}[\theta
w_{xyy}-a_{11}w_{x}-(a_{1}+s\partial_{x}a_{11})w-a_{2}\partial_{x}^{s-1}u_{y}
-a_{s-1}\partial_{x}^{s-1}u+f_{s-1}].
\end{equation*}
Substituting back into (4.4), we have
\begin{equation*}
\text{ }\text{ }
\frac{s\theta\partial_{x}a_{22}}{a_{22}}w_{xyy}+(s\partial_{x}a_{11}-
\frac{sa_{11}\partial_{x}a_{22}}{a_{22}})w_{x}+(a_{s}-a
-\frac{s\partial_{x}a_{22}}{a_{22}}(a_{1}-s\partial_{x}a_{11}))w+L_{\theta}w
\end{equation*}
\begin{equation*}
=\partial_{x}^{s}f-\frac{s(s\!-\!1)}{2}\partial_{x}^{2}a_{22}\partial_{x}^{s-2}u_{yy}
-\sum_{i=2}^{s-1}\!\sum_{j=2}^{i}
\left(\begin{array}{c}i\\j\end{array}\right)\!
\partial_{x}^{j+1}a_{22}\partial_{x}^{s-1-j}u_{yy}
-\sum_{i=0}^{s-1}\!\partial_{x}^{i}(\partial_{x}
a_{2}\partial_{x}^{s-1-i}u_{y})
\end{equation*}
\begin{equation*}
-\sum_{i=0}^{s-1}\partial_{x}^{i}(\partial_{x}
a_{s-1-i}\partial_{x}^{s-1-i}u)
+\frac{s\partial_{x}a_{22}}{a_{22}}[a_{2}\partial_{x}^{s-1}u_{y}
+a_{s-1}\partial_{x}^{s-1}u-f_{s-1}]\text{ }\text{ }\text{ }
\text{ }\text{ }\text{ }\text{ }\text{ }\text{ }\text{ }\text{
}\text{ }\text{ }\text{ }\text{ }\text{ }
\end{equation*}
\begin{equation*}
:=\widetilde{f}_{s}.\text{ }\text{ }\text{ }\text{ }\text{ }\text{
}\text{ } \text{ }\text{ }\text{ }\text{ }\text{ }\text{ }\text{
}\text{ } \text{ }\text{ }\text{ }\text{ }\text{ }\text{ }\text{
}\text{ }\text{ }\text{ }\text{ }\text{ }\text{ }\text{ }\text{
}\text{ }\text{ }\text{ }\text{ }\text{ }\text{ }\text{ }\text{
}\text{ }\text{ }\text{ }\text{ }\text{ }\text{ }\text{ }\text{
}\text{ }\text{ }\text{ }\text{ }\text{ }\text{ }\text{ }\text{
}\text{ }\text{ }\text{ }\text{ }\text{ }\text{ }\text{ }\text{
}\text{ }\text{ }\text{ }\text{ }\text{ }\text{ }\text{ }\text{
}\text{ }\text{ }\text{ }\text{ }\text{ }\text{ }\text{ }\text{
}\text{ }\text{ }\text{ }\text{ }\text{ }\text{ }\text{ }\text{
}\text{ }\text{ }\text{ }\text{ }\text{ }\text{ }\text{ }\text{
}\text{ }\text{ }\text{ }\text{ }\text{ }\text{ }\text{ }\text{
}\text{ }
\end{equation*}
If $\varepsilon=\varepsilon(s)$ and $\theta$ are sufficiently
small, we can apply lemma 4.1 to obtain
\begin{equation}
\parallel\partial_{x}^{s}u\parallel+\parallel\partial_{x}^{s}u_{y}\parallel
\leq M\parallel\widetilde{f}_{s}\parallel.
\end{equation}\par
   We now estimate each term of $\widetilde{f}_{s}$.  Using lemma 4.2
$(ii)$, lemma 2.2 $(iii)$, and the fact that $\partial_{x}a_{22}$
vanishes outside of $X$, produces
\begin{equation*}
\parallel\sum_{i=2}^{s-1}\sum_{j=2}^{i}
\left(\begin{array}{c}i\\j\end{array}\right)
\partial_{x}^{j+1}a_{22}\partial_{x}^{s-1-j}u_{yy}\parallel=
\parallel\sum_{i=2}^{s-1}\sum_{j=2}^{i}
\left(\begin{array}{c}i\\j\end{array}\right)
\partial_{x}^{j+1}a_{22}\partial_{x}^{s-1-j}u_{yy}\parallel_{0,X}
\text{ }\text{ }\text{ }\text{ }\text{ }\text{ }\text{ }\text{
}\text{ }
\end{equation*}
\begin{equation*}
\text{ }\text{ }\text{ }\text{ }\text{ }\text{ }\text{ }\text{ }
\text{ }\text{ }\text{ }\text{ }\text{ }\text{ }\text{ }\text{ }
\text{ } \text{ }\text{ }\text{ }\text{ }\text{ }\text{ }\text{
}\text{ }\text{ }\text{ }\text{ }\text{ }\text{ }\text{ }\text{
}\text{ }\text{ }\text{ }\text{ }\text{ }\text{ }\text{ }\text{
}\text{ }\text{ }\text{ }\text{ }\text{ }\text{ } \leq
M_{1}(|\partial_{x}^{3}a_{22}|_{\infty}\parallel
u\parallel_{s-1,X}+\parallel\partial_{x}^{3}a_{22}\parallel_{s-1,X}|u|_{\infty})
\end{equation*}
\begin{equation*}
\text{ }\text{ }\text{ }\text{ }\text{ }\text{ }\text{ }\text{ }
\text{ }\text{ }\text{ }\text{ }\text{ }\text{ }\text{ }\text{
}\text{ }\text{ }\text{ }\text{ }\text{ }\text{ } \text{ }\text{
}\text{ }\text{ }\text{ }\text{ }\text{ }\text{ }\text{ }\text{
}\leq M_{1}^{'}(\parallel u\parallel_{s-1,X}+\parallel
a_{22}\parallel_{s+2,X}|u|_{\infty}).
\end{equation*}
A calculation shows that
\begin{equation*}
\sum_{i=0}^{s-1}\partial_{x}^{i}(\partial_{x}
a_{2}\partial_{x}^{s-1-i}u_{y})=s\partial_{x}a_{2}\partial_{x}^{s-1}u_{y}
+\sum_{i=1}^{s-1}\sum_{j=1}^{i}\left(\begin{array}{c}i\\j\end{array}\right)
\partial_{x}^{j+1}a_{2}\partial_{x}^{s-1-j}u_{y}.
\end{equation*}
Then using the same procedure as above, we have
\begin{equation*}
\parallel\sum_{i=0}^{s-1}\partial_{x}^{i}(\partial_{x}
a_{2}\partial_{x}^{s-1-i}u_{y})\parallel\leq M_{2}\parallel
\partial_{x}^{s-1}u_{y}\parallel+M_{2}^{'}(\parallel u\parallel
_{s-1,X}+\parallel a_{2}\parallel_{s+2,X}|u|_{\infty}).
\end{equation*}
Furthermore, the following estimates are obtained in the same way:
\begin{equation*}
\parallel\sum_{i=0}^{s-1}\partial_{x}^{i}(\partial_{x}
a_{s-1-i}\partial_{x}^{s-1-i}u)\parallel+\parallel
\frac{s\partial_{x}a_{22}}{a_{22}}\sum_{i=0}^{s-2}\partial_{x}^{i}(\partial_{x}
a_{s-2-i}\partial_{x}^{s-2-i}u)\parallel
\end{equation*}
\begin{equation*}
\leq M_{3}(\parallel u\parallel_{s-1,X}+(\parallel
a\parallel_{s+2,X}+\parallel a_{1}\parallel_{s+2,X}+\parallel
a_{11}\parallel_{s+2,X})|u|_{\infty}),\text{ }\text{ }\text{
}\text{ }\text{ }\text{ }
\end{equation*}
and
\begin{equation*}
\parallel\frac{s\partial_{x}a_{22}}{a_{22}}\sum_{i=0}^{s-2}\partial_{x}^{i}(\partial_{x}
a_{2}\partial_{x}^{s-2-i}u_{y})\parallel\leq M_{4}(\parallel
u\parallel_{s-1,X}+\parallel a_{2}\parallel_{s+2,X}|u|_{\infty}).
\end{equation*}
Also, since
\begin{equation*}
\sum_{i=0}^{s-2}\partial_{x}^{i}(\partial_{x}
a_{22}\partial_{x}^{s-2-i}u_{yy})=(s-1)\partial_{x}a_{22}\partial_{x}^{s-2}u_{yy}
+\sum_{i=1}^{s-2}\sum_{j=1}^{i}\left(\begin{array}{c}i\\j\end{array}\right)
\partial_{x}^{j+1}a_{22}\partial_{x}^{s-2-j}u_{yy}
\end{equation*}
and $\partial_{x}a_{22}=O(\varepsilon)$, we find that
\begin{equation*}
\parallel\frac{s\partial_{x}a_{22}}{a_{22}}\sum_{i=0}^{s-2}\partial_{x}^{i}(\partial_{x}
a_{22}\partial_{x}^{s-2-i}u_{yy})\parallel \leq\varepsilon
s^{2}M_{5}\parallel\partial_{x}^{s-2}u_{yy}\parallel_{0,X}\text{ }
\text{ }\text{ }\text{ }\text{ }\text{ }\text{ }\text{ }\text{
}\text{ }
\end{equation*}
\begin{equation*}
\text{ }\text{ }\text{ }\text{ }\text{ }\text{ }\text{ }\text{
}\text{ }\text{ }\text{ }\text{ }\text{ }\text{ }\text{ }\text{
}\text{ }\text{ }\text{ }\text{ }\text{ }\text{ }\text{ }\text{
}\text{ }\text{ }\text{ }\text{ }\text{ } \text{ }\text{ }\text{
}\text{ }\text{ }\text{ }\text{ }\text{ }\text{ }\text{ }\text{
}\text{ }\text{ }\text{ }\text{ }\text{ }\text{ }\text{ }\text{
}\text{ }\text{ }\text{ } +M_{5}^{'}(\parallel
u\parallel_{s-1,X}+\parallel a_{22}\parallel_{s+2,X}|u|_{\infty}),
\end{equation*}
where $M_{5}$ is independent of $\varepsilon$ and $s$.\par
   Summing the above estimates produces:
\begin{equation}
\parallel\widetilde{f}_{s}\parallel\leq \!M_{6}(\parallel
f\parallel_{s}\!+\!\parallel u\parallel_{s-1,X}
\!+\!\parallel\partial_{x}^{s-1}u_{y}\parallel\! +\varepsilon
s^{2}\parallel\partial_{x}^{s-2}u_{yy}\parallel_{0,X}\!+\Lambda_{s+2}
|u|_{\infty}).
\end{equation}
Therefore, if we estimate
$\parallel\partial_{x}^{s-2}u_{yy}\parallel_{0,X}$ appropriately
and show that
\begin{equation*}
|u|_{\infty}\leq M_{7}\parallel f\parallel_{2,X},
\end{equation*}
the proof will be complete by induction.\par
   We now estimate
$\parallel\partial_{x}^{s-2}u_{yy}\parallel_{0,X}$.  Differentiate
the equation,
\begin{equation*}
\widetilde{L}u_{yy}:=-\theta u_{xxyy}+a_{22}u_{yy}=f-a_{11}u_{xx}
-a_{1}u_{x}-a_{2}u_{y}-au:=\widetilde{g},
\end{equation*}
with respect to $x$ $(s-2)$-times, then
\begin{equation*}
\widetilde{L}\partial_{x}^{s-2}u_{yy}=\partial_{x}^{s-2}\widetilde{g}
-\sum_{i=0}^{s-3}\partial_{x}^{i}(\partial_{x}a_{22}\partial_{x}^{s-3-i}u_{yy})
:=\widetilde{g}_{s-2}.
\end{equation*}
Multiply the above equation by $\partial_{x}^{s-2}u_{yy}$ and
integrate by parts in $X$ to obtain,
\begin{equation*}
\parallel\partial_{x}^{s-2}u_{yy}\parallel_{0,X}\leq M_{8}\parallel
\widetilde{g}_{s-2}\parallel_{0,X}.
\end{equation*}\par
   We now estimate $\parallel\widetilde{g}_{s-2}\parallel_{0,X}$.
Using the same methods as above, we have
\begin{equation*}
\parallel\partial_{x}^{s-2}(a_{1}u_{x}+a_{2}u_{y}+au)
+\sum_{i=0}^{s-3}\partial_{x}^{i}(\partial_{x}a_{22}\partial_{x}^{s-3-i}u_{yy})
\parallel_{0,X}\leq M_{9}(\parallel
u\parallel_{s-1,X}+\Lambda_{s+2}|u|_{\infty}).
\end{equation*}
Furthermore,
\begin{equation*}
\partial_{x}^{s-2}(a_{11}u_{xx})=a_{11}\partial_{x}^{s}u
+\sum_{i=1}^{s-2}\left(\begin{array}{c}s-2\\i\end{array}\right)
\partial_{x}^{i}a_{11}\partial_{x}^{s-2-i}u_{xx};
\end{equation*}
thus,
\begin{equation*}
\parallel\partial_{x}^{s-2}(a_{11}u_{xx})\parallel_{0,X}\leq
M_{10}(\parallel\partial_{x}^{s}u\parallel_{0,X}+\parallel
u\parallel_{s-1,X}+\Lambda_{s+2}|u|_{\infty}).
\end{equation*}
It follows that
\begin{equation}
\parallel\partial_{x}^{s-2}u_{yy}\parallel_{0,X}
\leq M_{11}(\parallel\partial_{x}^{s}u\parallel_{0,X}+\parallel
u\parallel_{s-1,X}+\Lambda_{s+2}|u|_{\infty}).
\end{equation}\par
   The coefficient of
$\parallel\partial_{x}^{s-2}u_{yy}\parallel_{0,X}$ in (4.6) is
$\varepsilon s^{2}M_{6}$.  If $\varepsilon=\varepsilon(s)$ is
chosen sufficiently small so that $\varepsilon
s^{2}MM_{6}M_{11}<\frac{1}{2}$, we can then bring $\varepsilon
s^{2}MM_{6}M_{11}\parallel\partial_{x}^{s}u\parallel_{0,X}$ from
(4.7) to the left-hand side of (4.5), so that by induction on $s$
\begin{equation*}
\parallel\partial_{x}^{s}u\parallel+\parallel\partial_{x}^{s}u_{y}\parallel
\leq M_{6}^{'}(\parallel f\parallel_{s}+\parallel
u\parallel_{s-1,X}+\Lambda_{s+2}(|u|_{\infty}+\parallel
f\parallel_{2,X})).
\end{equation*}\par
   We now estimate $|u|_{\infty}$ to complete the proof.  The
above methods can be used to show that
\begin{equation*}
\parallel
u\parallel_{2,X}\leq M_{12}\parallel f\parallel_{2,X}.
\end{equation*}
Then by the Sobolev lemma,
\begin{equation*}
\text{ }\text{ }\text{ }\text{ }\text{ }\text{ }\text{ } \text{
}\text{ }\text{ }\text{ }\text{ }\text{ }\text{ }\text{ }\text{
}\text{ }\text{ }\text{ }\text{ }\text{ }\text{ } \text{ }\text{
}\text{ }\text{ }\text{ }\text{ }\text{ }\text{ }\text{ }\text{
}\text{ }|u|_{\infty}\leq M_{13}\parallel u\parallel_{2,X}\leq
M_{13}^{'}\parallel f\parallel_{2,X}.\text{ }\text{ }\text{
}\text{ }\text{ }\text{ }\text{ }\text{ }\text{ }\text{ }\text{
}\text{ }\text{ }\text{ }\text{ }\text{ }\text{ }\text{ }\text{
}\text{ }\text{ }\text{ }\text{ }\text{ }\text{ }\text{ }\text{
}\text{ }\Box
\end{equation*}

   We now estimate the remaining derivatives.\medskip

\noindent\textbf{Proposition 4.2.}  \textit{Let $u$, $f$, $s$, and
$\varepsilon$ be as in proposition} 4.1.  \textit{Then}
\begin{equation*}
\parallel\partial_{x}^{\alpha}\partial_{y}^{\beta}u\parallel_{0,X}
\leq C_{s}(\parallel f\parallel_{s,X}+\parallel u\parallel_{s-1,X}
+\Lambda_{s+2}\parallel f\parallel_{2,X}),
\end{equation*}
\textit{for }$\alpha+\beta\leq s$, \textit{where $C_{s}$ is
independent of $\varepsilon$ and $\theta$.}\medskip

\noindent\textit{Proof.}  The cases $\beta=0,1,2$ follow from
(4.7) and proposition 4.1.  We proceed by induction on $\beta$.
Assume that the desired estimate holds for $0\leq\alpha\leq
s-\beta$, and $0\leq\beta\leq k-1$, for some $k\leq s$.\par
   Differentiate the equation,
\begin{equation*}
\widetilde{L}u_{yy}:=-\theta u_{xxyy}+a_{22}u_{yy}=f-a_{11}u_{xx}
-a_{1}u_{x}-a_{2}u_{y}-au:=\widetilde{g},
\end{equation*}
with respect to $\partial_{x}^{\alpha}\partial_{y}^{k-2}$ where
$0\leq\alpha\leq s-k$, then
\begin{eqnarray*}
\widetilde{L}\partial_{x}^{\alpha}\partial_{y}^{k}u&=&
\partial_{x}^{\alpha}\partial_{y}^{k-2}\widetilde{g}
-\sum_{i=0}^{\alpha-1}\partial_{y}^{k-2}\partial_{x}^{i}(\partial_{x}
a_{22}\partial_{x}^{\alpha-1-i}u_{yy})-\sum_{i=0}^{k-3}\partial_{y}^{i}
(\partial_{y}a_{22}\partial_{y}^{k-3-i}\partial_{x}^{\alpha}u_{yy})\\
&:=&\widetilde{g}_{\alpha,k-2}.
\end{eqnarray*}
Multiply the above equation by
$\partial_{x}^{\alpha}\partial_{y}^{k}u$, and integrate by parts
in X to obtain
\begin{equation*}
\parallel\partial_{x}^{\alpha}\partial_{y}^{k}u\parallel_{0,X}\leq
M\parallel\widetilde{g}_{\alpha,k-2}\parallel_{0,X}.
\end{equation*}\par
   We now estimate
$\parallel\widetilde{g}_{\alpha,k-2}\parallel_{0,X}$.  Using lemma
4.2 $(ii)$, we have
\begin{eqnarray*}
\parallel\partial_{x}^{\alpha}\partial_{y}^{k-2}(a_{11}u_{xx})\!\parallel_{0,X}\!\!\!\!
&\leq&\!\!\!
M_{1}(\parallel\partial_{x}^{\alpha+2}\partial_{y}^{k-2}u\parallel_{0,X}+
\!\!\!\sum_{p\leq\alpha,\text{ }q\leq k-2\atop
(p,q)\neq(0,0)}\!\!\!\parallel\partial_{x}^{p}\partial_{y}^{q}a_{11}
\partial_{x}^{\alpha-p}\partial_{y}^{k-2-q}u_{xx}\parallel_{0,X})\\
&\leq&\!\!\!
M_{1}^{'}(\parallel\partial_{x}^{\alpha+2}\partial_{y}^{k-2}u\parallel_{0,X}
+|a_{11}|_{C^{1}(\overline{X})}\parallel\!
u\!\parallel_{s-1,X}+\parallel\!
a_{11}\!\parallel_{s,X}\!|u|_{\infty})
\end{eqnarray*}
\begin{equation*}
\text{ }\text{ }\text{ }\text{ }\text{ }\text{ }\text{ } \text{
}\text{ }\text{ }\text{ }\text{ }\text{ }\text{ }\text{ }\text{ }
\text{ } \leq
M_{1}^{''}(\parallel\partial_{x}^{\alpha+2}\partial_{y}^{k-2}u\parallel_{0,X}
+\parallel u\parallel_{s-1,X}+\Lambda_{s+2}\parallel
f\parallel_{2,X}).
\end{equation*}
Furthermore, if $\alpha<s-k$ then
$\parallel\partial_{x}^{\alpha+2}\partial_{y}^{k-2}u\parallel_{0,X}\leq
\parallel u\parallel_{s-1,X}$, and if $\alpha=s-k$ the induction
assumption implies that
\begin{equation*}
\parallel\partial_{x}^{\alpha+2}\partial_{y}^{k-2}u\parallel_{0,X}\leq
M_{2}(\parallel f\parallel_{s,X}+\parallel u\parallel_{s-1,X}+
\Lambda_{s+2}\parallel f\parallel_{2,X}).
\end{equation*}
Thus,
\begin{equation*}
\parallel\partial_{x}^{\alpha}\partial_{y}^{k-2}(a_{11}u_{xx})\parallel_{0,X}
\leq M_{3}(\parallel f\parallel_{s,X}+\parallel
u\parallel_{s-1,X}+ \Lambda_{s+2}\parallel f\parallel_{2,X}).
\end{equation*}
Moreover, the methods of proposition 4.1 may be used to estimate
the remaining terms of
$\parallel\widetilde{g}_{\alpha,k-2}\parallel_{0,X}$ by
\begin{equation*}
M_{4}(\parallel u\parallel_{s-1,X}+\Lambda_{s+2}\parallel
f\parallel_{2,X}).
\end{equation*}
The desired conclusion now follows by combining the above
estimates.\qed
   From proposition 4.2, we obtain the following Moser estimate by induction on
$s$.\medskip

\noindent\textbf{Theorem 4.1.}  \textit{Let $u$ and $f$ be as in
theorem} 3.2. \textit{If $\varepsilon=\varepsilon(s)$ is
sufficiently small, then}
\begin{equation*}
\parallel u\parallel_{s,X}\leq C_{s}(\parallel f\parallel_{s,X}
+\Lambda_{s+2}\parallel f\parallel_{2,X}),
\end{equation*}
\textit{for $s\leq r-6$, where $C_{s}$ is independent of
$\varepsilon$ and $\theta$.}\medskip

   The estimate of theorem 4.1 is in terms of the variables
$(\xi,\eta)$ of lemma 2.2, which we have been denoting by $(x,y)$
for convenience.  We now swap notation and denote the original
variables of (2.1) by $(x,y)$, and the change of variables by
$(\xi,\eta)$.  Furthermore, let $\parallel\cdot\parallel_{s}$ and
$\parallel\cdot\parallel_{s}^{'}$ denote the $H^{s}(X)$ norm with
respect to the variables $(x,y)$ and $(\xi,\eta)$ respectively.
Similarly for $\Lambda_{s}$ and $\Lambda_{s}^{'}$.  We now obtain
the analogue of theorem 4.1 with respect to the variables $(x,y)$.
We will need the following lemma.\medskip

\noindent\textbf{Lemma 4.3.}  \textit{If
$\varepsilon=\varepsilon(s)$ is sufficiently small, then}
\begin{equation*}
\parallel\xi_{x}\parallel_{s}\leq C_{s}(\parallel
a_{12}\parallel_{s+3}+\parallel a_{22}\parallel_{s+5}),
\end{equation*}
\textit{for $s\leq r-7$, where $C_{s}$ is independent of
$\varepsilon$ and $\theta$.}\medskip

\noindent\textit{Proof.}  We prove the estimate by induction on
$s$. The case $s=0$ follows from the estimate,
\begin{equation*}
0<M_{1}\leq|\xi_{x}|\leq M_{2},
\end{equation*}
obtained in the proof of lemma 2.2.  Now assume that the estimate
holds for $s-1$.\par
   We first estimate the $x$-derivatives.  Differentiate the
equation
\begin{equation}
(\frac{a_{12}}{a_{22}})(\xi_{x})_{x}+(\xi_{x})_{y}=
-(\frac{a_{12}}{a_{22}})_{x}\xi_{x}
\end{equation}
with respect to $x$ $s$-times to obtain,
\begin{equation*}
(\frac{a_{12}}{a_{22}})(\partial_{x}^{s}\xi_{x})_{x}+
(\partial_{x}^{s}\xi_{x})_{y}=
-\partial_{x}^{s}[(\frac{a_{12}}{a_{22}})_{x}\xi_{x}]
-\sum_{i=0}^{s-1}\partial_{x}^{i}((\frac{a_{12}}{a_{22}})_{x}
\partial_{x}^{s-i}\xi_{x}):=h_{s}.
\end{equation*}
Then estimating $\partial_{x}^{s}\xi_{x}$ along the
characteristics of (4.8) as in the proof of lemma 2.2, we have
\begin{equation*}
|\partial_{x}^{s}\xi_{x}|_{C^{0}(\overline{X})}\leq
M_{3}|h_{s}|_{C^{0}(\overline{X})}.
\end{equation*}
Recalling that $a_{12}=O(\varepsilon)$, and using the analogue of
lemma 4.2 $(ii)$ for $C^{s}(\overline{X})$ norms in the same way
that the Sobolev version was used in proposition 4.1, produces
\begin{eqnarray*}
|h_{s}|_{C^{0}(\overline{X})}&\leq& \varepsilon(s+1)M_{4}
|\partial_{x}^{s}\xi_{x}|_{C^{0}(\overline{X})}\\
& &+M^{'}_{4}
(|(\frac{a_{12}}{a_{22}})_{xx}|_{C^{0}(\overline{X})}
|\xi_{x}|_{C^{s-1}(\overline{X})}+
|(\frac{a_{12}}{a_{22}})_{xx}|_{C^{s-1}(\overline{X})}|\xi_{x}|_{C^{0}(\overline{X})}).
\end{eqnarray*}
Therefore, if $\varepsilon$ is small enough to guarantee that
$\varepsilon(s+1)M_{3}M_{4}<\frac{1}{2}$, we can bring
$\varepsilon(s+1)M_{3}M_{4}|\partial_{x}^{s}\xi_{x}|_{C^{0}(\overline{X})}$
to the left-hand side:
\begin{equation}
|\partial_{x}^{s}\xi_{x}|_{C^{0}(\overline{X})}\leq M_{5}
(|\xi_{x}|_{C^{s-1}(\overline{X})}+|\frac{a_{12}}{a_{22}}|_{C^{s+1}(\overline{X})}).
\end{equation}\par
   We now estimate the remaining derivatives.  Assume that
\begin{equation}
|\partial_{x}^{\alpha}\partial_{y}^{\beta}\xi_{x}|_{C^{0}(\overline{X})}
\leq
M_{6}(|\xi_{x}|_{C^{s-1}(\overline{X})}+|\frac{a_{12}}{a_{22}}|_{C^{s+1}(\overline{X})})
\end{equation}
for all $0\leq\alpha\leq s-\beta$, $0\leq\beta\leq s-1$.  The case
$\beta=0$ is given by (4.9).  Differentiate (4.8) with respect to
$\partial_{x}^{\alpha-1}\partial_{y}^{\beta}$ to obtain
\begin{equation*}
\partial_{x}^{\alpha-1}\partial_{y}^{\beta+1}\xi_{x}=
-\partial_{y}^{\beta}[(\frac{a_{12}}{a_{22}})(\partial_{x}^{\alpha-1}\xi_{x})_{x}]
-\partial_{y}^{\beta}\partial_{x}^{\alpha-1}[(\frac{a_{12}}{a_{22}})_{x}\xi_{x}]
-\partial_{y}^{\beta}\sum_{i=0}^{s-1}\partial_{x}^{i}((\frac{a_{12}}{a_{22}})_{x}
\partial_{x}^{\alpha-1-i}\xi_{x}).
\end{equation*}
Using assumption (4.10) on the first term on the right-hand side,
and applying lemma 4.2 $(ii)$ to the remaining terms, we find
\begin{equation*}
|\partial_{x}^{\alpha-1}\partial_{y}^{\beta+1}\xi_{x}|_{C^{0}(\overline{X})}\leq
M_{7}(|\xi_{x}|_{C^{s-1}(\overline{X})}+
|\frac{a_{12}}{a_{22}}|_{C^{s+1}(\overline{X})}).
\end{equation*}
Thus, by induction on $\beta$, estimate (4.10) holds for all
$0\leq\alpha\leq s-\beta$, $0\leq\beta\leq s$.\par
   By induction on $s$, (4.10) implies that
\begin{equation*}
|\xi_{x}|_{C^{s}(\overline{X})}\leq M_{8}
|\frac{a_{12}}{a_{22}}|_{C^{s+1}(\overline{X})}.
\end{equation*}
Then the Sobolev lemma gives
\begin{equation*}
\parallel\xi_{x}\parallel_{s}\leq M_{9}
\parallel\frac{a_{12}}{a_{22}}\parallel_{s+3}.
\end{equation*}
Moreover, by lemma 4.2 $(ii)$ and $(iii)$ we have
\begin{eqnarray*}
\text{ }\text{ }\text{ }\text{ }\text{ }\text{ }\text{ }\text{
}\text{ }\text{ }\text{ }\text{ }\text{ }\text{ }\text{ }\text{
}\text{ }
\parallel\frac{a_{12}}{a_{22}}\parallel_{s+3}&\leq&
M_{10}(|a_{12}|_{\infty}\parallel\frac{1}{a_{22}}\parallel_{s+3}
+\parallel a_{12}\parallel_{s+3}|\frac{1}{a_{22}}|_{\infty})\\
&\leq&M_{11}(\parallel a_{22}\parallel_{s+5}+\parallel
a_{12}\parallel_{s+3}).\text{ }\text{ }\text{ }\text{ }\text{
}\text{ }\text{ }\text{ }\text{ }\text{ }\text{ }\text{ }\text{
}\text{ }\text{ }\text{ }\text{ }\text{ }\text{ }\text{ }\text{
}\text{ }\text{ }\text{ }\text{ }\text{ }\text{ }\text{ }\text{
}\text{ }\text{ }\text{ }\Box
\end{eqnarray*}

\noindent\textbf{Theorem 4.2.}  \textit{Let $u$ and $f$ be as in
theorem} 3.2. \textit{If $\varepsilon=\varepsilon(s)$ is
sufficiently small, then}
\begin{equation*}
\parallel u\parallel_{s}\leq C_{s}(\parallel f\parallel_{s}
+\Lambda_{s+11}\parallel f\parallel_{2}),
\end{equation*}
\textit{for $s\leq r-13$, where $C_{s}$ is independent of
$\varepsilon$ and $\theta$.}\medskip

\noindent\textit{Proof.}  Let $\sigma$ be a multi-index with
$|\sigma|\leq s$.  A calculation shows that
\begin{equation*}
\parallel\partial_{x,y}^{\sigma}u\parallel\text{ }\!\!\!
\leq M_{1}\parallel\sum_{|\gamma|\leq
s}G_{\gamma}\partial_{\xi,\eta}^{\gamma}u\parallel,
\end{equation*}
where $G_{\gamma}$ are polynomials in the variables
$x_{\xi}^{-1}=\xi_{x}$, $\partial_{\xi,\eta}^{\gamma_{1}}x_{\xi}$,
and $\partial_{\xi,\eta}^{\gamma_{2}}x_{\eta}$, such that
$\sum_{i}|\gamma_{i}|\leq s-|\gamma|$ for each term of
$G_{\gamma}$.  Then using lemma 4.2 $(ii)$ and $(iii)$, we find
that
\begin{equation*}
\parallel\partial_{x,y}^{\sigma}u\parallel\text{ }\!\!\!\leq M_{2}(\parallel u
\parallel_{s}^{'}+(\parallel x_{\xi}\parallel_{s+2}^{'}+\parallel
x_{\eta}\parallel_{s+2}^{'})|u|_{\infty}).
\end{equation*}
Similarly,
\begin{equation}
\parallel\partial_{\xi,\eta}^{\sigma}u\parallel\text{ }\!\!\!
\leq M_{3}(\parallel u\parallel_{s}+(\parallel
\xi_{x}\parallel_{s+2}+\parallel
\xi_{y}\parallel_{s+2})|u|_{\infty}).
\end{equation}
Then by theorem 4.1 and the Sobolev lemma, we have
\begin{equation}
\parallel\partial_{x,y}^{\sigma}u\parallel\text{ }\!\!\!
\leq M_{4}(\parallel f\parallel_{s}^{'}+\Lambda_{s+2}^{'}
\parallel f\parallel_{2}^{'})+M_{4}^{'}
(\parallel x_{\xi}\parallel_{s+2}^{'}+\parallel
x_{\eta}\parallel_{s+2}^{'})\parallel f\parallel_{2}.
\end{equation}\par
   We now estimate the terms on the right-hand side of (4.12).
Use lemma 4.2 $(ii)$, $(iii)$, and (4.11) to obtain
\begin{eqnarray*}
\parallel x_{\xi}\parallel_{s+2}^{'}\!\!\!\text{ }\!=\!\!\text{ }\!
\parallel\frac{1}{\xi_{x}}\parallel_{s+2}^{'}\!\!\!
&\leq& \!\!\!M_{5}\parallel\xi_{x}\parallel_{s+4}^{'}\\
&\leq& \!\!\!M_{6}(\parallel\xi_{x}\parallel_{s+4}+(\parallel
\xi_{x}\parallel_{s+6}+\parallel
\xi_{y}\parallel_{s+6})|\xi_{x}|_{\infty})\\
&\leq&\!\!\!
M_{7}(\parallel\xi_{x}\parallel_{s+6}+\parallel\frac{a_{12}}{a_{22}}\xi_{x}\parallel_{s+6})\\
&\leq& \!\!\!M_{8}(\parallel a_{12}\parallel_{s+9}+\parallel
a_{22}\parallel_{s+11}).
\end{eqnarray*}
Similarly,
\begin{equation*}
\parallel x_{\eta}\parallel_{s+2}^{'}\!\!\text{ }\!=\!\!\text{ }\!\parallel
\frac{\xi_{y}}{\xi_{x}}\parallel _{s+2}^{'}\!\!\text{ }\!\leq
M_{9}(\parallel a_{12}\parallel_{s+7}+\parallel
a_{22}\parallel_{s+9}).
\end{equation*}
Furthermore,
\begin{eqnarray*}
\parallel f\parallel_{s}^{'}\!\!&\leq& \!\!M_{10}(
\parallel f\parallel_{s}+(\parallel\xi_{x}\parallel_{s+2}+\parallel\xi_{y}
\parallel_{s+2})|f|_{\infty})\\
&\leq& \!\!M_{11}(\parallel f\parallel_{s}+(\parallel
a_{12}\parallel_{s+5}+\parallel a_{22}\parallel_{s+7})\parallel
f\parallel_{2}),
\end{eqnarray*}
and hence
\begin{equation*}
\parallel f\parallel_{2}^{'} \text{ }\!\!\!\leq M_{12}
(\parallel a_{12}\parallel_{7}+\parallel
a_{22}\parallel_{9})\parallel f\parallel_{2}\text{ }\!\!\!\leq
M_{13}\parallel f\parallel_{2}.
\end{equation*}
Also,
\begin{eqnarray*}
\parallel a_{ij}\parallel_{s+2}^{'}\!\!\!
&\leq& \!\!\!M_{14} (\parallel
a_{ij}\parallel_{s+2}+(\parallel\xi_{x}\parallel_{s+4}+\parallel\xi_{y}
\parallel_{s+4})|a_{ij}|_{\infty})\\
&\leq& \!\!\!M_{15}(\parallel a_{ij}\parallel_{s+2}+\parallel
a_{12}\parallel_{s+7}+\parallel a_{22}\parallel_{s+9}),
\end{eqnarray*}
so that
\begin{equation*}
\Lambda_{s+2}^{'}\leq M_{16}\Lambda_{s+9}.
\end{equation*}
Therefore, using the above estimates and summing over all
$|\sigma|\leq s$, (4.12) produces
\begin{equation*}
\text{ }\text{ }\text{ }\text{ }\text{ }\text{ }\text{ }\text{
}\text{ }\text{ }\text{ }\text{ }\text{ }\text{ }\text{ }\text{
}\text{ }\text{ }\text{ }\text{ }\text{ }\text{ }\text{ }\text{
}\text{ }\text{ }\text{ }\text{ }\text{ }\text{ }\text{ }\text{ }
\parallel u\parallel_{s}\text{ }\!\!\!\leq M_{17}(\parallel f\parallel_{s}
+\Lambda_{s+11}\parallel f\parallel_{2}). \text{ }\text{ }\text{
}\text{ }\text{ }\text{ }\text{ }\text{ }\text{ }\text{ }\text{
}\text{ }\text{ }\text{ }\text{ }\text{ }\text{ }\text{ }\text{
}\text{ }\text{ }\text{ }\text{ }\text{ }\text{ }\text{ }\text{
}\text{ }\Box
\end{equation*}

\bigskip
\noindent\textbf{5.  The Nash-Moser Procedure}
\setcounter{equation}{0} \setcounter{section}{5}
\bigskip\bigskip
\noindent In this section we will modify the Nash-Moser iteration
procedure to obtain a solution of
\begin{equation}
\Phi(w)=0\text{ }\text{ }\text{ }\text{ in }\text{ }\text{
}X_{\infty},
\end{equation}
where $X_{\infty}\subset X$ is a sufficiently small neighborhood
of the origin that will be defined below.  In order to accommodate
the requirement (theorem 3.2) that
$\partial_{x}^{\alpha}f|_{\partial\Omega}=0$, $\alpha\leq s-1$, we
will cut off the right-hand side of the modified linearized
equation,
\begin{equation*}
L_{\theta}u=f,
\end{equation*}
near $\partial X$ at each iteration, and then estimate the error
in a smaller domain at the next step.  Furthermore, the constant
$\theta$ will be chosen sufficiently small at each iteration, to
guarantee that the procedure converges.\par
   Let $\mu>5$.  Define a sequence of domains
$\{X_{n}\}_{n=1}^{\infty}$ by
\begin{equation*}
X_{1}=X,\text{ }\text{ }\text{ }\text{
}X_{n}=(1-\sum_{i=1}^{n-1}\mu^{-i})X,
\end{equation*}
where $\lambda X=\{\lambda x\mid x\in X\}$.  Then
$X_{\infty}=(1-\frac{1}{\mu-1})X$.  In addition, let
$\frac{3}{2}<\tau<2$ and define $\mu_{n}=\mu^{\tau^{n+n_{0}}}$,
where $n_{0}>0$ will be chosen sufficiently large.\par
   We now construct smoothing operators on $L^{2}(X_{n})$.  Fix
$\widehat{\psi}\in C^{\infty}_{c}(\mathbb{R}^{2})$ such that
$\widehat{\psi}\equiv 1$ in some neighborhood of the origin. Let
$\psi(x)=\int\int_{\mathbb{R}^{2}}\widehat{\psi}(\eta)e^{2\pi
i\eta\bullet x}d\eta$ be the inverse Fourier transform of
$\widehat{\psi}$.  Then $\psi$ is a Schwartz function and
satisfies $\int\int_{\mathbb{R}^{2}}\psi(x)dx=1$, and
$\int\int_{\mathbb{R}^{2}}x^{\alpha}\psi(x)dx=0$ for any
multi-index $\alpha\neq 0$.  If $g\in L^{2}(\mathbb{R}^{2})$ and
$\gamma\geq 1$, we define the smoothing operators
$S_{\gamma}^{'}:L^{2}(\mathbb{R}^{2})\rightarrow
H^{\infty}(\mathbb{R}^{2})$ by
\begin{equation*}
(S_{\gamma}^{'}g)(x)=\gamma^{2}\int\int_{\mathbb{R}^{2}}\psi(\gamma(x-y))g(y)dy.
\end{equation*}
Then we have (see [22]),\medskip

\noindent\textbf{Lemma 5.1.}  \textit{Let $a,b\in\mathbb{Z}_{\geq
0}$ and $g\in H^{a}(\mathbb{R}^{2})$, then}

$i)$ $\parallel
S_{\gamma}^{'}g\parallel_{H^{b}(\mathbb{R}^{2})}\leq
C_{a,b}\parallel g\parallel_{H^{a}(\mathbb{R}^{2})}$, $b\leq a$,

$ii)$ $\parallel
S_{\gamma}^{'}g\parallel_{H^{b}(\mathbb{R}^{2})}\leq
C_{a,b}\gamma^{b-a}\parallel g\parallel_{H^{a}(\mathbb{R}^{2})}$,
$a\leq b$,

$iii)$ $\parallel g-
S_{\gamma}^{'}g\parallel_{H^{b}(\mathbb{R}^{2})}\leq
C_{a,b}\gamma^{b-a}\parallel g\parallel_{H^{a}(\mathbb{R}^{2})}$,
$b\leq a$.\medskip

\noindent To complete the construction, we also need the following
extension theorem.\medskip

\noindent\textbf{Theorem 5.1 [23].}  \textit{Let $D$ be a bounded
convex domain in $\mathbb{R}^{2}$ with Lipschitz smooth boundary.
Then there exists a linear operator $T_{D}:L^{2}(D)\rightarrow
L^{2}(\mathbb{R}^{2})$ such that:}

$i)$  $T_{D}(g)|_{D}=g$,

$ii)$  $T_{D}:H^{a}(D)\rightarrow H^{a}(\mathbb{R}^{2})$
\textit{continuously for each $a\in\mathbb{Z}_{\geq 0}$}.\medskip

\noindent To obtain smoothing operators on $X_{n}$,
$S_{n}:L^{2}(X_{n})\rightarrow H^{\infty}(X_{n})$, we
set\linebreak $S_{n}g=(S^{'}_{\mu_{n}}T_{X_{n}}g)|_{X_{n}}$.
Furthermore, it is clear that the corresponding results of lemma
5.1 hold for each $S_{n}$.\par
   We now set up the iteration procedure.  A sequence of functions
$\{w_{n}\}_{n=1}^{\infty}$ will be shown to converge to a solution
of (5.1), and shall be defined inductively as follows.  Set
$w_{1}=0$ and suppose that $w_{j}$, $j\leq n$, are already defined
in $X_{j}$, then set\linebreak $w_{n+1}=w_{n}+S_{n}u_{n}$ in
$X_{n+1}$, where $u_{n}$ is defined in $X_{n}$ and will be
specified below. Set $f_{n}=-\Phi(w_{n})$ in $X_{n}$, and let
$\phi_{n}$ be a $C^{\infty}$ cut off function given by
\begin{equation*}
\phi_{n}(x)=\begin{cases}
1 & \text{if $x\in X_{n+1}$},\\
0 & \text{if $x\in X-X_{n}$,}
\end{cases}
\end{equation*}
such that
\begin{equation*}
|\phi_{n}|_{C^{s}(X_{n})}\leq M_{s}\mu^{sn}.
\end{equation*}
Let
\begin{equation*}
L(w_{n})=\sum_{i,j}a_{ij}(w_{n})\partial_{ij}+\sum_{i}a_{i}(w_{n})\partial_{i}
+a(w_{n})
\end{equation*}
denote the linearization of $\Phi(w)$ evaluated at $w_{n}$, and
let $\{\theta_{n}\}_{n=1}^{\infty}$ be a sequence of positive
numbers tending towards zero that will be specified later.  Then
define $u_{n}$ in $X_{n}$ by $u_{n}=v_{n}|_{X_{n}}$, where $v_{n}$
is the solution of
\begin{equation*}
L_{\theta_{n}}(w_{n})v_{n}=\phi_{n}f_{n}\text{ }\text{ }\text{
}\text{ in }\text{ }\text{ }X,
\end{equation*}
given by theorem 2.1.  Since $\mu>5$ we have $\frac{3}{4}X\subset
X_{\infty}$.  Therefore, it follows from the definition of
$\Phi(w)$ in (1.5) that the coefficients of
$L_{\theta_{n}}(w_{n})$ are well-defined in all of $X$, even
though $w_{n}$ is only defined in $X_{n}$.\par
   For simplicity we denote the Sobolev norms
$\parallel\cdot\parallel_{H^{s}(X_{n})}$ by
$\parallel\cdot\parallel_{s}^{n}$, and the
$C^{s}(\overline{X}_{n})$ norms by $|\cdot|_{s}^{n}$.  Let
$s_{*}\in\mathbb{Z}_{\geq 0}$ be fixed such that $\Phi(0)\in
H^{s_{*}}(X)$, and define
\begin{equation*}
\sigma=n(n+1)\tau^{-(n+1+n_{0})},\text{ }\text{ }\text{ }\text{
}\text{ } \delta=\frac{16}{\tau-1}.
\end{equation*}
The convergence of the sequence $\{w_{n}\}_{n=1}^{\infty}$ to a
solution of (5.1) will follow from the following four statements.
Each will be proven by induction on $j$, for some constants
$C_{1}$, $C_{2}$, and $C_{3}$ independent of $j$ and dependent on
$\mu$ and $s_{*}$.  We shall require that $s\leq
s_{*}-18-2\delta-\frac{6\tau}{2-\tau}$ and $s_{*}\geq 22+2\delta
+\frac{6\tau}{2-\tau}$.\bigskip

   I$_{j}$:   $\parallel w_{j}\parallel_{s+15}^{j}\leq\mu_{j}^{\sigma
    s+\delta}\parallel f_{1}\parallel^{1}_{s_{*}-15}$\bigskip

   II$_{j}$:   $\parallel u_{j-1}\parallel_{s}^{j-1}\leq C_{1}\mu_{j-1}
   ^{\tau^{-1}(s-s_{*}+18+2\delta)}\parallel
   f_{1}\parallel^{1}_{s_{*}-15}$\bigskip

   III$_{j}$:   $\parallel f_{j}\parallel_{s}^{j}\leq C_{2}\mu_{j}
   ^{\tau^{-1}(s-s_{*}+18+2\delta)}\parallel
   f_{1}\parallel^{1}_{s_{*}-15}$\bigskip

   IV$_{j}$:   $\parallel w_{j}\parallel^{j}_{14}\leq C_{3}$\bigskip

   To start the induction process observe that I$_{1}$, II$_{1}$,
and IV$_{1}$ are trivial, and that III$_{1}$ holds if we set
$C_{2}=\mu_{1}$.  Now assume that I$_{j}$,$\ldots$,IV$_{j}$ hold
for $1\leq j\leq n$. The next four propositions will prove the
induction step.  Note that the coefficients of $L(w_{j})$ satisfy
the conditions placed on (2.1) with $r=s_{*}-2$. Therefore, the
results of the previous sections apply to $L_{\theta_{j}}(w_{j})$,
$1\leq j\leq n$, as long as $\varepsilon(s_{*})$ and $\theta_{j}$
are sufficiently small and $s\leq s_{*}-15$.\medskip

\noindent\textbf{Proposition 5.1.}  \textit{If $s\leq s_{*}-15$
and $\mu(s_{*})$ is sufficiently large, then}
\begin{equation*}
\parallel w_{n+1}\parallel_{s+15}^{n+1}\leq\mu_{n+1}^{\sigma
s+\delta}\parallel f_{1}\parallel^{1}_{s_{*}-15}.
\end{equation*}\medskip

\noindent\textit{Proof.}  We have
\begin{equation*}
\parallel w_{n+1}\parallel^{n+1}_{s+15}\leq
\parallel w_{n}\parallel_{s+15}^{n}+\parallel S_{n}u_{n}\parallel
_{s+15}^{n}.
\end{equation*}
Furthermore, by theorem 4.2 and lemma 4.2 ($iii$),
\begin{eqnarray*}
\parallel S_{n}u_{n}\parallel_{s+15}^{n}&\leq&
M_{1}\mu_{n}^{15}\parallel u_{n}\parallel_{s}^{n}\\
&\leq&M_{2}\mu_{n}^{15}(\parallel\phi_{n}f_{n}\parallel_{s}^{n}+
\parallel
w_{n}\parallel_{s+15}^{n}\parallel\phi_{n}f_{n}\parallel_{2}^{n}).
\end{eqnarray*}\par
   Using lemma 4.2 ($ii$), we obtain
\begin{eqnarray*}
\parallel\phi_{n}f_{n}\parallel_{s}^{n}&\leq&
M_{3}(\parallel
f_{n}\parallel_{s}^{n}+\parallel\phi_{n}\parallel_{s}^{n}
|f_{n}|_{0}^{n})\\
&\leq&M_{4}(\parallel
f_{n}\parallel_{s}^{n}+\parallel\phi_{n}\parallel_{s}^{n}
\parallel f_{n}\parallel_{2}^{n})\\
&\leq&M_{5}\mu^{sn}\parallel f_{n}\parallel_{s}^{n}.
\end{eqnarray*}
Moreover, by definition of $f_{n}$ and lemma 4.2 ($iii$)
\begin{equation}
\parallel f_{n}\parallel_{s}^{n}\leq M_{6}(\parallel
f_{1}\parallel_{s_{*}-15}^{n}+\parallel w_{n}\parallel_{s+4}^{n}),
\end{equation}
so that,
\begin{equation*}
\parallel\phi_{n}f_{n}\parallel_{s}^{n}\leq M_{7}\mu^{sn}
(\parallel f_{1}\parallel_{s_{*}-15}^{n}+\parallel
w_{n}\parallel_{s+4}^{n}).
\end{equation*}
Similarly, using IV$_{n}$
\begin{equation*}
\parallel\phi_{n}f_{n}\parallel_{2}^{n}\leq M_{7}\mu^{2n}
(\parallel f_{1}\parallel_{s_{*}-15}^{n}+\parallel
w_{n}\parallel_{6}^{n})\leq M_{8}\mu^{2n}.
\end{equation*}\par
   We now have
\begin{equation*}
\parallel S_{n}u_{n}\parallel_{s+15}^{n}\leq M_{9}\mu_{n}^{16}
\mu^{sn}(\parallel f_{1}\parallel_{s_{*}-15}^{1}+\parallel
w_{n}\parallel_{s+15}^{n}).
\end{equation*}
Therefore,
\begin{eqnarray*}
\parallel w_{n+1}\parallel_{s+15}^{n+1}&\leq&2M_{9}
\mu_{n}^{16}\mu^{sn}(\parallel f_{1}\parallel_{s_{*}-15}^{1}+
\parallel w_{n}\parallel_{s+15}^{n})\\
&\leq&\mu_{n}^{16}\mu^{2sn}(\parallel
f_{1}\parallel_{s_{*}-15}^{1}+
\parallel w_{n}\parallel_{s+15}^{n}),
\end{eqnarray*}
where the last inequality holds if $\mu$ is chosen so large that
$2M_{9}\mu^{-1}\leq 1$.  It follows that
\begin{equation*}
\parallel w_{n+1}\parallel_{s+15}^{n+1}\leq (\prod_{i=1}^{n}
\mu_{i}^{16}\mu^{2si})M_{10}\parallel
f_{1}\parallel_{s_{*}-15}^{1},
\end{equation*}
where
\begin{equation*}
M_{10}=1+\mu_{1}^{-16}\mu^{-2s}+\cdots+\prod_{i=1}^{n-1}\mu_{i}^{-16}\mu^{-2si}
\leq 2,
\end{equation*}
if $\mu$ is large.  Hence
\begin{eqnarray*}
\parallel w_{n+1}\parallel_{s+15}^{n+1}&\leq&2\mu^{sn(n+1)+
\frac{16}{\tau-1}(\tau^{n+1+n_{0}}-\tau^{1+n_{0}})}\parallel f_{1}
\parallel_{s_{*}-15}^{1}\\
&\leq&\mu_{n+1}^{\sigma s+\delta}\parallel f_{1}
\parallel_{s_{*}-15}^{1},
\end{eqnarray*}
where $\sigma=n(n+1)\tau^{-(n+1+n_{0})}$ and
$\delta=\frac{16}{\tau-1}$.\qed

\noindent\textbf{Proposition 5.2.}  \textit{If $s\leq
s_{*}-20-2\delta$ and $n_{0}(s_{*})$ is sufficiently large, then}
\begin{equation*}
\parallel u_{n}\parallel_{s}^{n}\leq C_{1}\mu_{n}
   ^{\tau^{-1}(s-s_{*}+18+2\delta)}\parallel
   f_{1}\parallel^{1}_{s_{*}-15},
\end{equation*}
\textit{where $C_{1}$ depends on $\mu$ and $s_{*}$.}\medskip

\noindent\textit{Proof.}  By theorem 4.2
\begin{equation*}
\parallel u_{n}\parallel_{s_{*}-15}^{n}\leq M_{1}
(\parallel\phi_{n}f_{n}\parallel_{s_{*}-15}^{n}+\parallel w_{n}
\parallel_{s_{*}}^{n}\parallel\phi_{n}f_{n}\parallel_{2}^{n}),
\end{equation*}
where $M_{1}$ depends only on $s_{*}$.  By lemma 4.2 $(ii)$,
(5.2), and I$_{n}$
\begin{eqnarray*}
\parallel\phi_{n}f_{n}\parallel_{s_{*}-15}^{n}&\leq&M_{2}(
\parallel
f_{n}\parallel_{s_{*}-15}^{n}+\parallel\phi_{n}\parallel_{s_{*}-15}^{n}
\parallel f_{n}\parallel_{2}^{n})\\
&\leq&M_{3}(1+\mu^{(s_{*}-15)n})\mu_{n}^{\sigma(s_{*}-26)+\delta}\parallel
f_{1}\parallel_{s_{*}-15}^{1}\\
&\leq&M_{4}\mu_{n}^{2s_{*}\sigma+\delta}\parallel
f_{1}\parallel_{s_{*}-15}^{1},
\end{eqnarray*}
where $M_{3}$ depends only on $s_{*}$.  Similarly, III$_{n}$
yields
\begin{equation*}
\parallel\phi_{n}f_{n}\parallel^{n}_{2}\leq
M_{5}C_{2}\mu_{n}^{2\sigma+\tau^{-1}(20-s_{*}+2\delta)}\parallel
f_{1}\parallel_{s_{*}-15}^{1}.
\end{equation*}
Therefore, for some constant $M_{6}$ depending on $\mu$ and
$s_{*}$, we have
\begin{eqnarray}
\parallel u_{n}\parallel_{s_{*}-15}^{n}&\leq&M_{6}(\mu_{n}^{
2s_{*}\sigma+\delta}+\mu_{n}^{\sigma(s_{*}-15)+\delta}\mu_{n}^{
2\sigma+\tau^{-1}(20-s_{*}+2\delta)})
\parallel f_{1}\parallel_{s_{*}-15}^{1}\\
&\leq&2M_{6}\mu_{n}^{2s_{*}\sigma+\delta}\parallel
f_{1}\parallel_{s_{*}-15}^{1},\nonumber
\end{eqnarray}
since $s_{*}\geq 20+2\delta$.  Furthermore, lemma 2.3 and
III$_{n}$ produce
\begin{equation*}
\parallel u_{n}\parallel^{n}_{0}\leq M_{7}\parallel
f_{n}\parallel^{n}_{0}\leq
M_{7}C_{2}\mu_{n}^{\tau^{-1}(18-s_{*}+2\delta)}\parallel
f_{1}\parallel_{s_{*}-15}^{1}.
\end{equation*}
Then applying lemma 4.2 $(i)$, we find
\begin{eqnarray*}
\parallel u_{n}\parallel^{n}_{s}&\leq& M_{8}(\parallel
u_{n}\parallel^{n}_{0})^{1-\frac{s}{s_{*}-15}}(\parallel
u_{n}\parallel^{n}_{s_{*}-15})^{\frac{s}{s_{*}-15}}\\
&\leq& M_{9}\mu_{n}^{\tau^{-1}(18-s_{*}+2\delta)
(1-\frac{s}{s_{*}-15})+(2s_{*}\sigma+\delta)(\frac{s}{s_{*}-15})}\parallel
f_{1}\parallel_{s_{*}-15}^{1}\\
&\leq& M_{9}\mu_{n}^{\tau^{-1}(s-s_{*}+18+2\delta)}\parallel
f_{1}\parallel_{s_{*}-15}^{1},
\end{eqnarray*}
if $\sigma$ is sufficiently small.  Note that $\sigma$ may be made
arbitrarily small by choosing $n_{0}$ sufficiently large.  We then
set $C_{2}=M_{9}$ to obtain the desired result.\qed

\noindent\textbf{Proposition 5.3.}  \textit{If $s\leq
s_{*}-18-2\delta-\frac{6\tau}{2-\tau}$, $s_{*}\geq 22+2\delta
+\frac{6\tau}{2-\tau}$, $n_{0}(s_{*})$ and $\mu(s_{*})$ are
sufficiently large, and $\varepsilon(s_{*})$ is sufficiently
small, then}
\begin{equation*}
\parallel f_{n+1}\parallel_{s}^{n+1}\leq C_{2}\mu_{n+1}
^{\tau^{-1}(s-s_{*}+18+2\delta)}\parallel
f_{1}\parallel^{1}_{s_{*}-15}.
\end{equation*}\smallskip

\noindent\textit{Proof.}  Expanding $\Phi(w_{n+1})$ in a Taylor
series yields,
\begin{equation*}
f_{n+1}=f_{n}-L(w_{n})S_{n}u_{n}+Q_{n}=f_{n}-\theta_{n}
(S_{n}u_{n})_{\eta\eta\xi\xi}-L_{\theta_{n}}(w_{n})S_{n}u_{n}+Q_{n},
\end{equation*}
where $(\xi,\eta)$ are the change of variables given in section
$\S$2 by
\begin{equation*}
a_{12}(w_{n})\xi_{x}+a_{22}(w_{n})\xi_{y}=0\text{ }\text{ }\text{
in }\text{ }X, \text{ }\text{ } \xi(x,0)=x,\text{ }\text{ }\xi(\pm
x_{0},y)=\pm x_{0},\text{ }\text{ }\eta=y,
\end{equation*}
and where $Q_{n}$ is the quadratic error term given by
\begin{equation*}
Q_{n}=\int_{0}^{1}(t-1)\partial_{t}^{2}\Phi(w_{n}+tS_{n}u_{n})dt.
\end{equation*}
Since $L_{\theta_{n}}(w_{n})u_{n}=f_{n}$ in $X_{n+1}$, we have
\begin{equation}
f_{n+1}=L_{\theta_{n}}(w_{n})(u_{n}-S_{n}u_{n})-\theta_{n}(S_{n}u_{n})_{
\eta\eta\xi\xi}+Q_{n},
\end{equation}
in $X_{n+1}$.\par
   Each term of (5.4) shall be estimated separately.  First note
that $\theta_{n}$ may be chosen sufficiently small to guarantee
that,
\begin{equation*}
\parallel\theta_{n}(S_{n}u_{n})_{\eta\eta\xi\xi}\parallel_{s}^{n+1}
\leq\frac{1}{3}C_{2}\mu_{n+1}^{\tau^{-1}(s-s_{*}+18+2\delta)}
\parallel f_{1}\parallel_{s_{*}-15}^{1}.
\end{equation*}\par
   We now estimate $L_{\theta_{n}}(w_{n})(u_{n}-S_{n}u_{n})$.  By
lemma 4.2 and IV$_{n}$,
\begin{eqnarray*}
\parallel
L_{\theta_{n}}(w_{n})(u_{n}-S_{n}u_{n})\parallel_{s}^{n+1}\!\!\!&\leq&\!\!\!
\parallel
L_{\theta_{n}}(w_{n})(u_{n}-S_{n}u_{n})\parallel_{s}^{n}\\
&\leq&\!\!\! M_{1}(\parallel u_{n}-S_{n}u_{n}\parallel_{s+2}^{n}+
\parallel w_{n}\parallel_{s+4}^{n}|u_{n}-S_{n}u_{n}|_{0}^{n})\\
& &\!\!\!+O(\theta_{n})\\
&\leq&\!\!\! M_{2}(\parallel u_{n}-S_{n}u_{n}\parallel_{s+2}^{n}+
\parallel w_{n}\parallel_{s+4}^{n}\parallel u_{n}-S_{n}u_{n}\parallel_{2}^{n})\\
& &\!\!\!+O(\theta_{n})\\
&\leq&\!\!\! M_{3}(\mu_{n}^{s+17-s_{*}}\parallel
u_{n}\parallel_{s_{*}-15}^{n}+\mu_{n}^{17-s_{*}}\parallel
w_{n}\parallel_{s+4}^{n}\parallel u_{n}\parallel_{s_{*}-15}^{n})
\\& &\!\!\!+O(\theta_{n}).
\end{eqnarray*}
Furthermore, by (5.3)
\begin{equation*}
\parallel u_{n}\parallel_{s_{*}-15}^{n}\leq
M_{4}\mu_{n}^{2s_{*}\sigma+\delta}\parallel
f_{1}\parallel_{s_{*}-15}^{1}.
\end{equation*}
If $\theta_{n}$ and $\sigma$ are sufficiently small and $\mu$ is
sufficiently large, it follows that
\begin{eqnarray*}
\parallel
L_{\theta_{n}}(w_{n})(u_{n}-S_{n}u_{n})\parallel_{s}^{n+1}&\leq&
M_{5}\mu^{-1}(\mu_{n}^{3s_{*}\sigma+s-s_{*}+17+\delta}+
\mu_{n}^{3s_{*}\sigma-s_{*}+17+2\delta})\parallel
f_{1}\parallel_{s_{*}-15}^{1}\\
&\leq& \frac{1}{3}C_{2}\mu_{n+1}^{\tau^{-1}(s-s_{*}+18+2\delta)}
\parallel f_{1}\parallel_{s_{*}-15}^{1}.
\end{eqnarray*}\par
   We now estimate $Q_{n}$.  Apply lemma 4.2 ($ii$) to obtain,
\begin{eqnarray*}
\parallel Q_{n}\parallel_{s}^{n+1}&\leq&\parallel
Q_{n}\parallel_{s}^{n}\\
&\leq&\int_{0}^{1}\sum_{|\alpha|,|\beta|,|\rho|\leq 2}\parallel
\partial^{\rho}\Phi(w_{n}+tS_{n}u_{n})\partial^{\alpha}(S_{n}u_{n})
\partial^{\beta}(S_{n}u_{n})\parallel_{s}^{n}dt\\
&\leq&\int_{0}^{1}\sum_{|\alpha|,|\beta|,|\rho|\leq 2}M_{6}(|
\partial^{\rho}\Phi(w_{n}+tS_{n}u_{n})|_{0}^{n}\parallel\partial^{\alpha}(S_{n}u_{n})
\partial^{\beta}(S_{n}u_{n})\parallel_{s}^{n}\\
& &\text{ }\text{ }\text{ }\text{ }\text{ }\text{ }\text{ }\text{
}\text{ }\text{ }\text{ }\text{ }\text{ }\text{ }\text{ }\text{
}\text{ }\text{ }+\parallel
\partial^{\gamma}\Phi(w_{n}+tS_{n}u_{n})\parallel_{s}^{n}|\partial^{\alpha}(S_{n}u_{n})
\partial^{\beta}(S_{n}u_{n})|_{0}^{n})dt.
\end{eqnarray*}
Then the Sobolev lemma and the interpolation inequality $\parallel
u^{2}\parallel_{L^{2}}\leq C\parallel u\parallel_{H^{1}}^{2}$,
show that
\begin{eqnarray*}
\parallel Q_{n}\parallel_{s}^{n+1}&\leq&
\int_{0}^{1}\sum_{|\rho|\leq 2}M_{7}(\parallel
\partial^{\rho}\Phi(w_{n}+tS_{n}u_{n})\parallel_{2}^{n}
(\parallel S_{n}u_{n}\parallel_{s+3}^{n})^{2}\\
& &\text{ }\text{ }\text{ }\text{ }\text{ }\text{ }\text{ }\text{
}\text{ }\text{ }\text{ }\text{ }\text{ }+\parallel
\partial^{\rho}\Phi(w_{n}+tS_{n}u_{n})\parallel_{s}^{n}
(\parallel S_{n}u_{n}\parallel_{4}^{n})^{2})dt.
\end{eqnarray*}
Furthermore, by lemma 4.2 ($iii$), I$_{n}$, IV$_{n}$, and
proposition 5.2,
\begin{equation*}
\parallel Q_{n}\!\parallel_{s}^{n+1}\leq M_{8}
[(\parallel w_{n}\!\parallel_{6}^{n}\!+\mu_{n}^{2}\parallel
u_{n}\!\parallel_{4}^{n})(\mu_{n}^{3}\parallel
u_{n}\!\parallel_{s}^{n})^{2}\!+\!(\parallel
w_{n}\!\parallel_{s+4}^{n}\!+\mu_{n}^{4}\parallel
u_{n}\!\parallel_{s}^{n})(\parallel u_{n}\!\parallel_{4}^{n})^{2}]
\end{equation*}\vspace{-6mm}
\begin{equation*}
\leq (M_{9}\parallel f_{1}\parallel_{s_{*}-15}^{1})
[(1+\mu_{n}^{2+\tau^{-1}(-s_{*}+22+
2\delta)})\mu_{n}^{6+2\tau^{-1}(s-s_{*}+18+2\delta)}\text{ }\text{
}\text{ }\text{ }\text{ }\text{ }\text{ }\text{ }\text{ }
\end{equation*}
\begin{equation*}
\text{ }
+(\mu_{n}^{\sigma(s-11)+\delta}+\mu_{n}^{4+\tau^{-1}(s-s_{*}+18
+2\delta)})\mu_{n}^{2\tau^{-1}(-s_{*}+22+2\delta)}]
\parallel f_{1}\parallel_{s_{*}-15}^{1}
\end{equation*}
\begin{equation*}
\leq (M_{10}\parallel
f_{1}\parallel_{s_{*}-15}^{1})\mu_{n}^{s-s_{*}+18+2\delta}\parallel
f_{1}\parallel_{s_{*}-15}^{1},\text{ }\text{ }\text{ }\text{
}\text{ }\text{ }\text{ }\text{ }\text{ }\text{ }\text{ }\text{
}\text{ }\text{ }\text{ }\text{ }\text{ }\text{ }\text{ }\text{
}\text{ }\text{ }\text{ }\text{ }\text{ }\text{ }\text{ }\text{
}\text{ }\text{ }
\end{equation*}
since $s\leq s_{*}-18-2\delta-\frac{6\tau}{2-\tau}$ and $s_{*}\geq
22+2\delta+\frac{6\tau}{2-\tau}$.  If $\varepsilon(s_{*})$ is
sufficiently small to guarantee that $M_{10}\parallel
f_{1}\parallel_{s_{*}-15}^{1}\leq\frac{1}{3}C_{2}$, then
\begin{equation*}
\parallel Q_{n}\parallel_{s}^{n+1}\leq
\frac{1}{3}C_{2}\mu_{n+1}^{\tau^{-1}(s-s_{*}+18+2\delta)}
\parallel f_{1}\parallel_{s_{*}-15}^{1}.
\end{equation*}
By combining the estimates for each term of (5.4) we obtain the
desired result.\qed

\noindent\textbf{Proposition 5.4.}  \textit{If $n_{0}(s_{*})$ is
sufficiently large, then}
\begin{equation*}
\parallel w_{n+1}\parallel_{14}^{n+1}\leq C_{3},
\end{equation*}
\textit{where $C_{3}$ depends on $\mu$ and $s_{*}$}.\medskip

\noindent\textit{Proof.}  Let $a=14+\tau^{-1}(18+2\delta-s_{*})$,
and note that since $s_{*}\geq 22+2\delta+\frac{6\tau}{2-\tau}$,
$\tau\geq\frac{3}{2}$, we have $a<0$.  If $n_{0}$ is sufficiently
large, we may apply proposition 5.2 to obtain,
\begin{eqnarray*}
\text{ }\text{ }\text{ }\text{ }\text{ }\text{ }\text{ }\text{
}\text{ }\text{ }\text{ }\text{ }\text{ }\text{ }\text{ }\text{
}\text{ }\text{ }\text{ }\text{ }\text{ }\text{ }\text{ }
\parallel w_{n+1}\parallel_{14}^{n+1}&\leq&
\sum_{i=1}^{n}\parallel S_{i}u_{i}\parallel_{14}^{i}\\
&\leq&\sum_{i=1}^{n}\mu_{i}^{14}\parallel u_{i}\parallel_{0}^{i}\\
&\leq&\sum_{i=1}^{\infty}\mu_{i}^{a}\parallel
f_{1}\parallel_{s_{*}-15}^{1}:=C_{3}.\text{ }\text{ }\text{
}\text{ }\text{ }\text{ }\text{ }\text{ }\text{ }\text{ }\text{
}\text{ }\text{ }\text{ }\text{ }\text{ }\text{ }\text{ }\text{
}\text{ }\text{ }\text{ }\text{ }\text{ }\text{ }\text{ }\text{
}\text{ }\text{ }\text{ }\Box
\end{eqnarray*}\par
   In order to obtain the largest value for $s$ and smallest
lower bound for $s_{*}$ which satisfy the conditions of
propositions 5.1 - 5.4, we choose $\tau=1.6$ so that $s_{*}\geq
100$ and $s\leq s_{*}-96$.  We now establish two corollaries which
will complete the proof of theorem 1.3.\medskip

\noindent\textbf{Corollary 5.1.}  $w_{n}\rightarrow w$\textit{ in
}$H^{s_{*}-96}(X_{\infty})$.\medskip

\noindent\textit{Proof.}  If $s\leq s_{*}-96$, then by II$_{n}$
\begin{eqnarray*}
\parallel w_{i}-w_{j}\parallel_{s}^{\infty}&\leq&
\sum_{k=j}^{i}\parallel u_{k}\parallel_{s}^{k}\\
&\leq&
C_{1}\sum_{k=j}^{i}\mu_{k}^{\tau^{-1}(s-s_{*}+18+2\delta)}\parallel
f_{1}\parallel_{s_{*}-15}^{1}.
\end{eqnarray*}
Hence, $\{w_{n}\}$ is Cauchy in $H^{s}(X_{\infty})$ for all $s\leq
s_{*}-96$ since $18+2\delta<96$.\qed

\noindent\textbf{Corollary 5.2.}  $\Phi(w_{n})\rightarrow
0$\textit{ in }$H^{s_{*}-96}(X_{\infty})$.\medskip

\noindent\textit{Proof.}  If $s\leq s_{*}-96$, then by III$_{n}$
\begin{eqnarray*}
\text{ }\text{ }\text{ }\text{ }\text{ }\text{ }\text{ }\text{
}\text{ }\text{ }\text{ }\text{ }\text{ }\text{ }\text{ }\text{
}\text{ }\text{ }\text{ }\text{ }\text{ }\text{ }\text{ }
\parallel \Phi(w_{n})\parallel_{s}^{\infty}&\leq&\parallel f_{n}
\parallel_{s}^{n}\\
&\leq& C_{2}\mu_{n}^{\tau^{-1}(s-s_{*}+18+2\delta)}\parallel
f_{1}\parallel_{s_{*}-15}^{1}\rightarrow 0.\text{ }\text{ }\text{
}\text{ }\text{ }\text{ }\text{ }\text{ }\text{ }\text{ }\text{
}\text{ }\text{ }\text{ }\text{ }\text{ }\text{ }\text{ }\Box
\end{eqnarray*}\par
   Since $s_{*}\geq 100$, it follows that $w_{n}\rightarrow w$ in
$C^{2}(\overline{X}_{\infty})$.  Therefore
$\Phi(w_{n})\rightarrow\Phi(w)$, showing that $w$ is a solution of
(5.1).  Furthermore, if $l$ is as in theorem 1.3, then we have
$w\in C^{l-98}$, $l\geq 100$.  This completes the proof of theorem
1.3.

\bigskip\bigskip
\noindent\textbf{References}
\bigskip\bigskip

\noindent1.\hspace{.08in} Birkhoff, G., Rota, G.-C.: Ordinary
Differential Equations.  Blaisdell Publishing,
\par\hspace{.04in}London, 1969.

\noindent2.\hspace{.08in} Friedrichs, K. O.: The identity of weak
and strong extensions of differential oper-
\par\hspace{.04in}ators.  Trans. Amer. Math. Soc. \textbf{55},
132-151 (1944).

\noindent3.\hspace{.08in} Gallerstedt, S.: Quelques probl\`{e}mes
mixtes pour l'\'{e}quation
$y^{m}z_{xx}+z_{yy}=0$.\par\hspace{.04in}Arkiv f\"{o}r Matematik,
Astronomi och Fysik \textbf{26A} (3), 1-32 (1937).

\noindent4.\hspace{.08in} Han, Q.:  On the isometric embedding of
surfaces with Gauss curvature changing\par\hspace{.04in}sign
cleanly.  Comm. Pure Appl. Math. \textbf{58}, 285-295 (2005).

\noindent5.\hspace{.08in} Han, Q.: Local isometric embedding of
surfaces with Gauss curvature changing\par\hspace{.04in}sign
stably across a curve. Cal. Var. \& P.D.E. \textbf{25}, 79-103
(2006).

\noindent6.\hspace{.08in} Han, Q.: Smooth local isometric
embedding of surfaces with Gauss
curvature\par\hspace{.04in}changing sign cleanly. Preprint.

\noindent7.\hspace{.08in} Han, Q., Hong, J.-X.:  Isometric
Embedding of Riemannian Manifolds in Eu-\par\hspace{.04in}clidean
Spaces. Mathematical Surveys and Monographs, Vol. 130, AMS,
Prov-\par\hspace{.04in}idence, RI, 2006.

\noindent8.\hspace{.08in} Han, Q., Hong, J.-X., Lin, C.-S.: Local
isometric embedding of surfaces with
non-\par\hspace{.04in}positive curvature. J. Differential Geom.
\textbf{63}, 475-520 (2003).

\noindent9.\hspace{.08in} Han, Q., Khuri, M.:  On the local
isometric embedding in $\mathbb{R}^{3}$ of surfaces
with\par\hspace{.04in}Gaussian curvature of mixed sign. Preprint.

\noindent10. Jacobowitz, H.: Local isometric embeddings. Seminar
on Differential Geometry,\par\hspace{.04in}edited by S.-T. Yau,
Annals of Math. Studies \textbf{102}, 1982, 381-393.

\noindent11. Khuri, M.: The local isometric embedding in
$\mathbb{R}^{3}$ of two-dimensional
Riemannian\par\hspace{.04in}manifolds with Gaussian curvature
changing sign to finite order on a curve.\par\hspace{.04in}J.
Differential Geom., to appear.

\noindent12. Khuri, M.: Counterexamples to the local solvability
of Monge-Amp\`{e}re equations\par\hspace{.04in}in the plane. Comm.
PDE \textbf{32}, 665-674 (2007).

\noindent13. Ladyzenskaja, O. A., Solonnikov, V. A., Ural'ceva, N.
N.: Linear and Quasi-\par\hspace{.04in}Linear Equations of
Parabolic Type. Translations of Mathematical
Mono-\par\hspace{.04in}graphs \textbf{23}, 1968.

\noindent14. Lax, P. D., Phillips, R. S.: Local boundary
conditions for dissipative symmetric\par\hspace{.04in}linear
differential operators. Comm. Pure Appl. Math. \textbf{13},
427-455 (1960).

\noindent15. Lin, C.-S.: The local isometric embedding in
$\mathbb{R}^{3}$ of 2-dimensional
Riemannian\par\hspace{.04in}manifolds with nonnegative curvature.
J. Differential Geom. \textbf{21}, 213-230 (1985).

\noindent16. Lin, C.-S.: The local isometric embedding in
$\mathbb{R}^{3}$ of two-dimensional
Riemannian\par\hspace{.05in}manifolds with Gaussian curvature
changing sign cleanly. Comm. Pure Appl.\par\hspace{.05in}Math.
\textbf{39}, 867-887 (1986).

\noindent17. Nadirashvili, N., Yuan, Y.: Improving Pogorelov's
isometric embedding coun-\par\hspace{.05in}terexample. Preprint.

\noindent18. Pogorelov, A. V.: An example of a two-dimensional
Riemannian metric not\par\hspace{.05in}admitting a local
realization in $E_{3}$. Dokl. Akad. Nauk. USSR \textbf{198}, 42-43
(1971).

\noindent19. Poznyak, E. G.: Regular realization in the large of
two-dimensional metrics of\par\hspace{.05in}negative curvature.
Soviet Math. Dokl. \textbf{7}, 1288-1291 (1966).

\noindent20. Poznyak, E. G.: Isometric immersions of
two-dimensional Riemannian metrics in\par\hspace{.03in}Euclidean
space. Russian Math. Surveys \textbf{28}, 47-77 (1973).

\noindent21. Peyser, G.: On the identity of weak and strong
solutions of differential equations\par\hspace{.04in}with local
boundary conditions. Amer. J. Math. \textbf{87}, 267-277 (1965).

\noindent22. Schwartz, J. T.: Nonlinear Functional Analysis. New
York University, New York,\par\hspace{.02in}1964.

\noindent23. Stein, E.: Singular Integrals and Differentiability
Properties of Functions. Prince-\par\hspace{.04in}ton University
Press, Princeton, 1970.

\noindent24. Taylor, M. E.: Partial Differential Equations III.
Springer-Verlag, New York,\par\hspace{.04in}1996.

\noindent25. Weingarten, J.: \"{U}ber die theorie der Aubeinander
abwickelbarren Oberfl\"{a}chen.\par\hspace{.04in}Berlin, 1884.

\end{document}